\documentclass[12pt]{amsart}

\usepackage{mathrsfs}

\usepackage{caption}

\usepackage{amsthm}
\usepackage{hyperref}
\usepackage{amssymb}
\usepackage{latexsym}
\usepackage{amsmath}
\usepackage{amsfonts}
\usepackage{mathabx}
\usepackage[dvips]{graphicx}
\usepackage[utf8]{inputenc}

\usepackage[LGR,T1]{fontenc}%
\usepackage[french,english]{babel}
\usepackage{url}
\usepackage{enumerate}
\usepackage{hyperref}
\usepackage{color}

\usepackage{comment}
\usepackage{cancel}

\usepackage{fancyhdr}
\usepackage{bbm}
\usepackage{tikz}

\newtheorem{Theorem}{Theorem}[part]
\newtheorem{Definition}{Definition}[part]
\newtheorem{Proposition}{Proposition}[part]

\newtheorem{Assumption}{Assumption}[part]
\newtheorem{Lemma}{Lemma}[part]
\newtheorem{Corollary}{Corollary}[part]
\newtheorem{Remark}{Remark}[part]

\makeatletter \@addtoreset{equation}{section}

\@addtoreset{Definition}{section}

\@addtoreset{Theorem}{section}

\@addtoreset{Proposition}{section}

\@addtoreset{Property}{section}

\@addtoreset{Assumption}{section}

\@addtoreset{Corollary}{section}

\@addtoreset{Lemma}{section}

\@addtoreset{Remark}{section}

\@addtoreset{Example}{section}

\newcommand{\cD}{\mathcal{D}}
\newcommand{\cE}{\mathcal{E}}
\newcommand{\cF}{\mathcal{F}}

\newcommand{\cL}{\mathcal{L}}
\newcommand{\cM}{\mathcal{M}}

\newcommand{\cS}{\mathcal{S}}

\newcommand{\E}{\mathbb{E}}

\renewcommand{\P}{\mathbb{P}}

\newcommand{\R}{\mathbb{R}}

\def \proof{{\noindent \bf Proof. }}
\def \eproof{\hbox{ }\hfill$\Box$}

\newcommand{\ud}{\mathrm{d}}

\newcommand{\1}{{\bf 1}}
    
\newcommand{\set}[1]
    {\ensuremath{\{ #1 \}}}
    
\newcommand{\HP}[1] 
    {\ensuremath{\mathscr{H}^{#1}}}

\newcommand{\SP}[1]{\ensuremath{\mathscr{S}^{#1}}}

\newcommand{\esp}[1]{\ensuremath{\mathbb{E} \!\! \left[#1\right] }}

\renewcommand{\Xi}[1]{X_{i #1}}

\newcommand{\BP}[1]{\ensuremath{\mathscr{B}^{#1}}}
\newcommand{\NB}[2]{\ensuremath{\left\| #2 \right\|_{\mathscr{B}^{#1}}} }
\newcommand{\NL}[2]{\ensuremath{\left\| #2 \right\|_{\mathscr{L}^{#1}}} }

\newcommand{\bone}{\mathbf{1}}
\newcommand{\EE}{\mathbb{E}}

\definecolor{jcg}{HTML}{0aa344}

\begin{document}
\title[Martingales on a manifold with a boundary and reflected BSDEs]{Martingales on a Euclidean manifold with a boundary and reflected BSDEs in non-convex domains}
\author[M. ARNAUDON]{Marc ARNAUDON}
\address{Université de Bordeaux, IMB, UMR 5251, F-33400 Talence, France.}
\email{marc.arnaudon@math.u-bordeaux.fr}
\author[J.-F. CHASSAGNEUX]{Jean-François CHASSAGNEUX}
\address{ENSAE-CREST \& Institut Polytechnique de Paris.}
\email{chassagneux@ensae.fr}
\author[S. NADTOCHIY]{Sergey NADTOCHIY}
\address{Department of Mathematical Sciences, Carnegie Mellon University, Pittsburgh, PA 15213.}
\email{{snadtoch@andrew.cmu.edu}}
\author[A. RICHOU]{Adrien RICHOU}
\address{Université de Bordeaux, IMB, UMR 5251, F-33400 Talence, France.}
\email{adrien.richou@math.u-bordeaux.fr}
\thanks{Partial support from the ANR Project ReLISCoP (ANR-21-CE40-0001)is acknowledged.}

\begin{abstract}
{The purpose of this paper is twofold. First, we introduce the notion of a $\Gamma$-martingale on a Euclidean manifold with a boundary (i.e., the closure of an open connected domain in $\R^d$), we provide its equivalent characterization through the $\Gamma$-convex functions, and we establish its connection with the reflected backward stochastic differential equations (BSDEs) in the associated domain. Second, we show how the tools of stochastic geometry can be used to develop a new method for proving existence and uniqueness of solutions to reflected BSDEs. We implement this method and obtain a well-posedness result for reflected BSDEs in any bounded, two-dimensional, simply-connected domain that is locally $C^2$-diffeomorphic to a convex set. This work extends the results of \cite{Chassagneux-Nadtochiy-Richou-22} and \cite{Gegout-Petit-Pardoux-96}.}
\end{abstract}

\maketitle




Mathematics Subject Classification: 60D05, 60G65, 60J60

\section{Introduction}
\subsection{Motivation and main contributions}
{
Backward stochastic differential equations (BSDEs), first introduced in the linear case by \cite{Bismut-76} and later in a general framework by \cite{Pardoux-Peng-90}, have since been extensively studied due to their wide range of applications, particularly in stochastic control, financial mathematics, and their connections with semi-linear partial differential equations (PDEs); see, e.g., \cite{Zhang-17} for an overview.
BSDEs can be viewed as a nonlinear generalization of conditional expectations in $\mathbb{R}^d$ with respect to a given filtration. In this context, the resulting notion is referred to as a nonlinear expectation or $g$-expectation (see \cite{Peng-97}). Notably, a solution to a BSDE with a zero generator recovers the classical conditional expectation of the terminal condition.
Now, suppose that a terminal condition takes values in the closure $\bar{\mathcal{D}}$ of a domain $\mathcal{D} \subset \mathbb{R}^d$. It is then natural to ask whether one can find a natural extension of the notion of conditional expectation that remains in $\bar{\mathcal{D}}$. This leads to the concept of a reflected BSDE (with a zero generator).

The theory of reflected BSDEs is well understood in spatial dimension one (i.e., when the $Y$-component of the solution evolves in an interval and is reflected at its boundary); see, for example, \cite{el1997reflected, cvitanic1996backward, dumitrescu2016generalized, Grigorova-Imkeller-Offen-Ouknine-Quenez-17, Grigorova-Imkeller-Ouknine-Quenez-18}. However, the multidimensional case poses serious additional challenges, notably due to the lack of a comparison principle. Most of the existing well-posedness results in higher dimensions have been established under the convexity assumption on the domain; see, e.g., \cite{Gegout-Petit-Pardoux-96, Klimsiak-Rozkosz-Slominski-15, Chassagneux-Richou-17, Fakhouri-Ouknine-Ren-18}. To our knowledge, the only result in a non-convex setting is found in \cite{Chassagneux-Nadtochiy-Richou-22}.
To see the importance of convexity, notice that, when the generator is zero and $\mathcal{D}$ is convex, the $Y$-component of a solution coincides with the conditional expectation, which automatically remains in $\bar{\mathcal{D}}$, hence no reflection is needed. The latter is not the case in a non-convex case. In other words, for convex domains, the reflection only needs to counter the drift term arising from the generator, whereas, in a non-convex case, it may also need to take into account the martingale term. This observation explains the challenge of extending the analysis from a convex to a non-convex domain. It also illustrates that, contrary to most works on BSDEs, the zero-generator case is already non-trivial and mathematically rich: it corresponds to choosing a notion of conditional expectation that is constrained to stay in $\bar{\mathcal{D}}$.
\newline

The work \cite{Chassagneux-Nadtochiy-Richou-22} proves several existence and uniqueness results for reflected BSDEs in a fairly restrictive setting -- assuming the weak star-shape property and excessive smoothness of the domain $\mathcal{D}$, as well as a smallness property of the terminal data (though the latter is not needed in a Markovian framework) -- which notably does not fully cover the setting of \cite{Gegout-Petit-Pardoux-96}.
A particularly insightful remark in \cite{Chassagneux-Nadtochiy-Richou-22} connects the zero-generator case with the theory of $\Gamma$-martingales on manifolds. Specifically, when the terminal condition lies in a sufficiently ``concave'' part of $\partial \mathcal{D}$ (as viewed from inside the domain), the solution remains on the boundary and becomes a $\Gamma$-martingale on the boundary manifold (see Proposition 5.1 in \cite{Chassagneux-Nadtochiy-Richou-22}). On the other hand, if the terminal condition lies within a convex subset of $\mathcal{D}$, the solution is a classical martingale in $\mathbb{R}^d$.
These observations naturally lead to the following questions:
\begin{enumerate}
\item Can one define a natural notion of a $\Gamma$-martingale on a manifold with a boundary (here, $\bar{\mathcal{D}}$ is viewed as a Euclidean manifold with a boundary)?
\item Is a solution of a reflected BSDE with zero generator a $\Gamma$-martingale?
\end{enumerate}
The first goal of this paper is to provide positive answers to both questions.
\newline

To the best of our knowledge, until now, the notion of a $\Gamma$-martingale has only been introduced for manifolds without boundaries; see, e.g., \cite{Meyer-81,Darling-82,Emery-89}.
There are several ways to define (or characterize) the notion of a $\Gamma$-martingale in a (smooth) manifold without boundary (imbedded in $\R^d$). One possibility is to define it as a process that lives on the manifold and is a sum of a (usual) Euclidean martingale and a finite-variation process whose velocity is orthogonal to the manifold. Alternatively, one can characterize a $\Gamma$-martingale as a process $X$ such that, for any $\Gamma$-convex function $\psi$ (i.e. a function $\R^d \to \R$ that is convex along the geodesics of the manifold), $\psi(X)$ is (locally) a Euclidean sub-martingale.
In this article, we follow the same approach to define a $\Gamma$-martingale on a manifold with a boundary, but we intentionally avoid working in the local coordinates of the manifold and present all statements and derivations using the coordinates of $\R^d$, into which our manifold ($\bar{\cD} \subset \R^d$) is embedded, as this is natural for the connection we make with the reflected BSDEs, and because it makes the paper accessible for a reader without a background in differential geometry.
The proposed notion of a $\Gamma$-martingale is given in Definition \ref{def:gamma-mart}, while its characterization through the $\Gamma$-convex functions is established in Propositions \ref{prop gamma-function} and \ref{prop caracterisation gamma martingale}.
\newline

The second goal of this paper is to leverage the tools of stochastic geometry in order to improve the existence and uniqueness results of \cite{Chassagneux-Nadtochiy-Richou-22}.
It is mentioned in the latter paper (see the discussion in the second half of Section 5 in \cite{Chassagneux-Nadtochiy-Richou-22}) that, for $d\geq 3$, the uniqueness of a solution to a reflected BSDE may fail in general -- i.e., with a general non-Markovian terminal data and a general domain $\cD$ -- even if the domain is infinitely smooth and the generator is equal to zero.
In view of this observation, it appears natural to restrict our analysis to $d=2$, if the goal is to find a qualitative, as opposed to quantitative, condition on the domain that would guarantee the well-posedness of the associated reflected BSDEs and would include non-convex domains. 
However, even for $d=2$, the uniqueness of solutions to the reflected BSDEs associated with a given domain is expected to fail if the domain $\cD$ is not simply connected. 
Indeed, consider a domain $\cD$ with a circular hole inside and any two opposite points on this circle. The two half-circles connecting these points form two different geodesics in $\bar\cD$. Then, choosing a terminal condition that is supported on these two points, one can easily construct two $\Gamma$-martingales which evolve on the aforementioned geodesics and which both have the prescribed terminal value.
Thus, a general well-posedness result for $d=2$ can only be expected for simply-connected domains $\cD$.
In Theorems \ref{thm existence bsde} and \ref{Thm-Uniqu1} herein, we prove existence and uniqueness of solutions to general reflected BSDEs assuming that $\bar{\cD}$ is two-dimensional, bounded, simply connected and locally $C^2$-diffeomorphic to a convex set. In particular, for $d=2$, the latter theorems generalize the results of \cite{Chassagneux-Nadtochiy-Richou-22} and \cite{Gegout-Petit-Pardoux-96}.

The approach employed herein for the proof of existence differs substantially from that of \cite{Chassagneux-Nadtochiy-Richou-22}. Indeed, \cite{Chassagneux-Nadtochiy-Richou-22} follows the algorithm that is standard for this type of problems, by considering a penalized version of the associated BSDE without reflection and making the penalty term tend to $+\infty$. In the present article, on the contrary, we adapt the intrinsic method used by Kendall to construct martingales on manifolds in \cite{Kendall-90}. Namely, we notice that any bounded two-dimensional simply-connected $\bar\cD$, equipped with its geodesic distance, is a CAT(0) space (a.k.a. Hadamard space). We then consider the classical notion of a mean in metric spaces, known as a Fréchet mean, and apply a backward recursion to construct a dynamic version of this mean. However, in order to cover the reflected BSDEs with non-zero exogenous generators (i.e., any generator defined as a given stochastic process, as opposed to a feedback function of the solution), we have to modify the latter construction: between any two steps of the aforementioned recursion, we add a transport step in the direction prescribed by the generator. This construction is implemented in the proofs of Proposition \ref{prop:Gamma.mtg.drift.Markov} and Theorem \ref{thm existence bsde exogeneous}. The general case is, then, treated by using a Picard iteration scheme and the tools of BMO martingales. 

The proofs of both existence and uniqueness results (Proposition \ref{prop:Gamma.mtg.drift.Markov} and Theorems \ref{Thm-Uniqu1}, \ref{Thm-Uniqu2}, \ref{thm existence bsde exogeneous}) rely crucially on the fact that the squared geodesic distance, viewed as a function on $\bar\cD\times\bar\cD$, is convex along any two geodesics. In particular, the proof of uniqueness is based on using the squared geodesic distance as a Lyapunov's function for the associated reflected BSDE. This convexity property of the squared geodesic distance is enjoyed by all CAT(0) spaces, which explains why the connection to CAT(0) is so important and provides another explanation of why the restriction to two-dimensional simply-connected domains $\cD$ is natural. Note that simply connected domains in higher dimensions are not, in general, CAT(0) spaces, and neither are the two-dimensional domains that are not simply-connected. 


We conjecture that our main results (Theorems \ref{Thm-Uniqu1}, \ref{Thm-Uniqu2}, \ref{thm existence bsde exogeneous} and \ref{thm existence bsde}) remain valid in the setting where $\cD$ is a $2$-dimensional Cartan-Hadamard manifold (i.e., simply connected with a non-positive sectional curvature), satisfying the regularity condition of Assumption \ref{ass:main}. The reason is that Cartan-Hadamard manifolds are CAT(0) spaces, implying that the square of their geodesic distance function is convex and smooth. In addition, it is well known that, in these manifolds, the Fréchet means of compactly supported probability measures exist, are unique and depend smoothly on the measures. To establish such an extension, in all the proofs herein, one would need to replace the Euclidean lines with geodesics and to perform linearizations via parallel translations along geodesics. 
\newline

The remainder of this paper is organized as follows. Section 1.2 sets the notations, and Section 1.3 states the main assumption (Assumption \ref{ass:main}) which holds throughout the paper, as well as several corollaries of this assumption. Section 2 introduces the definition of $\Gamma$-martingales (with drifts) in a Euclidean manifold with a boundary (Definitions \ref{def:gamma-mart-0} and \ref{def:gamma-mart}), establishes their equivalent characterization through the $\Gamma$-convex functions (Proposition \ref{prop gamma-function} and Proposition \ref{prop caracterisation gamma martingale}), and describes their connection to the reflected BSDEs (Remark \ref{rem link bsde}). Section 3 establishes further properties of $\bar\cD$ and of its geodesic distance, under the additional assumption that $\cD$ is two-dimensional and simply connected. In particular, an Itô's formula for the squared geodesic distance is given in Corollary \ref{Ito-formula-Psi}, and a uniqueness and stability result for $\Gamma$-martingales with prescribed drifts and terminal conditions is stated in Theorem \ref{Thm-Uniqu1} for a general continuous filtration. Theorem \ref{thm existence bsde} in Section 4 states the desired existence and uniqueness result for solutions to the reflected BSDEs, assuming a Brownian filtration and under the same assumptions on $\cD$ as in Section 3 (this implies a corresponding existence and uniqueness result for the $\Gamma$-martingales with drifts, in a Brownian filtration, via Remark \ref{rem link bsde}). {Finally, in Section 5, we state and prove several auxiliary results related to Fr\'echet mean.}

}

\subsection{Notations}
\smallskip
We consider a complete probability space $(\Omega, \mathcal{F},\mathbb{P})$ equipped with a continuous filtration $(\mathcal{F}_t)_{t \geqslant 0}$. Some results are obtained in a Brownian setting: in this case, $(\mathcal{F}_t)_{t \geqslant 0}$ will denote the augmented natural filtration of a Brownian motion $(W_t)_{t \geqslant 0}$ in $\mathbb{R}^{d'}$. We set a terminal time $T>0$.

For $p\geq1$, we denote by $\mathcal{L}^p$ the space of (classes of equivalence of)\footnote{We drop this clarification in further definitions.} $\mathcal{F}_T$-measurable random variables $\xi$ (with values in a Euclidean space), s.t. $\|\xi\|_{\mathcal{L}^p}:=\esp{|\xi|^p}^{1/p}<\infty$. The space $\mathcal{L}^\infty$ stands for all $\mathcal{F}_T$-measurable essentially bounded random variables.
We define $\SP{p}$ as the space of continuous adapted process (with values in a Euclidean space) $Y$, s.t. $\|Y\|_{\SP{p}}:=\| \sup_{t \in [0,T]} |Y_t| \|_{\mathcal{L}^p} <\infty$.
We define $\SP{\infty}$ as the space of continuous adapted processes (with values in a Euclidean space) $Y$, s.t. $\|Y\|_{\SP{\infty}}:=\| \sup_{t \in [0,T]} |Y_t| \|_{\mathcal{L}^\infty} <\infty$.
We also define $\HP{p}$ as the space of progressively measurable processes (with values in a Euclidean space) $Z$, s.t. $\|Z\|_{\HP{p}}:=\esp{ \int_0^T |Z_t|^p \ud t }^{1/p}<\infty$, while $\HP{\infty}$ is the space of progressively measurable processes $Z$, s.t. $\|Z\|_{\HP{\infty}}:=\| \int_0^T |Z_t|^2 \ud t \|_{\mathcal{L}^\infty}^{1/2}<\infty$.
Next, for $p\geq1$, we define $\mathcal{M}^p$ as the space of all continuous local martingales $M$ with $\|M\|_{\mathcal{M}^p}:=\esp{\langle M\rangle_T^{p/2}}^{1/p}<\infty$.
We also denote by $\mathrm{Var}_{t}(K)$ the variation of a process $K_\cdot$ (with values in a Euclidean space) on the time interval $[0,t]$ and by $\mathscr{K}^p$, for $p \in [1,\infty]$, the set of finite-variation process $K$ such that $\NL{p}{\mathrm{Var}_{[0,T]}(K)} < \infty$ and $K_0=0$.
Finally, we denote by $\BP{2}$ the set of processes $V \in \HP{2}$, satisfying
\begin{align*}
 \NB{2}{V} := \NL{\infty}{\mathrm{sup}_{t \in [0,T]} \esp{\int_{t}^T |V_s|^2 \ud s | \cF_{t}}}^\frac12<+\infty.  \;
\end{align*}
Let us remark that $V \in \BP{2}$ implies that the martingale $\int_0^. V_s \ud W_s$ is a BMO martingale, and $\NB{2}{V}$ is the BMO norm of $\int_0^. V_s \ud W_s$. We refer to \cite{Kazamaki-94} for further details about BMO martingales.

\subsection{Framework}

Let the domain $\cD$ be a bounded non-empty open connected subset of $\mathbb{R}^d$, and denote by $\partial \cD$ its boundary. Without loss of generality, we assume that $0 \in \cD$. For any $x \in \partial \cD$, we denote by $\mathfrak{n}(x)$ the normal exterior cone, i.e.,  the polar of the tangent cone, of $\cD$ at $x$, and denote by $\mathfrak{n}_u(x)$ its subset consisting of all unit vectors. We set $\mathfrak{n}(x)=\{0\}$ for all $x \in \cD$. 
We also denote by $P_{\bar{\cD}}$ the set-valued projection function onto $\bar{\cD}$ and define, for any $r >0$, the open $r$-neighborhood
\begin{equation*}
    \cD_r := \{ x \in \R^d | |x-P_{\bar \cD}(x)|<r\}.
\end{equation*}

\noindent In what follows, we need to refer to the following regularity property for a bounded domain $\mathfrak{D}\subset\R^d$:

\smallskip
\noindent \textbf{(R)} For all $x \in \partial \mathfrak{D}$, there exists a neighborhood of $x$ (in $\bar{\mathfrak{D}} $) that is $C^2$-diffeomorphic to a convex set of $\R^d$.

\smallskip
Importantly, the regularity property \textbf{(R)} spreads to product spaces as stated in the following proposition given without proof.
\begin{Proposition}\label{pr reg for product}
If $\mathfrak{D}$ and $\mathfrak{D}'$ are two domains in $\R^d$ that possess the regularity property \textbf{(R)}, then the domain $\mathfrak{D} \times \mathfrak{D}'$ in $\R^{2d}$ also possesses the regularity property \textbf{(R)}.
\end{Proposition}
\smallskip

\noindent \emph{Throughout the paper}, we impose the above regularity property on the domain $\cD$. Namely the following assumption holds true throughout the paper, even if not cited explicitly.
 \begin{Assumption}
    \label{ass:main}
    The domain $\cD$ has the regularity property \textbf{(R)}.
\end{Assumption}


\begin{Remark}
    (i) Assumption \ref{ass:main} is clearly fulfilled if for example $\cD$ is a $C^2$ domain (i.e. $\partial \cD$ is a $C^2$ manifold of dimension $d-1$), or if $\cD$ is a convex set. 
    \\
    (ii) According to Proposition \ref{pr reg for product}, $\cD\times\cD$ satisfies \textbf{(R)}. Note that, even if $\cD$ is a $C^2$-domain, the set $\cD \times \cD$ is not a $C^2$-domain.
\end{Remark}

\noindent Assumption \ref{ass:main} has important corollaries for $\cD$.

\begin{Proposition}
    \label{prop:extsphere:interiorcone}
    $ $
    \begin{itemize}
        \item  For all $x \in \partial \cD$, $\mathfrak{n}(x) \supsetneq \{0\}$.
        \item $\cD$ satisfies the interior cone condition: for all $x \in \partial \cD$, there exists $\varepsilon>$ and a closed cone $\mathcal{K}$, centered at $x$ and having non-empty interior, such that $\mathcal{K} \cap B(x,\varepsilon) \subset \bar{\cD}$.
        \item {There exist $\alpha>0$, $R>0$,  $y_1,\ldots,y_n\in \partial\cD$ and $a_1,\ldots, a_n\in\R^d$, such that: $|a_i|=1$ for all $i$, $\partial\cD\subset \bigcup_{i=1}^n B(y_i,R)$, and $\xi\cdot a_i\geq\alpha$ for all $\xi\in\mathfrak{n}_u(y)$, all $y\in \partial\cD\cap B(y_i,2R)$ and all $i$.}
        \item $\cD$ satisfies the exterior sphere property: i.e., there exists $R_0$ such that, for all $x \in \partial \cD$, $u \in \mathfrak{n}_u(x)$ and $x' \in \bar{\cD}$, 
        \begin{equation}
            \label{ineq:extsphere} 
            (x-x')\cdot u +\frac{1}{2R_0} |x-x'|^2 \geqslant 0.
        \end{equation}
    \end{itemize}
\end{Proposition}


\begin{Remark}
    Let us remark that \eqref{ineq:extsphere} is equivalent to saying that $B(x+R_0 u,R_0)\cap \bar \cD = \emptyset$. In other terms, we can roll a ball of radius $R_0$ all around $\bar \cD$.
\end{Remark}

\begin{Remark}\label{re on diff and func} Let $\mathfrak{D}$ be a domain with the regularity property \textbf{(R)}.\\
(i) In the remainder of the paper, a function $f:\bar{\mathfrak{D}}  \to \R$ is said to be differentiable at $x \in \partial \mathfrak{D}$ if there exists a linear operator $U_x : \R^d \to \R$ such that $|f(y)-f(x)-U_x(y-x)| = o(|y-x|)$ for all $y \in \bar{\mathfrak{D}}$. Then, the interior cone property implies the uniqueness of $U_x$ and allows us to define properly the class of $C^1$ functions on $\bar{\mathfrak{D}}$. By the same token, we define the class of functions $C^k$ on $\bar{\mathfrak{D}}$, or on $\bar{\mathfrak{D}} \cap U $ with an open $U\subset\R^d$, for all $k \in \mathbb{N}$. Note that we always consider the Euclidian topology in these definitions.
\\
(ii) Using the Whitney extension theorem, we can extend any function $\psi \in C^2(\bar{\mathfrak{D}})$ to a $C^2$ function on $\R^d$, thanks to the interior cone property (see, e.g., \cite{whitney1934differentiable}). This justifies the application of It\^o's formula to such functions $\psi$, which appears later in the paper.
\end{Remark}

\noindent \textbf{Proof of Proposition \ref{prop:extsphere:interiorcone}}.
Let us prove the first point.
We consider $x \in\partial \cD$, $U$ a neighborhood of $x$ and $\phi : U \to \phi(U)$ a $C^2$ diffeomorphism such that $\phi(\bar \cD \cap U)$ is convex. For this convex we know that $\mathfrak{n}(\phi(x)) \supsetneq \{0\}$. Let $\mathfrak{n}^\flat(\phi(x)) $ be the set of linear forms on $\R^d$ of the form $v\mapsto u\cdot v$, for $u\in\mathfrak{n}(\phi(x))$. Denoting by $u^\flat$ this linear form, the map $u\mapsto u^\flat$ from $\mathfrak{n}(\phi(x))$ to $\mathfrak{n}^\flat(\phi(x)) $ is a linear bijection.  Moreover $\alpha\in \mathfrak{n}^\flat(\phi(x)) $ if and only if for all interior direction $v$ of $\phi(\bar \cD \cap U)$ at $\phi(x)$, $\alpha(v)\le 0$ ($v$ is the speed at $\phi(x)$ of a $C^1$ curve staying for some time in $\phi(\bar \cD \cap U)$).   Denoting $T_{\phi(x)}\phi^{-1}$ the tangent map of $\phi^{-1}$ at $\phi(x)$, we have that $v$ is an interior direction of $\phi(\bar \cD \cap U)$ at $\phi(x)$ 
 if and only if $T_{\phi(x)}\phi^{-1}(v)$ is an interior direction of $\bar\cD$ at $x$. Consequently, $(\alpha\circ T_x\phi)\left(T_{\phi(x)}\phi^{-1}(v)\right)\le 0$ since $(T_{\phi(x)}\phi^{-1})^{-1} = T_x\phi$. Let us denote by $\alpha\mapsto \alpha^\sharp$ the inverse bijection of $u\mapsto u^\flat$. We proved that $\mathfrak{n}(x)=\left(\mathfrak{n}^\flat(\phi(x))\circ T_x\phi\right)^\sharp$. In particular, this set contains nonzero vectors.

For the second point, again we start from the fact that the property is satisfied for $\phi(\bar \cD \cap U)$ at $\phi(x)$ and get a cone  $\mathcal{K}'$ with center $\phi(x)$ such that $\mathcal{K}'\cap B(\phi(x),\varepsilon')\subset \phi(\bar \cD \cap U)$. Then  $\phi^{-1}(\mathcal{K}'\cap B(\phi(x),\varepsilon'))\subset \bar\cD$, and with compactness arguments and regularity of $\phi^{-1}$ it is easy to find $\mathcal{K}$ and $\varepsilon$ with the desired properties such that $\mathcal{K} \cap B(x,\varepsilon)\subset \phi^{-1}(\mathcal{K}'\cap B(\phi(x),\varepsilon'))$.

Let us now prove the third point. Since $\partial\cD$ is compact, it is sufficient to prove that for all $x\in\partial\cD$, there exists $\alpha>0$, $R>0$, $a\in \R^d$ such that $|a|=1$ and $\xi\cdot a\ge \alpha$ for all $\xi\in \mathfrak{n}_u(y)$, all $y\in \partial\cD\cap B(x, R)$. With the same notation as before, it is a well-known fact that by convexity of $\phi(\bar \cD \cap U)$,  there exists a closed cone $K\subset\R^d$ with center $0$, nonempty interior and spherical section, and a  neighborhood $V$ of $x$ in $\partial\bar\cD$ such that for all $y\in V$, all elements of $K$ are interior directions for $\phi(\bar \cD \cap U)$ at $\phi(y)$. As in the first step of the proof, we deduce that all vectors of $T_{\phi(y)}\phi^{-1}(K)$ are interior directions for $\cD$ at $y$. Moreover, since $\phi$ is $C^2$ and so $T\phi^{-1}$ is $C^1$, one easily checks that possibly by reducing $V$ there exists a  closed cone $K'\subset\R^d$ with center $0$, nonempty interior and spherical section,  such that for all $y\in V$, all elements of $K'$ are contained in $T_{\phi(y)}\phi^{-1}(K)$, and consequently are interior directions for $\bar \cD $ at $y$. If $u$ is for instance the unitary central vector of $K'$, then for all $y\in V$ and $v\in \mathfrak{n}(y)$, $u\cdot v\le 0$. Consequently, $a:=-u$ answers the question. 

It remains to prove the exterior sphere property. Once again, Since $\partial\cD$ is compact it is sufficient to prove that for all $x \in \partial \cD$, there exists $R_0>0$ and $U$ a neighborhood of $x$ such that, for all $y \in \partial \cD \cap U$, $\forall u \in \mathfrak{n}_u(y)$, $\forall x' \in \bar{\cD}$, we have
\begin{equation}
    \label{ineq:extsphere:proof}
(x-x')\cdot u +\frac{1}{2R_0} |x-x'|^2 \geqslant 0.
\end{equation}
We set $x \in \partial \cD$ and we consider $V$ a neighborhood of $x$ and $\phi : V \to \phi(V)$ a $C^2$ diffeomorphism such that $\phi(\bar \cD \cap V)$ is convex. We set $R$ such that $\bar{B}(\phi(x),2R)\subset \phi(V)$. Since a convex domain satisfies the exterior sphere property for any radius, we have that for all $z \in \phi(\partial \cD) \cap \bar B(\phi(x),R)$, $\forall u \in \mathfrak{n}_u(z)$, $\forall x' \in \phi(\bar \cD \cap V)$, \eqref{ineq:extsphere:proof} is satisfied. In particular, for all $z \in \phi(\partial \cD) \cap \bar B(\phi(x),R)$, $\forall u \in \mathfrak{n}_u(z)$, $B(z+R u,R) \subset \phi(V \setminus \bar\cD)$ and $\bar{B}(z+R u,R) \subset \phi(V \setminus \cD)$. 
Moreover, $\phi^{-1}(B(z+R u,R))$ is a $C^2$ compact domain since $\phi^{-1}$ is $C^2$. In particular it means that it satisfies the interior sphere property, with a positive radius denoted $R_{z}>0$. Moreover, some elementary but tedious computations show that we can choose any $R_{z} \in ]0,M_z]$ where $M_z$ is a continuous function of $R$ and second derivatives of $\phi^{-1}$ at $z$. Since second derivatives of $\phi^{-1}$ are continuous, we can set $R_0:=\inf_{z \in \phi(\partial \cD) \cap \bar B(\phi(x),R)} M_z>0$. Thus, we can set $U:= \phi^{-1}(B(\phi(x),R))$ which is a neighborhood of $x$ and we have that for all $y \in \partial \cD \cap U$, $\forall u \in \mathfrak{n}_u(y)$, $\forall x' \in \bar \cD \cap U $, \eqref{ineq:extsphere:proof} is satisfied. Since $\bar B(y+uR_0,R_0) \subset U \setminus \cD$, \eqref{ineq:extsphere:proof} is also satisfied for all $x'$ that are in $\bar \cD \setminus U$ which concludes the proof. \eproof

\smallskip

The exterior sphere property has useful corollaries.
\begin{Lemma}
    \label{lem:proj:lip}
    $P_{\bar{\cD}}$ is a single-valued function on $\cD_{R_0}$. Moreover, for all $r \in (0,R_0)$, $P_{\bar{\cD}}$ is Lipschitz on $\cD_r$, with the Lipschitz constant $\frac{R_0}{R_0-r}$.
\end{Lemma}

\proof
The first part of the Lemma is direct, see e.g. Corollary 2.1 in \cite{Chassagneux-Nadtochiy-Richou-22}. Then, Theorem 4.1 in \cite{Poliquin-Rockafellar-Thibaullt-00} and Theorem 4.8 in \cite{Clark-Stern-Wolenski-95} allow us to conclude.
\eproof

\smallskip

Next, since the set $\bar{\cD}$ is flat, we define the length of any absolutely continuous curve in $\bar\cD$ as the standard Euclidean setting, i.e., as an integral of the absolute velocity of this curve.
Then, we define a geodesic between two points $x$ and $y$ in $\bar{\cD}$ as an absolutely continuous path $\gamma : [0,1] \rightarrow \bar{\cD}$ such that $\gamma_0=x$, $\gamma_1=y$, $\gamma$ has a constant speed (i.e., $|\dot{\gamma}|$ is constant and equal to the inverse of the length of $\gamma$), and such that, locally, $\gamma$ is a distance-minimizing curve.

\smallskip 

We make the following observations for the space $\bar{\cD}$, again implied by Assumption \ref{ass:main}.
\begin{Theorem} $\bar{\cD}$ is a \emph{geodesic space}: namely,
for any $x,y \in \bar{\cD}$, there exists at least one minimizing geodesic between $x$ and $y$. Moreover, all minimizing geodesic between $x$ and $y$ are $C^1$.
\end{Theorem}
We refer to Corollary III page 48 in \cite{Wolter-79} for a proof.

\begin{Remark}
    \label{rem:geodesic:reg}
The geodesics of $\bar{\cD}$ are not necessarily $C^2$, even if the boundary $\partial\cD$ is $C^2$. If the boundary is $C^2$, the geodesics can be decomposed into 
 \begin{itemize}
  \item geodesic segments of $\cD$, whose acceleration vanishes,
  \item  geodesic segments of the boundary $\partial\cD$, whose acceleration is outwardly normal to $\partial \cD$,
  \item switch points, where the geodesic switches from a boundary segment to an interior segment and vice-versa,
  \item intermittent points, which are the accumulation points of the switch points.
 \end{itemize}
 We refer to \cite{Alexander-Berg-Bishop-86,Albrecht-Berg-91} for the examples of pathological and good behavior behavior of the geodesics under stronger assumptions on the boundary. 
\end{Remark}

\section{Martingales, with and without a drift, on a Euclidean manifold with a boundary}

In order to make use of the stochastic geometry tools, we view the domain $\bar{\cD}$ as a Euclidean manifold with a boundary. Indeed, we have the following properties: 
\begin{itemize}
\item For any point $x \in \cD$, there exists a neighborhood of $x$ that is equal, hence trivially isometric, to an open subset of $\mathbb{R}^d$. 
\item  For any point $x \in \partial\cD$, there exists a neighborhood of $x$ that is homeomorphic to an open subset of $\mathbb{R}^{d-1}\times \mathbb{R}^+$. 
\end{itemize}
We remark that the assumptions on $\bar{\cD}$ are not sufficient to replace the homeomorphism by a diffeomorphism: in other words, the manifold is not necessarily a differential manifold, due to a potential lack of regularity of the boundary.

\smallskip
\noindent Let us us start by adapting the notion of a martingale on a manifold, classically called a $\Gamma$-martingale, to our setting.


\begin{Definition}
    \label{def:gamma-mart-0}
 Let $X$ be a continuous (Euclidean) semimartingale with values in $\bar{\cD}$. Then, $X$ is a $\Gamma$-martingale on $\bar{\cD}$ if 
 $X=X_0+K+M$, with $M$ being an $\mathbb{R}^d$-valued local continuous martingale and with $K$ being an $\mathbb{R}^d$-valued continuous process such that
 $$K_t = \int_0^t  k_s \ud\!\text{Var}_s (K),$$
 where $(k_s)_{s \in [0,T]}$ is a progressively measurable process satisfying
 $$k_t \in \mathfrak{n}(X_t), \quad \text{for a.e. }  t \in [0,T].$$ 
\end{Definition}

In order to treat general reflected BSDEs, we also need to define the notion of a $\Gamma$-martingale with a drift on a Euclidean manifold with a boundary.

\begin{Definition}
\label{def:gamma-mart}
 Let $f$ be an element of $\HP{1}$, and let $X$ be a continuous (Euclidean) semimartingale with values in $\bar{\cD}$. Then, $X$ is a $\Gamma$-martingale on $\bar{\cD}$ with the drift $f$ if $X=X_0-F+ K+M$, with $F=\int_0^\cdot f_s\,\ud s$, with $M$ being an $\R^d$-valued local continuous martingale, and with $K$ being an $\mathbb{R}^d$-valued adapted continuous process such that
 $$K_t = \int_0^t k_s\ud \!\text{Var}_s (K),$$
 where $(k_s)_{s \in [0,T]}$ is a progressively measurable process satisfying
 $$k_t \in \mathfrak{n}(X_t), \quad \text{for a.e. }  t \in [0,T].$$ 
\end{Definition}

\begin{Remark}
    \label{rem link bsde}
        Let us consider the following reflected BSDE 
        \begin{align}
            \label{eq reflected bsde}
             \left\{ \begin{aligned} &(i)\;Y_t = \xi+\int_t^T f(s,Y_s,Z_s)\ud s-\int_t^T \ud K_s-\int_t^T Z_s \ud W_s, \; 0 \leqslant t \leqslant T,\\
            & (ii)\; Y_. \in \bar{\cD} \text{ a.s.}, \quad { \dot{K}_\cdot \in \mathfrak{n}(Y_\cdot) \;\ud t \otimes \ud \P-a.e.},\quad \int_0^T \1_{\set{Y_s \notin \partial \cD}} \ud \mathrm{Var}_s(K) = 0,
                    \end{aligned}
             \right.
        \end{align}
        where $\xi \in \mathscr{L}^{\infty}$, $f(.,y,z)$ is a progressively measurable process for all $(y,z) \in \bar{\cD} \times \mathbb{R}^{d \times d'}$. By definition, a solution of the reflected BSDE \eqref{eq reflected bsde} is the triple of processes $(Y,Z,K)$ that satisfies:
        \begin{enumerate}
         \item $Y \in \mathscr{S}^{\infty}$,
         \item $\int_0^T |Z_s|^2 \ud s<+\infty$ a.s.,
         \item $K \in \mathscr{K}^1$,
         \item $\int_0^T |f(s,Y_s,Z_s)|\ud s<+\infty$ a.s.
        \end{enumerate}
        According to Definition \ref{def:gamma-mart}, for any solution $(Y,Z,K)$ of the reflected BSDE \ref{eq reflected bsde}, the process $Y$ is indeed a $\Gamma$-martingale on $\bar{\cD}$, with the drift $(f(s,Y_s,Z_s))_{s\in[0,T]}$ and with the terminal value $\xi$.  
\end{Remark}

\medskip

As in the case of manifolds without boundary, we introduce the notion of a $\Gamma$-convex function in order to characterize $\Gamma$-martingales.


\begin{Definition}
    \label{def:gamma-convex}
Consider $\psi:\bar\cD \rightarrow \R$ and let $U$ be an open subset of $\mathbb{R}^d$ s.t. $ U\cap\bar\cD \neq \emptyset$.
We say that $\psi_{|\bar{\cD} \cap U}$  is a $\Gamma$-convex function if, for any geodesic curve $\gamma$ on $\bar{\cD} \cap U$, the function $\psi\circ\gamma:\,[0,1]\rightarrow\R$ is convex.
\end{Definition}

\noindent We also consider an important subset of $\Gamma$-convex functions.

\begin{Definition}
    \label{de spe gamma func}
  Consider $\psi \in C^1(\bar{\cD})$ and let $U$ be an open subset of $\mathbb{R}^d$ s.t. $ U\cap\bar\cD \neq \emptyset$.  
We say that $\psi_{|\bar{\cD} \cap U}$ is a \emph{special $\Gamma$-convex function} if it is $\Gamma$-convex and, for all \( x \in \partial \cD \cap U \) and all \( v \in \mathfrak{n}(x) \), it holds that  
\[
\nabla \psi(x) \cdot v \ge 0.
\]
We say that $\psi$ is a \emph{global special $\Gamma$-convex function} if the above property holds with \( U = \mathbb{R}^d \).  
\end{Definition}


\begin{Remark}\label{rem:global.vs.local.GammaConvex}
    \begin{enumerate}
    \item For a convex function $\psi : \mathbb{R}^d \rightarrow \mathbb{R}$, the restriction $\psi_{|\bar \cD}$ is not necessarily a $\Gamma$-convex function. {For example, consider $U=B(0,1)$ and $\bar \cD\cap U = U\setminus B(x_c,1)$, where $x_c=(0,...,0,1)$, and set $\psi(x)=\sum_{i=1}^d \alpha_i x_i$. Then, by considering the geodesic $t \mapsto (1-\sqrt{1-(t/2)^2},0,...,0,t/2)$, we deduce that $\psi$ is not a $\Gamma$-convex function if $\alpha_d<0$.}
    \item A special $\Gamma$-convex function is of course a $\Gamma$-convex function, but the converse is not true in general. {Indeed, let us consider $U=B(0,\varepsilon)$ and $\bar \cD\cap U = U\setminus B(x_c,1)$. Then, one can show that $\psi(x)=-x_d+\sum_{i=1}^{d-1} x_i^2$ is a $\Gamma$-convex function for $\varepsilon>0$ small enough, but it is not a special $\Gamma$-convex function since $\nabla \psi(0)=(0,...,0,-1)$ and $(0,...,0,1)\in \mathfrak{n}(0)$.}
    \item 
    $\Gamma$-convex functions and special $\Gamma$-convex functions play the same role as classical convex functions in the Euclidean space $\mathbb{R}^d$. They are used in the subsequent parts of the paper as test functions to characterize $\Gamma$-martingales, in analogy with the Euclidean case: in $\mathbb{R}^d$, a process $X$ is a martingale if and only if, for every convex function $\psi$, the process $\psi(X)$ is a real submartingale (under suitable integrability assumptions).

A key difference from the Euclidean setting is that the $\Gamma$-convex and special $\Gamma$-convex functions are defined only locally. This localization arises because, in some situations, there are not enough global test functions to fully characterize a $\Gamma$-martingale. For instance, if the domain $\cD$ contains a ``hole" (e.g. a missing ball), one can show that all global $\Gamma$-convex and special $\Gamma$-convex functions must be constant on the boundary of that ball. Such functions are therefore insufficient to determine the reflection direction of a semimartingale in $\bar{\cD}$ (see the proof of Proposition~\ref{prop caracterisation gamma martingale}).

On the other hand, it is clearly more convenient to work with global $\Gamma$-convex and special $\Gamma$-convex functions, that is, by taking $U = \mathbb{R}^d$ in the definition above. For some classes of domains $\cD$, such a restriction is indeed sufficient to characterize $\Gamma$-martingales; see Corollary~\ref{global gamma mart char}.
    
    \end{enumerate}
\end{Remark}

\noindent Special $\Gamma$-convex functions play a key role in the characterization of $\Gamma$-martingales.

\begin{Proposition}
\label{prop gamma-function}
{
Consider $\psi \in C^2(\bar{\cD})$ and let $U$ be an open subset of $\mathbb{R}^d$ s.t. $ U\cap\bar\cD \neq \emptyset$. Assume that $\psi_{|U\cap\bar{\cD}}$ is a special $\Gamma$-convex function according to Definition \ref{de spe gamma func}.}
 Then $\nabla^2 \psi \geqslant 0$ on $\bar{\cD} \cap U$. Moreover, if $X$ is a $\Gamma$-martingale with a drift $f$, then the {finite-variational component} of the real semimartingale $\psi(X_t)+\int_0^t \nabla \psi(X_s)\cdot f_s \,\ud s$ is {a.s.}  non-decreasing in the random  open  set $\{t:\, X_t \in U \cap \bar{\cD}\}$.
\end{Proposition}

\proof
Let us consider $x \in \cD$. Then there exists $\varepsilon >0$ such that $\bar{B}(x,\varepsilon) \subset \cD$, and for all $y \in \bar{B}(x,\varepsilon)$, $\gamma : [0,1] \rightarrow \bar{B}(x,\varepsilon)$  given by $\gamma_t = x+(y-x)t$ is a geodesic. Then $\varphi(.):=\psi(\gamma_.)$ is a $C^2$ convex function and we get the result since $$\varphi''(0) = (y-x)^\top \nabla^2 \psi(x) (y-x) \geqslant 0.$$
{If $x \in \partial \cD$, we just have to use the continuity of $\nabla^2 \psi$ in order to conclude.} It proves the first part of the proposition.

Let us now consider $X$ a $\Gamma$-martingale with drift $f$. It\^o formula\footnote{Have in mind Remark \ref{re on diff and func}(ii).
} yields us 
$$d\psi(X_t) +\nabla \psi(X_t). f_t dt =\nabla \psi(X_t) dK_t +\nabla \psi(X_t) dM_t + \frac{1}{2} \langle \nabla^2\psi(X_t) dM_t,dM_t \rangle,$$ 
and the compensator satisfies
\begin{align*}
    & \1_{X_t \in \bar{\cD} \cap U}\left( \nabla \psi(X_t) dK_t  + \frac{1}{2} \langle \nabla^2\psi(X_t) dM_t,dM_t \rangle \right)\\
    =&  \1_{X_t \in \bar{\cD} \cap U} \nabla \psi(X_t).k_t \ud \text{Var}_t (K)+ \1_{X_t \in \bar{\cD} \cap U} \frac{1}{2} \langle \nabla^2\psi(X_t) dM_t,dM_t \rangle \geqslant 0.
\end{align*}
\eproof

Our goal now is to obtain a converse statement, i.e. to characterize $\Gamma$-martingales through the special $\Gamma$-convex functions.

\begin{Proposition}
\label{prop caracterisation gamma martingale}
 Let $X$ be a {continuous adapted process} with values in $\bar{\cD}$, and let $(f_s)_{s \geqslant 0}$ {be an element of $\HP{1}$}. Assume that, 
 {for any open set $U \subset \R^d$ s.t. $U\cap\bar{\cD}\neq \emptyset$ and any $C^2$ function $\psi$ on $\bar{\cD}$, such that $\psi_{|U\cap\bar{\cD}}$ is a special $\Gamma$-convex function}, the {finite-variational component} of the real semimartingale $\psi(X_t)+\int_0^t \nabla \psi(X_s)\cdot f_s \,\ud s$ is {a.s.}   non-decreasing in the random  open set $\{t:\, X_t \in U \cap \bar{\cD}\}$. Then, $X$ is a $\Gamma$-martingale with the drift $f$. 
\end{Proposition}

\proof
{
{\bf Step 1.} We claim that, for any $x \in \partial \cD$, there exists a non-empty open neighborhood $\mathcal{O}_x$ of $x$ and a basis $\{e_j^x\}_{1 \le j \le d}$ (not necessarily orthogonal) such that, for any $1\le j \le d$, there exists a completion $(e_1,\ldots,e_{d-1},e_j^x)$ of $e^x_j$ to an orthonormal basis, in which the set $\bar{\cD}\cap\mathcal{O}_x$ can be represented as a sub-graph of a function of the first $d-1$ coordinates.

To prove this claim, we fix an arbitrary $x \in \partial \cD$ and consider a diffeomorphism $\Phi$ and a non-empty open neighborhood $\tilde{\mathcal{O}}$ of $x$ such that $\Phi$ maps $\tilde{\mathcal{O}}\cap\cD$ into a convex subset $\tilde D$ of $\R^d$ (the existence of such $\Phi$ and $\tilde{\mathcal{O}}$ is guaranteed by Assumption \ref{ass:main}). Without loss of generality (since an affine transformation does not change the convexity properties) we assume that the Jacobian of $\Phi$ at $x$ equals the identity, that $\Phi(x)=x$ and that $x=0$. We denote by $K$ the smallest cone centered at the origin that contains $\tilde D$ (it is well defined as $\tilde D$ is convex). Then, we define $E$ as the intersection of the interior of $K$ (which is non-empty since $\tilde D$ is open) with the negative of the normal exterior cone of $\tilde D$ at the origin. It is easy to see that the latter intersection is non-empty and that, for any $z_0\in E$, we can find an orthonormal basis $(\tilde e_1,\ldots,\tilde e_d)$ of $\R^d$ such that $\tilde e_d$ points from the origin to $z_0$ and such that, in the coordinates associated with this basis, the image of $\partial\cD$ under $\Phi$ can be viewed locally as a graph of a convex function:
\begin{align*}
\Phi(\partial\cD\cap\tilde{\mathcal{O}}) = \{(\bar x,y)\in\R^d:\,y=\tilde F(\bar x),\,|\bar x|<\tilde\epsilon\},
\end{align*}
where we have reduced $\tilde{\mathcal{O}}$, if needed, and introduced a convex function $\tilde F:\,\{\bar x\in \R^{d-1}:\,|\bar x|<\tilde\epsilon\}\rightarrow\R$ along with a constant $\tilde\epsilon>0$. Further reducing $\tilde{\mathcal{O}}$ and $\tilde\epsilon>0$, we can assume that $\tilde F$ is uniformly Lipschitz in its domain.

Next, consider a unitary linear operator $A$ on $\R^d$ with the property that $\|A-I\| \leq \varepsilon$, for some $\varepsilon>0$. Denoting by $(\bar x',y')$ the image of a point $(\bar x,y)$ under such a mapping $A$, we deduce the existence of $a_0\in\R$, $\bar a,\bar b\in\R^{d-1}$ and $\bar A \in \R^{(d-1)\times(d-1)}$ such that
\begin{align*}
& y = y' + a_0\,y' + \bar a^\top\,\bar x',\quad \bar x = \bar A\,\bar x' + \bar b\,y',
\end{align*}
and the norms of $(a_0,\bar a,\bar b)$ can be made arbitrarily small by the choice of $\varepsilon>0$.
Then, $\Phi(\partial\cD\cap\tilde{\mathcal{O}})$ is described via
\begin{align*}
&y' = \tilde F(\bar A\,\bar x' + \bar b\,y') - a_0\,y' - \bar a^\top\,\bar x',
\end{align*}
for $|\bar A\,\bar x' + \bar b\,y'|<\tilde\epsilon$.
Choosing a sufficiently small $\varepsilon>0$, we ensure that the Lipschitz coefficient of the right hand side of the above equation, viewed as a function of $y'$, is small enough, so that there exists $\hat\epsilon>0$ such that, for any fixed $\bar x'$ with $|\bar x'|<\hat\epsilon$, there exists exactly one $y'$ that satisfies the above equation. The latter means that there exists a function $\hat F$ such that
\begin{align*}
A\circ \Phi(\partial\cD\cap\tilde{\mathcal{O}}) = Q := \{(\bar x',y')\in\R^d:\,y'=\hat F(\bar x'),\,|\bar x'|<\hat\epsilon\},
\end{align*}
where we reduce $\tilde{\mathcal{O}}$ as needed.
Using the Lipschitz property of $F$ and the smallness of $(a_0,\bar a,\bar b)$, it is easy to deduce that $\hat F$ is also Lipschitz.
Next, we choose operators $A_1,\ldots,A_d$, with the associated $\varepsilon>0$ being small enough, so that the above representation of the image of $\Phi(\partial\cD\cap\tilde{\mathcal{O}})$ under each $A_i$ holds with a Lipschitz function $\hat F_j$, a set $Q_j$, and with, possibly, smaller $\tilde{\mathcal{O}}$ and $\tilde\epsilon>0$, so that each $A_j$ is invertible, and so that $\{A_j^{-1}\,\tilde e_d\}_{j=1}^d$ are linearly independent. For each $j=1,\ldots,d$, we define $e_j^x:=A_j^{-1}\,\tilde e_d$ and complete it (to form an orthonormal basis) with $\{A_j^{-1}\,\tilde e_i\}_{i=1}^{d-1}$.

\smallskip

It remains to show that each $e^x_j$ has the desired properties. Without loss of generality, we only consider $j=1$.
Next, we notice that $\tilde\Phi:=A_1\circ\Phi\circ A^{-1}_1$ is a $C^2$ diffeomorphism that maps the origin into itself and whose Jacobian at the origin is the identity (i.e., it inherits these properties from $\Phi$). Then, using the representation $\tilde\Phi=(\tilde\Phi_1,\tilde\Phi_2)$, we obtain
\begin{align*}
&A_1(\partial\cD\cap\tilde{\mathcal{O}})=A_1\circ\Phi^{-1}\circ A^{-1}_1(Q_1) = \tilde\Phi^{-1}(Q_1) \\
&= \{(\bar x,y)\in\R^d:\,\tilde\Phi_2(\bar x,y)=\hat F_1(\tilde\Phi_1(\bar x,y)),\,|\tilde\Phi_1(\bar x,y)|<\hat\epsilon\}.
\end{align*}
Since $\tilde\Phi(0)=0$ and the Jacobian of $\tilde\Phi$ at the origin equals identity, we have
\begin{align*}
& \nabla_y\tilde\Phi_1(\bar x,y) = O(|\bar x|+|y|),\quad  \nabla_y\tilde\Phi_2(\bar x,y) = 1 + O(|\bar x|+|y|).
\end{align*}
Using the above and the Lipschitz property of $\hat F_1$, we conclude that there is a small enough $\epsilon>0$ such that, for any $|\bar x|<\epsilon$, there exists a unique fixed point of the mapping $y$ to the equation
\begin{align*}
y\mapsto y + \hat F_1(\tilde\Phi_1(\bar x,y)\wedge \hat\epsilon) - \tilde\Phi_2(\bar x,y),
\end{align*}
with $\tilde\Phi_1(\bar x,y)< \hat\epsilon$.
The above yields the existence of a non-empty open neighborhood $\mathcal{O}_x$ of $x$ and a function $F_1:\,B(0,\epsilon)\rightarrow \R$ such that
\begin{align*}
&A_1(\partial\cD\cap\mathcal{O}_x) = \{(\bar x,y)\in\R^d:\,y=F_1(\bar x),\,|\bar x|<\epsilon\}.
\end{align*}

\smallskip

Finally, to complete Step 1, we show that
\begin{align*}
&\text{either}\quad A_1(\cD\cap\mathcal{O}_x) = H_1:= \{(\bar x,y)\in\R^d:\,y<F_1(\bar x),\,|\bar x|<\epsilon\}\cap A_1(\mathcal{O}_x)\\
&\text{or}\quad A_1(\cD\cap\mathcal{O}_x) = H_2:= \{(\bar x,y)\in\R^d:\,y>F_1(\bar x),\,|\bar x|<\epsilon\}\cap A_1(\mathcal{O}_x).
\end{align*}
To this end, we notice that $A_1(\cD\cap\mathcal{O}_x)$ has a non-empty intersection with either $H_1$ or $H_2$. Without loss of generality, we assume that $A_1(\cD\cap\mathcal{O}_x)\cap H_1\neq\emptyset$.
Reducing $\mathcal{O}_x$, we can assume that the latter set is a small enough open cylinder centered around $e_j^x:=A_j^{-1}\,\tilde e_d$.
Then, it is easy to see that $H_1$ is connected. 
We claim that $H_1\subset A_1(\cD\cap\mathcal{O}_x)$. Indeed, if there exists $z_1\in H_1\setminus A_1(\cD\cap\mathcal{O}_x)$, then we can connect it to $z_2\in A_1(\cD\cap\mathcal{O}_x)\cap H_1$ via a continuous curve that stays inside $H_1$. It is easy to see that this curve must intersect the boundary of $A_1(\cD\cap\mathcal{O}_x)\cap H_1$. Since this intersection point, denoted $z_3$, cannot belong to the $\partial H_1\cup\partial A_1(\mathcal{O}_x)$, it must belongs to the boundary of $A_1(\cD)$, which coincides with $A_1(\partial\cD)$. Since $z_3\in H_1\subset A_1(\mathcal{O}_x)$, we conclude that $z_3\in A_1(\partial\cD\cap\mathcal{O}_x)\cap H_1=\emptyset$, which is a desired contradiction.
Similarly, we deduce that $H_2\subset A_1(\mathcal{O}_x\setminus \bar{\cD})$, thus, completing Step 1.

\smallskip

{\bf Step 2.} We claim that, for any $1\le j \le d$, the function $u_j^x:\,y \mapsto \langle y,e_j^x \rangle$ is a special $\Gamma$-convex function in $\mathcal{O}_x$. Without loss of generality, we consider $j=1$ and assume $x=0$.

First, we show that $u_1^x$ is a Gamma-convex function. To this end, we recall the local representation of $\cD$ obtained in Step 1 and deduce the existence of an orthonormal basis $(e_1,\ldots,e_{d-1},e_1^x)$, such that the set $\cD\cap\mathcal{O}_x$ (with $\mathcal{O}_x$ being a small cylinder centered around $e_1^x$), written in the coordinates induced by this basis, is given by
\begin{align}
&\{(\bar x,y)\in\R^d:\,C_1<y<F_1(\bar x),\,|\bar x|<\epsilon\},\label{eq.prop2.2.pf.new.1}
\end{align}
with a constant $C_1\in\R$.
Consider an arbitrary geodesic curve $\gamma$ in $\bar{\cD}\cap\mathcal{O}_x$ and notice that, in the new coordinates, we have $\langle \gamma_t,e_1^x \rangle$ equals the last coordinate of $\gamma_t=(\gamma^1_t,\ldots,\gamma^d_t)$. Then, it suffices to show that the function $t\mapsto \gamma^d_t$ is convex. We argue by contradiction and assume that $t\mapsto \gamma^d_t$ is not convex. Notice that, since $\gamma$ is locally distance-minimizing, there exists $\varepsilon>0$ such that, for any $t\leq s\in[0,1]$ with $|t-s|<\varepsilon$, the arc $\gamma_{[t,s]}$ is a minimizing geodesic connecting $\gamma_t$ and $\gamma_s$. Then, the lack of convexity of $t\mapsto \gamma^d_t$ implies the existence of $0\leq t_1<t_2 < t_3\leq 1$, with $|t_3-t_1|<\varepsilon$, such that { $\gamma^d_{t_2}>\gamma^d_{t_1}+\frac{\gamma^d_{t_3}-\gamma^d_{t_1}}{t_3-t_1}(t_2-t_1)$  }. Let us consider a new curve $\tilde\gamma$ that coincides with $\gamma$ outside of $[t_1,t_3]$, and for any $t\in[t_1,t_3]$, we have {$\tilde\gamma_t:=(\gamma^1_t,\ldots,\gamma^{d-1}_t,\gamma_{t_1}^d+\frac{\gamma^d_{t_3}-\gamma^d_{t_1}}{t_3-t_1}(t-t_1))$}.
Using the representation \eqref{eq.prop2.2.pf.new.1}, we conclude that $\tilde\gamma$ is in $\bar{\cD}\cap\mathcal{O}_x$. On the other hand, it is clear that the length of $\tilde\gamma_{[t_1,t_3]}$ is strictly smaller than that of $\gamma_{[t_1,t_3]}$, which contradicts the fact that $\gamma_{[t,s]}$ is a minimizing geodesic connecting $\gamma_t$ and $\gamma_s$, and completes the proof that $u_1^x$ is a Gamma-convex function.

\smallskip

To conclude Step 2, it remains to show that, for any $z\in \partial\cD\cap\mathcal{O}_x$ and any $v\in\mathfrak{n}(z)$, we have $e_1^x\cdot v\geq0$. We work in the coordinates induced by the basis $(e_1,\ldots,e_{d-1},e_1^x)$ {and we consider $z\in \partial\cD\cap\mathcal{O}_x$ and $v=(v_1,\ldots,v_d)\in\mathfrak{n}(z)$.} 
Since the point $z$ admits an exterior sphere (see Proposition \ref{prop:extsphere:interiorcone}), we must have {$(z-x')\cdot v + \frac{1}{R_0} |z-x'|^2\geq0$ for any $x'\in\bar\cD$. It remains to notice that, due to the representation of $\cD\cap\mathcal{O}_x$ via \eqref{eq.prop2.2.pf.new.1}, there exists $\varepsilon_0>0$ such that $x':=(z_1,\ldots,z_{d-1},z_d-\varepsilon)\in \bar\cD$ for all $\varepsilon \in (0,\varepsilon_0)$, which yields:
\begin{align*}
0\leq(z-x')\cdot v + \frac{1}{R_0} |z-x'|^2 = \varepsilon \left(v_d+ \frac{\varepsilon}{R_0}|v_d|^2\right).
\end{align*}
Then, by taking $\varepsilon \to 0$, we conclude that $v_d \geq 0$.
}
}

\smallskip

{\bf Step 3.} Let us prove that $X$ is a semimartingale in $\mathbb{R}^d$. 
\\
{\bf Step 3.a} By compactness we can extract a finite set of neighborhood, denoted $(\mathcal{O}_i) _{1\le i \le I}$ by a slight abuse of notation, such that  $\partial \cD \subset\bigcup_{1 \le i\le I} \mathcal{O}_i$. By the same slight abuse of notation, we denote $(u_j^i)_{1 \le j \le d}$ the special  $\Gamma$-convex functions associated to $\mathcal{O}_i$. We set $O_0 := \cD$, $(e^0_i)_{1 \le i\le d}$ the canonical orthonormal basis and $(u_j^0)_{1 \le j \le d}$ the special  $\Gamma$-convex functions associated. 
\\
{\bf Step 3.b} {Let us consider now a sequence of stopping time $(\tau_n)_{n \in \mathbb{N}}$ such that $X_t \in \mathcal{O}_{i^n}$ for all $t \in [\tau_n,\tau_{n+1})$. Then, for all $n \in \mathbb{N}$, $u_{j}^{i^n}(X)+\int_0^. \nabla u_j^{i^n}(X_s) f_s \ud s$ is a semimartingale on $[\tau_n,\tau_{n+1})$ for all $1 \le j \le d$ which implies that $X$ is a semimartingale on $[\tau_n,\tau_{n+1})$ and gives us the announced result.}

\smallskip

{\bf Step 4.} We now prove that $X$ is a $\Gamma$-martingale  with drift $f$.  Let us denote by $K$ the {finite-variational component} of the real semimartingale $\psi(X_t)+\int_0^t \nabla \psi(X_s).f_s ds$.
\\
{\bf Step 4.a} We start by taking $U = \cD$. Then, for all $u \in \mathbb{R}^d$, $x \rightarrow u \cdot x$ is a special $\Gamma$-convex function on $U$. In particular, $\1_{X_t \in \cD} dK^i_t$ and $-\1_{X_t \in \cD} dK^i_t$ are increasing for all $1 \leqslant i \leqslant d$ which allows to conclude that $dK_t = \1_{X_t \in \partial \cD} dK_t$.
\\
{\bf Step 4.b} By writing $K = \int_0^. k_s \ud \text{Var}_s (K)$, it remains to prove that $k_t \in \mathfrak{n}(X_t)$ for a.e. $t \in [0,T]$, in order to conclude. The compactness of $\bar{\cD}$, gives us that, for all $n \in \mathbb{N}^*$, there exists a finite family $(x_i^n)_{i \in I_n}$ of elements of $\partial \cD$ such that $\partial \cD \subset \bigcup_{i \in I_n} B(x_i^n,1/n)$. If we take $U = B(x_i^n, 1/n)$, then $x \mapsto u \cdot x$ is a special $\Gamma$-convex function on $U$  as soon as
$$u \in \mathcal{C}^{i,n}:=\{ u \in \mathbb{R}^d | u \cdot v \geqslant 0, \forall x \in \partial \cD \cap B(x_i^n, 1/n), \forall v \in \mathfrak{n}(x)\}.$$
It implies that, for all $n\in \mathbb{N}^*$, for all $t \in [0,T]$ , $k_t \in (\tilde{\mathcal{C}}^{X_t,n})^*$ when $X_t \in \partial \cD$,
with 
$$\tilde{\mathcal{C}}^{x,n} :=\bigcap_{i \in I_n, x \in B(x_i^n,1/n)} \mathcal{C}^{i,n}, \quad \forall x \in \partial \cD,\forall n \in \mathbb{N}^*,$$
and where the superscript $^*$ denotes the dual cone, i.e. for a cone $\mathcal{C}$
$$(\mathcal{C})^* := \{v \in \mathbb{R}^d | v \cdot u \geqslant 0, \forall u \in \mathcal{C} \}.$$
Hence $k_t \in \bigcap_{n \in \mathbb{N}^*} (\tilde{\mathcal{C}}^{X_t,n})^*$. Let us remark that 
$$\tilde{\mathcal{C}}^{x,n} \supset {\mathcal{C}}^{x,n} := \{ u \in \mathbb{R}^d | u \cdot v \geqslant 0, \forall y \in \partial \cD \cap B(x, 2/n), \forall v \in \mathfrak{n}(y)\}=\left(\bigcup_{y \in B(x, 2/n)} \mathfrak{n}(y) \right)^*  ,$$
Then, it just remains to prove that, for all $x \in \partial \cD$, 
\begin{equation}
    \label{equality dual cone}
    \bigcap_{n \in \mathbb{N}^*} \left(\bigcap_{y \in B(x,2/n)} (\mathfrak{n}(y))^*\right)^* =\mathfrak{n}(x) .    
\end{equation}
 
Firstly, we have for all $n \in \mathbb{N}^*$, $\bigcap_{y \in B(x,2/n)} (\mathfrak{n}(y))^* \subset \mathfrak{n}(x)^*$ and then
$$\bigcap_{n \in \mathbb{N}^*} \left(\bigcap_{y \in B(x,2/n)} (\mathfrak{n}(y))^*\right)^* \supset \mathfrak{n}(x) .$$ Now, let us prove the other inclusion. 
To do it, it is sufficient to assume that $\bar{\cD} \cap B(x,2/n)$ is convex, at least for $n$ large enough. Indeed, if it is not the case, by assumption, there exists a $C^2$-diffeomorphism that sends $\bar{\cD} \cap B(x,2/n)$ to a convex set. We have 
\begin{equation}
    \label{equality dual cone 2}
    \bigcap_{n \in \mathbb{N}^*} \left(\bigcap_{y \in B(x,2/n)} (\mathfrak{n}(y))^*\right)^* = \bigcap_{n \in \mathbb{N}^*} \left(\left(\bigcup_{y \in B(x,2/n)} \mathfrak{n}(y)\right)^*\right)^* = \bigcap_{n \in \mathbb{N}^*} \bar{\text{co}} \left(\bigcup_{y \in B(x,2/n)} \mathfrak{n}(y)\right).
\end{equation}
 
By using Corollary 24.5.1 in \cite{Rockafellar-70}, for all $\varepsilon>0$, there exists $n$ large enough such that
$$\bigcup_{y \in B(x,2/n)} \mathfrak{n}(y) \subset \mathfrak{n}(x) + B(0,\varepsilon).$$
Then, \eqref{equality dual cone 2} becomes

$$\bigcap_{n \in \mathbb{N}^*} \left(\bigcap_{y \in B(x,2/n)} (\mathfrak{n}(y))^*\right)^* \subset \bigcap_{\varepsilon >0} \left( \mathfrak{n}(x) + \bar{B}(0,\varepsilon) \right) =\mathfrak{n}(x)$$ 
which proves \eqref{equality dual cone}.
\eproof



\medskip

The next Lemma shows how to extend a (local) special $\Gamma$-convex function to a global one.  For a $\Gamma$-convex function $\psi$ on  $U\cap \bar{\cD}$, where $U$ is an open subset of $\R^d$, we define $\nabla^2\psi$ along a geodesic $t\mapsto \gamma(t)$, for almost every $t$ in an interval $I$ of non-zero length, as $\nabla^2\psi(\dot\gamma(t),\dot\gamma(t)):=(\psi\circ\gamma)''(t)$, in the sense of sigma-additive (nonnegative) measures on~$I$. Thus, we are able to compare $\nabla^2\psi^1$ and $\nabla^2\psi^2$ for any two $\Gamma$-convex functions $\psi^1$ and $\psi^2$. In particular, the inequality $\nabla^2\psi \ge c I$ means that $(\psi\circ\gamma)''(t)\ge c\|\dot\gamma(t)\|^2$ for a.e. $t$, for any geodesic $\gamma$.
Finally, for any $C^2$ function $g$, we define $\nabla^2(g\circ\psi):=(g'\circ\psi)\nabla^2\psi+(g''\circ\psi)d\psi\otimes d\psi$ along the geodesics.

\begin{Lemma}
\label{extension}
Assume that, for every point $o\in \bar{\cD}$, there exist {an open subset $U_o$ of $\R^d$ containing $o$ and a nonnegative special $\Gamma$-convex function $\psi|_{U_o\cap\R^d}$}, vanishing only at $o$ and satisfying $\nabla^2\psi_o\ge cI$ for some $c=c(o)>0$. Then, for any $o\in \bar{\cD}$, any open $U\subset\R^d$ containing $o$, and any $C^2$ function~$\psi$ on $\bar{\cD}$, such that $\psi|_{U \cap \bar{\cD}}$ is a special $\Gamma$-convex function, there exists an open subset $U'$ of $\R^d$ containing $o$ and a global special $\Gamma$-convex function $\bar\psi$ on $\bar{\cD}$ coinciding with $\psi$ on $\bar \cD \cap U'$. 
\end{Lemma}

\proof

First consider the case $o\in \cD$. There exists $\varepsilon>0$ such that $B(o,\varepsilon)\subset U\cap \cD$. Consider a smooth nonincreasing  function $\eta : \R_+\to \R$ such that $\eta(r)=1$ if $r\le \varepsilon/2$ and $\eta(r)=0$ if $r\ge \varepsilon$. Let $\alpha>0$ such that $\psi_o\ge \alpha$ outside $B(o,\varepsilon/2)$ ($\alpha$ exists by compactness of $\bar\cD\backslash B(o,\varepsilon/2)$). Let now $g : \R_+\to \R$ be  a smooth nondecreasing convex function satisfying $g(r)=0$ for $r\le \alpha/2$ and $g(r)=r^p$ for some $p\ge 3$ for $r\ge \alpha$.
 Define for $M>0$ the function 
\begin{align*}
\bar\psi : &\bar{\cD}\to \R\\
x&\mapsto \eta(|x-o|)\psi(x)+ M(g\circ\psi_o)(x).
\end{align*}
We have $\bar\psi=\psi$ on $U':=\{\psi_o<\alpha/2\}$ (which is a neighborhood of $o$ included in $U\cap\cD$) since it is included in $B(o,\varepsilon/2))$.

Outside $B(o,\varepsilon)$ we have $\bar\psi=M(g\circ\psi_o)$ which is special $\Gamma$-convex since 
\begin{equation}\label{nabla}
\nabla(g\circ\psi_o)=(g'\circ \psi_o)\nabla\psi_o
\end{equation}
and
\begin{equation}\label{Hess}
 \nabla^2(g\circ\psi_o)=(g'\circ\psi_o) \nabla^2\psi_o+(g''\circ\psi_o)d\psi_o\otimes d\psi_o.
\end{equation}

Finally, letting $\eta_o(x)=\eta(|x-o|)$, we have outside U'
\begin{equation}\label{Hess2}
\begin{split}
 \nabla^2\bar\psi\\=&
M(g'\circ\psi_o) \nabla^2\psi_o+M(g''\circ\psi_o)d\psi_o\otimes d\psi_o
\\&+\eta_o \nabla^2\psi+\psi \nabla^2\eta_o+ d\eta_o\otimes d\psi +d\psi\otimes d\eta_o\\
&\ge \left(Mp(\alpha/2)^{p-1}c-\left\|\psi  \nabla^2\eta_o+d\eta_o\otimes d\psi +d\psi\otimes d\eta_o\right\|_\infty\right)I.
\end{split}
\end{equation}
Choosing 
\begin{equation}\label{M}
M\ge \frac{2^{p-1}\left\|\psi  \nabla^2\eta_o+d\eta_o\otimes d\psi +d\psi\otimes d\eta_o\right\|_\infty}{p\alpha^{p-1}c}
\end{equation}
we get the result for $o\in \cD$.

Consider now the case $o\in\partial\cD$. Without loss of generality we can assume that $\psi\ge 0$. The proof will be very similar, after we have constructed a cut-off  function $\eta_o$ whose gradient has positive scalar product with all outward normal vectors: a special cut-off function. We can assume that $U\cap \bar{\cD}$ is $C^2$-diffeomorphic to a convex set in $\R^d$. Let us call $\phi$ the diffeomorphism.  
Restricting again $U$, we can assume that $\phi(\bar U\cap \partial{\cD})$ is the graph of a convex function $f: \bar B_{d-1}(0,r)\to [0,m']$ with $B_{d-1}(0,r)$ the Euclidean ball in $\R^{d-1}$, $o=0_{\R^d}=(0,f(0))$, $0$ is a minimum for $f$, and $m'=\max\{f(x), \ |x|=r\}$. Possibly restricting $U$ and changing $\phi$, we can also assume that $f$ vanishes only at $0$, that $\phi(\bar U)= \bar B_{d-1}(0,r)\times [0,m']$  and that $\phi(\bar U\cap \bar{\cD})=\{(x,y)\in \bar B_{d-1}(0,r)\times [0,m'],\ y\ge f(x)\}$. Consequently, letting 
$m:=\min\{f(x), \ |x|=r/2\},$ we have by convexity of $f$ : $2m\le \min\{f(x), \ |x|=r\}\le m'$. Define the cut-off function
\begin{align*}
\tilde\eta_o :\R^{d-1}\times \R&\to \R\\
(x,y)&\mapsto \eta(y)
\end{align*}
with $\eta :\R\to \R$ smooth and nonincreasing, $\eta(y)=1$ if $y\le m$, $\eta(y)=0$ if $y\ge 2m$. 

For $(x,f(x))\in \partial{\cD}$, any element of the normal unitary exterior cone $\mathfrak{n}_u((x,f(x)))$ has the form
$\displaystyle
v=\frac1{\sqrt{\alpha^2+1}}(\alpha, -1)
$
with $\alpha\in \partial f(x)$ the subdifferential of $f$ at $x$. Then $$ v\cdot \nabla \tilde\eta_o((x,f(x)))=\frac{(-1)\eta'(f(x))}{\sqrt{\alpha^2+1}}\ge 0$$ implying that $\tilde\eta_o$ is a special cut-off function. 

For obtaining a special cut-off function in the original $\bar\cD$, just compose with the $C^2$ diffeomorphism $\phi$. We let $\eta_o =\tilde\eta_o\circ \phi$.

The rest of the proof is similar to the first part. The positive real number $\alpha>0$ is now defined such that $\psi_o\ge \alpha$ outside $\phi^{-1}\left(B_{d-1}(0,r/2)\times [0,m)\right)\cap\bar\cD$ and again, $U'=\{x\in \bar\cD, \ \psi_o(x)<\alpha/2\}$. The function $\bar \psi$  will be defined for $x\in\bar\cD$ as 
$$
\bar\psi(x)=\eta_o(x)\psi(x)+M(g\circ\psi_o)(x).
$$
The outward gradient property for $\bar\psi$ is directly obtained from the formula
$$
\nabla\bar\psi=\eta_o\nabla\psi+ \psi\nabla\eta_o+M(g'\circ \psi_o)\nabla\psi_o
$$
and the fact that we have assumed that $\psi\ge0$, together with $\eta_o\ge 0$ and $g'\circ \psi_o\ge 0$. For the positive Hessian property, the proof is similar to the first part (with $B(o,\varepsilon)$ replaced by $\phi^{-1}\left(B_{d-1}(0,r)\times [0,2m)\right)\cap\bar\cD$)  and the details are left to the reader.
\eproof

\smallskip

As a corollary of Lemma~\ref{extension}, we obtain a global characterization of $\Gamma$-martingales with drifts via (global or not) special $\Gamma$-convex functions.

\begin{Corollary}
\label{global gamma mart char}
Under the assumptions of Lemma~\ref{extension}, let $X$ be a {continuous adapted process} with values in $\bar{\cD}$ and let $f$ be an element of $\HP{1}$. Then, $X$ is a $\Gamma$-martingale with drift $f$ if and only if for all special $\Gamma$-convex functions {(or, equivalently, for all global special $\Gamma$-convex functions)} $\psi$ on $\bar\cD$, the {finite-variational component} of the real semimartingale $\psi(X_t)+\int_0^t \nabla \psi(X_s)\cdot f_s \,\ud s$ is nondecreasing.
\end{Corollary}
\proof
The proof is analogous to the one of Proposition~\ref{prop caracterisation gamma martingale}.
\eproof


\section{The case of dimension two}

\subsection{Properties of \texorpdfstring{$\bar{\mathcal{D}}$}{D}}

The case $d=2$ has an important property that is crucial for the main results of this work: namely, if $\bar{\cD}\subset\R^2$ is simply connected, then it is a $CAT(0)$ geodesic space, which means, roughly speaking, that triangles in $\bar{\cD}$ are thinner than in $\mathbb{R}^2$.

\smallskip

Most results in the subsequent part of the paper rely on the following assumption (which, nevertheless, is cited explicitly in each formal statement). 
\begin{Assumption}\label{ass:3}
$d=2$ and $\cD$ is simply connected.
\end{Assumption}

We begin by recalling the definitions and several properties of $CAT(0)$ spaces. The following definition comes from Section 2.1 in \cite{Alexander-Kapovitch-Petrunin-19}.
\begin{Definition}
    \label{def:CAT(0)}
    Consider a metric space $(\mathcal{X},d_{\mathcal{X}})$. Then, $\mathcal{X}$ is a $CAT(0)$ space, a.k.a. Alexandrov non-positively curved space (shortly, NPC space), if, for all $(x,y,p,q) \in \mathcal{X}^4$ that admit minimizing geodesics $\gamma^{x,y}$, $\gamma^{x,p}$, $\gamma^{p,y}$, $\gamma^{x,q}$ and $\gamma^{q,y}$, we have, for all $z \in \gamma^{x,y}_{[0,1]}$,
\begin{align*} 
d_{\mathcal{X}}(p,q) \leqslant  d_{\mathbb{R}^2}(\tilde{p},\tilde{z}) + d_{\mathbb{R}^2}(\tilde{q},\tilde{z}),
\end{align*}
where $(\tilde{x},\tilde{y},\tilde{p})$ (resp., $(\tilde{x},\tilde{y},\tilde{q})$) is a triangle in $\mathbb{R}^2$ whose edges have Euclidean lengths equal to $d_{\mathcal{X}}(x,y)$, $d_{\mathcal{X}}(x,p)$, $d_{\mathcal{X}}(y,p)$ (resp., $d_{\mathcal{X}}(x,y)$, $d_{\mathcal{X}}(x,q)$, $d_{\mathcal{X}}(y,q)$), and $\tilde{z} \in [\tilde{x},\tilde{y}]$ is such that $d_{\mathbb{R}^2}(\tilde{x},\tilde{z})=d_{\mathcal{X}}(x,z)$. 
\end{Definition}

In the case of geodesic spaces,
we have a simpler and more intuitive characterization of $CAT(0)$ spaces, see Section 2.2.2 in \cite{Alexander-Kapovitch-Petrunin-19}.

\begin{Proposition}
    \label{prop:def-alternative:CAT(0)}
    A geodesic space $\mathcal{X}$ is a $CAT(0)$ space if and only if all its
triangles are thin, i.e., for all $(x,y,p) \in \mathcal{X}$ and $z \in \gamma^{x,y}_{[0,1]}$, we have $d_{\mathcal{X}}(p,z) \leqslant d_{\mathbb{R}^2}(\tilde{p},\tilde{z})$, {where $(\tilde{x},\tilde{y},\tilde{p})$ is a triangle in $\mathbb{R}^2$ whose edges have Euclidean lengths equal to $d_{\mathcal{X}}(x,y)$, $d_{\mathcal{X}}(x,p)$, $d_{\mathcal{X}}(y,p)$, and $\tilde{z} \in [\tilde{x},\tilde{y}]$ is such that $d_{\mathbb{R}^2}(\tilde{x},\tilde{z})=d_{\mathcal{X}}(x,z)$.}
\end{Proposition}



\begin{Theorem}\label{thm:CAT(0)}
    Let Assumption \ref{ass:3} hold. Then, $\bar{\cD}$ is a geodesic $CAT(0)$ space: i.e., it is a geodesic space and, equipped with its geodesic metric $d_{\bar{\cD}}$, it is a CAT(0) space.  
\end{Theorem}

We refer to Theorem 4.4.1 in \cite{Alexander-Kapovitch-Petrunin-19} (see also \cite{Bishop-08}) for the complete proof of Theorem \ref{thm:CAT(0)}. 
Nevertheless, for a pedagogical purpose, we sketch the proof of the theorem.

\smallskip

\noindent \textbf{Sketch of Proof of Theorem \ref{thm:CAT(0)}}

From the previous Section, we already know that $\bar{\cD}$ is a geodesic space. Thus, we only need to prove that it is $CAT(0)$.
Let us consider a triangle with vertices $(A,B,C) \in \bar{\cD}^3$ and edges $\gamma^{A,B}$, $\gamma^{B,C}$ and $\gamma^{C,A}$. If the triangle is flat, the result is obvious, so we assume in the following that the triangle is not flat. Let us define $A' = \gamma^{A,B}_{t_{A,B}} = \gamma^{C,A}_{t_{A,C}}$ with 
$$t_{A,B} = \inf \{ t \in [0,1] | \gamma^{A,B}_t \neq \gamma^{C,A}_{t'}, \quad \forall t' \in [0,1] \}$$
and
$$t_{A,C} = \sup \{ t \in [0,1] | \gamma^{C,A}_t \neq \gamma^{A,B}_{t'}, \quad \forall t' \in [0,1] \}.$$
By the same token, we define $B'$, $C'$, $t_{B,A}$, $t_{B,C}$, $t_{C,A}$ and $t_{C,B}$.
By the uniqueness of minimal geodesics, we get that $\gamma^{A,B}_{[0,t_{A,B}]}$ and $\gamma^{C,A}_{[t_{A,C},1]}$ (resp. $\gamma^{B,C}_{[0,t_{B,C}]}$ and $\gamma^{A,B}_{[t_{B,A},1]}$, $\gamma^{C,A}_{[0,t_{C,A}]}$ and $\gamma^{B,C}_{[t_{C,B},1]}$) coincide, and the concatenation of $\gamma^{A',B'}$, $\gamma^{B',C'}$ and $\gamma^{C',A'}$ (i.e. the triangle with vertices $(A',B',C')$) is a Jordan curve. Then, the interior of this Jordan curve is necessarily in $\cD$ and this Jordan curve can reach the boundary of $\bar{\cD}$ only from the exterior. In particular, it implies that, if we follow one edge of the triangle $(A',B',C')$, then we can only turn to the exterior direction: indeed, if we turn to the interior direction, then we can find an alternative straight line shortcut which is a contradiction with the fact that edges are minimal geodesics. Finally, this Jordan curve is a triangle with concave edges and it is possible to show that this is a thin triangle.
\eproof

\smallskip



\smallskip

Many interesting properties of $\bar\cD$ are now intrinsically inherited from its $CAT(0)$ nature. We collect them in the next proposition, and we refer to Propositions 2.2.3 and 2.2.7 in \cite{Alexander-Kapovitch-Petrunin-19} for the proof.
\begin{Proposition}
    \label{Prop:CAT(0):cD}
    {Let Assumption \ref{ass:3} hold.} Then, we have:
    \begin{enumerate}[i)]
        \item Minimizing geodesics are unique, and $\gamma^{x,y}$ depends continuously on $(x,y)$ in the sense of the uniform topology, i.e., with respect to
        $$d(\gamma^{x,y},\gamma^{x',y'}) := \sup_{t \in [0,1]} |\gamma^{x,y}_t - \gamma^{x',y'}_t|.$$
        \item The geodesic distance is a $\Gamma$-convex function on $\bar{\cD} \times \bar{\cD}$, i.e., for all minimizing geodesics $\gamma^1$ and $\gamma^2$,
        $t \mapsto d_{\bar{\cD}}(\gamma^1_t,\gamma^2_t)$ is a (real) convex function.
        \item Any geodesic is a minimizing geodesic. 
    \end{enumerate}
\end{Proposition}

\begin{Remark}
    $ $
    \begin{itemize}
        \item The assumption that $\bar\cD$ is simply connected is necessary for the $CAT(0)$ property and for Proposition \ref{Prop:CAT(0):cD} to hold. Indeed, consider a domain $\cD$ with a circular hole in it and choose any three distinct points on this circle. The resulting triangle is clearly not thin. Moreover, one can easily deduce that the properties of Proposition \ref{Prop:CAT(0):cD} are not satisfied in this counter-example.
        \item The assumption $d=2$ is also crucial. Indeed, in higher dimensions, $\bar{\cD}$ can only remain a $CAT(0)$ space under very strong additional assumptions (see Theorem 4.3.1. and Proposition 4.2.6 in \cite{Alexander-Kapovitch-Petrunin-19}). Moreover, in higher dimensions, it is easy to construct a smooth domain $\cD$ such that the uniqueness of minimal geodesics does not hold: for example, consider a domain whose boundary contains a hemisphere (see the end of Section 5 in \cite{Chassagneux-Nadtochiy-Richou-22}).
    \end{itemize}
\end{Remark}

 
{

We now introduce some notations that are needed in order to derive a version of It\^o's formula tailored to our setting.
\begin{Definition}\label{def:re de rot-angle}  
    Let $\mathfrak{D}$ be an open, bounded and connected domain in a Euclidean space.
Assume that, for any two points in $\bar{\mathfrak{D}}$, there exists a unique minimizing geodesic between them. Then, for any $x,y \in \bar{\mathfrak{D}}$, we denote by $\overrightarrow{xy}$ the vector $\dot{\gamma}_0$, where $\gamma:[0,1] \rightarrow \bar{\mathfrak{D}}$ is the minimizing geodesic between $x$ and $y$.

If $\mathfrak{D}\subset\R^2$, then, for any $(x,y) \in \bar{\mathfrak{D}}\times \bar{\mathfrak{D}}$ with $x\neq y$, we denote by $\theta(x,y) \in (-\pi,\pi]$ the (unique) angle such that
    \begin{align}\label{eq e rot-angle}
        -\overrightarrow{yx} = R(\theta(x,y))\, \overrightarrow{xy},
    \end{align}
    where $R(\theta)$ denotes the rotation matrix of angle $\theta$. By convention, we set $R(\theta(x,x))=I$.

\end{Definition}
}

\medskip

{
For later use, we define $\Psi : \bar{\cD} \times \bar{\cD} \rightarrow \mathbb{R}^+$ via 
\begin{align}\label{eq.Psi.d2.def}
\Psi(x,y) = d_{\bar{\cD}}^2(x,y), 
\end{align}
where $ d_{\bar{\cD}}$ is the geodesic distance in $\bar{\cD}$.
}

\smallskip

By combining Assumptions \ref{ass:main} and \ref{ass:3}, together with Proposition \ref{prop:extsphere:interiorcone}, we deduce the additional properties of geodesics in $\bar{\cD}$.

\begin{Proposition}
    \label{prop:geodesics}
    {Let Assumption \ref{ass:3} hold and let $\gamma$ be a geodesic in $\bar\cD$.} Then, there exists a finite $N$ and a partition $0=t_1 <\cdots<t_{2N}=1$ such that one of the following two statements holds:
    \begin{itemize}
        \item For all integers $k$ equal to $1$ modulo $4$, the curve $\gamma$ does not turn to the left on $[t_k,t_{k+1}]$.
        \item For all even integers $k$, the curve $\gamma$ is a straight line on $[t_k,t_{k+1}]$.
        \item For all integers $k$ equal to $3$ modulo $4$, the curve $\gamma$ does not turn to the right on $[t_k,t_{k+1}]$.
    \end{itemize}
    or
    \begin{itemize}
        \item For all integers $k$ equal to $1$ modulo $4$, the curve $\gamma$ does not turn to the right on $[t_k,t_{k+1}]$.
        \item For all even integers $k$, the curve $\gamma$ is a straight line on $[t_k,t_{k+1}]$.
        \item For all integers $k$ equal to $3$ modulo $4$, the curve $\gamma$ does not turn to the left on $[t_k,t_{k+1}]$.
    \end{itemize}
\end{Proposition}

\proof 
Using Assumption \ref{ass:main} and the fact that $\cD$ is simply connected, we deduce that $\partial D$ is a Jordan curve for which we can find an orientation. Moreover, apart from the initial and terminal points, at any point $\gamma_t$ where the geodesic curve touches the boundary, it is tangent to the boundary and, consequently, $\dot \gamma_t$ either points in the same direction as the orientation (in this case, we denote $c(\gamma_t):=1$) or in the opposite direction (in this case, we denote $c(\gamma_t):=-1$). In the other cases, we denote $c(\gamma_t):=0$.

1. We first assume that $\gamma_0 \in \cD$ (or $\gamma_1 \in \cD$ and we reverse the time). Then we define 
$$r_1 := \inf \{t \in (0,1), c(\gamma_t) \neq 0 \} \quad \text{and} \quad r_k := \inf \{t \in (0,1), c(\gamma_t) \notin \{c(\gamma_{r_{k-1}}),0\}\}, \quad \forall k \geqslant 1$$
with the convention $\inf \emptyset = +\infty$. We can show that there is no finite accumulative point for the sequence $(r_k)_{k \geqslant 1}$. Indeed, if we have an accumulative point it means that $r_k \rightarrow \tilde{t}$ and so $\gamma_{r_k} \rightarrow \gamma_{\tilde{t}} \in \partial \cD$. Moreover, $(\gamma_{r_{2k}})_{k \in \mathbb{N}^*}$ converges to $\gamma_{\tilde{t}}$ from one side while $(\gamma_{r_{2k+1}})_{k \in \mathbb{N}}$ converges from the other side. In particular, we necessarily have $\tilde{t}=1$. Moreover, if we reverse time from $1$, we easily see that the only allowed direction from $\gamma_1$ to get a geodesic is $-\dot \gamma_1$ which is in contradiction with the interior cone property. So we get a finite sequence $(r_k)_{k \leqslant N}$ with $N \geqslant 0$. If $N=0$, it is sufficient to set $t_1=0$ and $t_2=1$. If not, we can set
 $t_{2k-1} = r_{k}$ for all $1 < k \leqslant K$, $t_{2k} := \sup \{t \in [r_k,r_{k+1}], c(\gamma_t) \neq 0\}$, for all $1 \geqslant k < N$ and $t_1=0$, $t_{2N}=1$. We can remark that this sequence is stricly increasing, that is to say: $t_{2k} <r_{k+1}$ for all $1 \geqslant k < N$ since $\partial \cD$ is a Jordan curve.

 2. If there exists $t \in (0,1)$ such that $\gamma_t \in \cD$, then we can do the same reasoning for $\gamma_{[0,t]}$ and $\gamma_{[t,1]}$ and then concatenate sequences obtained which show the result. 
 
 3. Finally, the last possibility is the case where $\gamma$ lives in $\partial \cD$. By continuity of $\gamma$ and the Jordan curve, we easily get that $\gamma$ always turn to the same side: if not, one side allows some straight line shortcuts which is in contradiction with the fact that $\gamma$ is a geodesic.
%
\eproof

\subsection{Properties of \texorpdfstring{$\Psi$}{Psi}, \texorpdfstring{$\theta$}{theta} and It\^o formula}

\begin{Proposition}
\label{Prop-equivalence-Psi-Euclidean}
{Let Assumption \ref{ass:3} hold.}
Then, the function $\Psi$ defined in \eqref{eq.Psi.d2.def} is continuous (with respect to the Euclidean topology), and there exists a constant $C \geqslant 1$ such that 
\begin{equation}
\label{equivalence normes}
 |x-y|^2 \leqslant \Psi(x,y) \leqslant C|x-y|^2, \quad \forall x,y, \in \bar{\cD}.
 \end{equation}
\end{Proposition}

\proof
The left-hand side of the inequality is obvious. so we focus on the right-hand side and the continuity.
Let us take $(x,y) \in \bar{\cD}^2$ such that $|x-y|\leqslant R_0/2$. Then the line segment $[x,y] \subset \cD_{R_0/2}$. By using Lemma \ref{lem:proj:lip}, we get that the projection of $[x,y]$ onto $\bar{\cD}$ is a path in $\bar{\cD}$ between $x$ and $y$, with length bounded by $2|x-y|$: in other words, 
$$\Psi(x,y) \leqslant 4 |x-y|^2,$$
which gives us the continuity of $\Psi$ in a neighborhood of the diagonal set $\{(x,x) | x \in \bar{\cD}\}$. We easily get the continuity of $\Psi$ in $\bar{\cD}$ by recalling that 
 $$|\Psi^{1/2}(x,y)-\Psi^{1/2}(x',y')|\leqslant \Psi^{1/2}(x,x') + \Psi^{1/2}(y,y').$$ 
 To conclude, we just remark that $(x,y) \mapsto \Psi(x,y) |x-y|^{-2}$ is a continuous function on the compact set $\bar{\cD} \cap \{(x,y)| |x-y| \geqslant R_0/2\}$ and then it is bounded.
\eproof

\begin{Proposition}
    \label{Prop-Rotation}  
    {Let Assumption \ref{ass:3} hold.}
    Then, for any geodesic $\gamma$, its velocity function $\dot \gamma$ {(which is well-defined for every $t$)} is $\frac{|\gamma_0-\gamma_1|^2}{R_0}$-Lipschitz and hence absolutely continuous on $[0,1]$, and $|\ddot\gamma| \leqslant \frac{|\gamma_0-\gamma_1|^2}{R_0}$ where $R_0$ is defined in \eqref{ineq:extsphere}.
    
    In addition, for any $(x,y) \in \bar{\cD}\times \bar{\cD}$, we have $|\theta(x,y)| \leqslant \frac{|x-y|}{R_0}$, and there exists $C>0$ such that, for all $(x,y) \in \bar{\cD}\times \bar{\cD}$,
     \begin{equation}
     \label{ineq-Rtheta}
     \| I- R(\theta(x,y))\| \leqslant C|x-y|
     \end{equation}
     and
     \begin{equation}
      \label{expansion-log2}
      \gamma_t^{x,y} = \gamma_s^{x,y} +\overrightarrow{xy}(t-s) + O(|y-x|^2), \quad \forall \, 0 \leqslant s \leqslant t \leqslant 1,
     \end{equation}
    where $O()$ is uniform in $(s,t)$. In particular, the above estimate implies that 
     \begin{equation}
     \label{expansion-log}
     \overrightarrow{xy}=y-x + O(|y-x|^2),
     \end{equation}
     {uniformly in $(x,y)\in \bar\cD\times\bar\cD$.}
     
     Finally, the following Taylor expansion holds:
     \begin{equation}
        \label{expansion-order2}
        \gamma_t^{x,y} = x + \overrightarrow{xy}t + \langle \gamma_t^{x,y}-x,v \rangle v + O(t^3),
     \end{equation}
     where $v$ is a unit vector orthogonal to $\overrightarrow{xy}$.
\end{Proposition}
    
\proof
    Let us consider a geodesic $\gamma$. Thanks to Proposition \ref{prop:geodesics}, $\gamma$ is locally a graph of a convex function. 
    Let us consider an interval $I \subset [0,1]$ such that $\gamma$ is the graph of a convex function on $I$ with finite derivative. We take $t<t'$ such that $(t,t') \in I^2$. If $\dot \gamma_t=\dot \gamma_{t'}$ then trivially 
    $|\dot \gamma_t-\dot \gamma_{t'}|\leqslant \frac{|\gamma_{0}-\gamma_{1}|^2}{R_0}|t-t'|$.
    If not, we consider 
    $$r := \inf\{s\geqslant t, \dot \gamma_s \neq \dot \gamma_t\}, \quad r' := \sup\{s \leqslant t', \dot \gamma_s \neq \dot \gamma_t\}$$
    which means that $\gamma_r \in \partial \cD$ and $\gamma_{r'}\in \partial \cD$. We denote $(u_r,v_r)$ (resp. $(u_{r'},v_{r'})$) a direct orthonormal basis such that $\dot \gamma_r$ (resp. $\dot \gamma_{r'})$ is positively colinear to $u_r$ (resp. $u_{r'})$. Since $\gamma_r$ is not colinear to $\gamma_{r'}$, there is a unique point $M$ at the intersection between the line orthogonal to $u_r$ containing $\gamma_{r}$ and the line orthogonal to $u_{r'}$ containing $\gamma_{r'}$. 
    We easily compute that 
    $$M = \gamma_r + \frac{\langle \gamma_{r'}-\gamma_r,u_{r'} \rangle}{\langle v_r,u_{r'}\rangle} v_r.$$
    Since we have the exterior sphere property (with radius $R_0$), then we must have $|M-\gamma_r| \geqslant R_0$ which gives us
    $$|\langle v_r,u_{r'}\rangle| \leqslant \frac{1}{R_0} |\langle \gamma_{r'}-\gamma_r,u_{r'} \rangle| \leqslant \frac{|\gamma_0-\gamma_1|}{R_0}|r'-r| \leqslant \frac{|\gamma_0-\gamma_1|}{R_0}|t'-t|.$$
    By using this inequality, we get
    \begin{align*}
        |\dot \gamma_t -\dot \gamma_{t'}|^2 &= |\dot \gamma_r -\dot \gamma_{r'}|^2 = 2|\gamma_0-\gamma_1|^2 (1-\langle u_r,u_{r'} \rangle)\\
        &= 2|\gamma_0-\gamma_1|^2 \left(1-\sqrt{1-|\langle v_r,u_{r'}\rangle|^2}\right)\\
        &\leqslant 2|\gamma_0-\gamma_1|^2 \left(1-\sqrt{1-\frac{|\gamma_0-\gamma_1|}{R_0}|t'-t|}\right) 
    \end{align*}
    as soon as $\frac{|\gamma_0-\gamma_1|}{R_0}|t'-t| \leqslant 1$ which is fulfilled if we take $I$ small enough. Since the right hand side of the previous inequality is asymptotically equivalent to $\frac{|\gamma_0-\gamma_1|^4 }{R_0^2}|t'-t|^2$ when $|t'-t|$ tends to $0$, we finally get that $\dot \gamma$ is Lipschitz with the constant Lipschitz $\frac{|\gamma_0-\gamma_1|^2 }{R_0}$.

    By the same computations, we also have, for $(t,t') \in I^2$,
    $$|\sin(\theta(\gamma_t,\gamma_{t'}))|= |\sin(\theta(\gamma_r,\gamma_{r'}))| = |\langle v_r,u_{r'}\rangle| \leqslant \frac{|\gamma_0-\gamma_1|}{R_0}|t'-t|$$
    which gives us $|\theta(\gamma_0,\gamma_1)| \leqslant \frac{|\gamma_0-\gamma_1|}{R_0}$ and \eqref{ineq-Rtheta} by using Proposition \ref{Prop-equivalence-Psi-Euclidean} and the boundedness of $\|I-R(\theta(x,y))\|$.

    Finally, for all $(x,y) \in \bar{\cD}\times \bar{\cD}$ and $0 \leqslant s \leqslant t \leqslant 1$, we have 
    $$\gamma^{x,y}_t = \gamma^{x,y}_s + \int_s^t \dot \gamma^{x,y}_u du = x + \overrightarrow{xy}(t-s) + \int_s^t \dot \gamma^{x,y}_u - \dot \gamma^{x,y}_0 du$$
    and, using The Lipschitz property of $\dot \gamma$ and Proposition \ref{Prop-equivalence-Psi-Euclidean}, 
    \begin{eqnarray*}
    \left|\int_s^t \dot \gamma^{x,y}_u - \dot \gamma^{x,y}_0 \ud u\right| \ &\leqslant& \int_t^s |u|\frac{|x-y|^2}{R_0} \ud u \leqslant C|x-y|^2 
    \end{eqnarray*}
    which gives us \eqref{expansion-log} and \eqref{expansion-log2}.
    Moreover, we have 
    \begin{align*}
    \gamma_t^{x,y} - x - \overrightarrow{xy}t - \langle \gamma_t^{x,y}-x,v \rangle v &= \langle \gamma_t^{x,y} - x - \overrightarrow{xy}t,\frac{\overrightarrow{xy}}{|\overrightarrow{xy}|}\rangle \frac{\overrightarrow{xy}}{|\overrightarrow{xy}|}\\
    &= \int_0^t \langle \dot\gamma_u^{x,y} -\dot\gamma_0^{x,y},\frac{\overrightarrow{xy}}{|\overrightarrow{xy}|}\rangle \ud u \frac{\overrightarrow{xy}}{|\overrightarrow{xy}|}\\
    &= \int_0^t \left(\cos(\theta(\gamma_0^{x,y},\gamma_u^{x,y}))-1\right) |x-y| \ud u \frac{\overrightarrow{xy}}{|\overrightarrow{xy}|}
    \end{align*}
    which gives us, by using the Lipschitz property of $\theta$,
    \begin{align*}
        |\gamma_t^{x,y} - x - \overrightarrow{xy}t - \langle y-x,v \rangle v| &\leqslant |x-y| \int_0^t \frac{|\theta(\gamma^{x,y}_u,\gamma^{x,y}_0)|^2}{2}\ud u \\
        &\leqslant |x-y|^3 \int_0^t \frac{|u|^2}{2R_0^2}\ud u  = \frac{|x-y|^3}{6R_0^2}t^3.
    \end{align*}
    \eproof

    



\medskip

\begin{Proposition}
    \label{Prop-Psi}
    {Let Assumption \ref{ass:3} hold.}
    {Then, the function $\Psi$ defined in \eqref{eq.Psi.d2.def} is an element of $C^1(\bar\cD)$, and $\nabla \Psi(x,y) = -2(\overrightarrow{xy},\overrightarrow{yx})$, where $\overrightarrow{xy}$ is given in Definition \ref{def:re de rot-angle}.}
    In addition, $\Psi$ is a special $\Gamma$-convex function vanishing precisely on the diagonal
    $\Delta = \{(x,x) | x \in \bar{\cD}\}$: in particular, for all $(x,y)\in\bar\cD\times\bar\cD$, $(u,v) \in \mathfrak{n}(x) \times \mathfrak{n}(y)$, we have $\langle \nabla \Psi(x,y),(u,v) \rangle \geqslant 0$. Moreover, $\Psi$ is strongly $\Gamma$-convex, in the sense that, for all geodesics $\gamma^1$ and $\gamma^2$, and for all $t\in [0,1]$, we have
    \begin{align}
    \label{ineq:psi:prop}
    \Psi(\gamma_t^1,\gamma_t^2) \geqslant& \Psi(\gamma_0^1,\gamma_0^2)+ \langle \nabla \Psi(\gamma^1_0,\gamma^2_0), (\dot{\gamma}^{1}_0,\dot{\gamma}^{2}_0) \rangle t \\
    \nonumber &+ 2\int_0^t (t-s) \left| \dot{\gamma}^1_s -R(\theta(\gamma_{s}^1,\gamma_{s}^2)) \dot{\gamma}^2_s \right|^2 \ud s.
    \end{align}
\end{Proposition}


\proof
Proposition \ref{Prop:CAT(0):cD} already gives us the $\Gamma$-convexity of $\Psi$. To prove the other statements of the proposition, we start by proving the following directional second-order expansion of $\Psi$: for all $(x,y)\in \bar{\cD}^2$, for all geodesics $\gamma^1$ and $\gamma^2$ starting from $x$ and $y$, we have 
\begin{align}
    \label{ineq:psi:final3} \Psi(\gamma_{\varepsilon}^1,\gamma_{\varepsilon}^{2})  \geqslant & \Psi(x,y)  - \varepsilon 2\langle \overrightarrow{xy}, \dot{\gamma}^{1}_0 \rangle - \varepsilon 2\langle \overrightarrow{yx},\dot{\gamma}^{2}_0 \rangle  +\varepsilon^2 |\dot{\gamma}^{1}_0- R(\theta(x,y))\dot{\gamma}^{2}_0|^2
    +  o(\varepsilon^2)
\end{align}

{\bf Step 1}. Let us consider $(x,y) \in \cD$ and $(u,\tilde{u})$ some vectors such that $|u|=|\tilde{u}| \leqslant r$ and $\bar{B}(x,r) \cup \bar{B}(y,r) \subset \cD$. In particular, $\gamma^{x,x+ u}$ and $\gamma^{y,y+\tilde{u}}$ are some straight lines. 

{\bf Step 1.a}. We first study the case where $\gamma^{x,y}$ is a straight line. Euclidean geometry gives us directly, for all $\varepsilon \in [0,1]$,
\begin{align*}
& \Psi(\gamma_{\varepsilon}^{x,x+ u},\gamma_{\varepsilon}^{y,y+ \tilde{u}})\\
 \geqslant& |x+\varepsilon u-y-\varepsilon\tilde{u}|^2 = |x-y|^2 + \varepsilon 2\langle x-y,u-\tilde{u} \rangle + \varepsilon^2|u-\tilde{u} |^2\\
=& \Psi(x,y) - 2\varepsilon \langle \overrightarrow{xy}, \dot{\gamma}^{x,x+u}_0 \rangle - 2\langle \overrightarrow{yx},\dot{\gamma}^{y,y+\tilde{u}}_0 \rangle + \epsilon^2|\dot{\gamma}^{x,x+ u}_0- R(\theta(x,y))\dot{\gamma}^{y,y+\tilde{u}}_0|^2
\end{align*}
since $R(\theta(x,y))=I$.




{\bf Step 1.b}. It remains now to study the case where $\gamma^{x,y}$ is not a straight line which implies in particular that this geodesic touches the boundary. 

We denote $I_x^{\varepsilon, u}$ (resp. $I_y^{ \varepsilon, \tilde{u}}$) the confluence point of geodesics $\gamma^{x,y}$ and $\gamma^{\gamma_{\varepsilon}^{x,x+u},y}$ (resp. $\gamma^{x,y}$ and $\gamma^{x,\gamma_{\varepsilon}^{y,y+\tilde{u}}}$), i.e. the unique point where these geodesics meet:  indeed, there exists a unique pair $(t,t') \in [0,1]^2$ such that $I_x^{u} = \gamma^{x,y}_t = \gamma^{\gamma_{\varepsilon}^{x,x+u},y}_{t'}$ and, for all $\eta >0$, there exists $t'' \in [(t-\eta) \vee 0,t)$ (if the latter interval is not empty) such that $\gamma^{x,y}_{t''} \notin \gamma^{\gamma_{\varepsilon}^{x,x+u},y}_{[t',1]}$ or there exists $t'' \in (t',(t'+\eta) \wedge 1]$ (if the latter interval is not empty) such that $\gamma^{\gamma_{\varepsilon}^{x,x+u},y}_{t''} \notin \gamma^{x,y}_{[0,t]}$. 
The following lemma will be proved after.
\begin{Lemma}
    \label{lem:confluence points}
    {Let Assumption \ref{ass:3} hold, and let $(x,y)$ be an arbitrary element of $\bar{\cD}\times\bar\cD$.} Assume that $\gamma^{x,y}$ is not a straight line.
    Then, there exists $\eta>0$ such that, for all $u,\tilde{u}$ satisfying $|u|,|\tilde{u}| \leqslant \eta$, $x+u \in \bar{\cD}$, $x+\tilde{u} \in \bar{\cD}$, and for all $\varepsilon \in [0,1]$, there exist $t_{x,\varepsilon, u} \leqslant t_{y,\varepsilon,\tilde{u}}$ such that $I_x^{\varepsilon, u} =\gamma^{x,y}_{t_{x,\varepsilon, u}}$ and $I_y^{\varepsilon, \tilde{u}} = \gamma^{x,y}_{t_{y,\varepsilon,\tilde{u}}}$. In other words, $\gamma^{x,y}$ and $\gamma^{\gamma_{\varepsilon}^{x,x+u},\gamma_{\varepsilon}^{y,y+\tilde{u}}}$ coincide between $I_x^{\varepsilon, u}$ and $I_y^{\varepsilon, \tilde{u}}$.
\end{Lemma}

Thus, we assume that $r \leqslant\eta$, with $\eta$ given by Lemma \ref{lem:confluence points}. We denote by $(u_x,v_x)$ an orthonormal basis such that $u_x$ is positively colinear with $\overrightarrow{xy}$, and denote by $(u_y,v_y)$ an orthonormal basis such that $u_y$ is positively colinear with $-\overrightarrow{yx}$. We use following notations: $\gamma_{\varepsilon}^{x,x+u}-x= \alpha u_x+\beta v_x$ and $\gamma_{\varepsilon}^{y,y+\tilde{u}}-y = \delta u_x+\rho v_x$.

Then, Lemma \ref{lem:confluence points} gives us
\begin{align*}
    \Psi(\gamma_{\varepsilon}^{x,x+u},\gamma_{\varepsilon}^{y,y+\tilde{u}}) = \left(\Psi^{1/2}(\gamma_{\varepsilon}^{x,x+u},I_x^{\varepsilon, u})+\Psi^{1/2}(I_x^{\varepsilon, u},I_y^{\varepsilon,\tilde{u}})+ \Psi^{1/2}(I_y^{\varepsilon,\tilde{u}},\gamma_{\varepsilon}^{y,y+\tilde{u}})\right)^2.
\end{align*}
We can easily show that $x+ \alpha u_x$ is the projection of $\gamma_{\varepsilon}^{x,x+u}$ on $\gamma^{x+\alpha u_x,I_x^u}_{[0,1]}$. Indeed, $\gamma_{\varepsilon}^{x,x+u} -(x+\alpha u_x)$ is orthogonal to $u_x$, $\gamma^{x,x+\alpha u_x}$ is a straight line and $\gamma^{x+ \alpha u_x,I_x^{\varepsilon,u}}$ is ``behind'' the line passing through $x+\alpha u_x$ along the direction $u_x$. 
Then, the Pythagorean inequality in $CAT(0)$ spaces (see, e.g., Theorem 2.3.3 in \cite{Jost-97}) yields 
\begin{align}
    \label{Pythagorean-inequality}
\Psi(\gamma_{\varepsilon}^{x,x+u},I_x^{\varepsilon, u}) \geqslant \Psi(x+\alpha u_x,I_x^{\varepsilon, u}) + \beta^2.
\end{align}
Applying the same reasoning on the $y$-side of the geodesic $\gamma^{x,y}$, we obtain
\begin{align*}
    & \Psi(\gamma_{\varepsilon}^{x,x+u},\gamma_{\varepsilon}^{y,y+\tilde{u}})\\
     \geqslant &\left(\sqrt{\Psi(x+\alpha u_x,I_x^{\varepsilon, u})+\beta^2}+\Psi^{1/2}(I_x^{\varepsilon, u},I_y^{\varepsilon, \tilde{u}})+\sqrt{\Psi(I_y^{\varepsilon, \tilde{u}},y+ \delta u_y)+\rho^2}\right)^2 \\
    =& \left(\sqrt{(\Psi^{1/2}(x,I_x^{\varepsilon, u})-\alpha)^2+\beta^2}+\Psi^{1/2}(I_x^{\varepsilon, u},I_y^{\varepsilon,\tilde{u}})+\sqrt{(\Psi^{1/2}(I_y^{\varepsilon, \tilde{u}},y)+ \delta)^2+\rho^2}\right)^2.
\end{align*}
Thus, a series expansion of the right hand side of the previous inequality gives us
\begin{align}
    \nonumber &\Psi(\gamma_{\varepsilon}^{x,x+u},\gamma_{\varepsilon}^{y,y+\tilde{u}})\\ 
    \nonumber \geqslant & \Psi(x,y) + \alpha^2+\delta^2+2\Psi^{1/2}(x,y)(\delta-\alpha)\\
    \nonumber & +\Psi^{1/2}(x,y)\left( \frac{\beta^2}{\Psi^{1/2}(x,I_x^{\varepsilon,u})} + \frac{\rho^2}{\Psi^{1/2}(I_y^{\varepsilon,\tilde{u}},y)}\right) + o(\alpha^2+\beta^2+\delta^2+\rho^2)\\
    \nonumber \geqslant & \Psi(x,y) + \alpha^2+\delta^2+2\Psi^{1/2}(x,y)(\delta-\alpha)\\
    \nonumber & +\left( 1+ \frac{\Psi^{1/2}(I_y^{\varepsilon,\tilde{u}},y)}{\Psi^{1/2}(x,I_x^{\varepsilon,u})}\right)\beta^2 + \left(1+\frac{\Psi^{1/2}(x,I_x^{\varepsilon,u})}{\Psi^{1/2}(I_y^{\varepsilon,\tilde{u}},y)}\right)\rho^2 + o(\alpha^2+\beta^2+\delta^2+\rho^2)\\
    \label{ineq:Psi}
    \geqslant & \Psi(x,y) + 2\Psi^{1/2}(x,y)(\delta-\alpha)  +(\delta-\alpha)^2+ (\beta-\rho)^2 + o(\alpha^2+\beta^2+\delta^2+\rho^2)
\end{align}
where have used that $\Psi^{1/2}(x,y) \geqslant \Psi^{1/2}(x,I_x^{\varepsilon,u})+ \Psi^{1/2}(I_y^{\varepsilon,\tilde{u}},y)$ and Young inequality. 
By replacing $\alpha$, $\beta$, $\delta$ and $\rho$ by their values
\begin{align*}
\alpha = \varepsilon \langle u,u_x \rangle, \quad \beta = \varepsilon \langle u,v_x \rangle, \quad \delta = \varepsilon \langle \tilde{u},u_y \rangle, \quad \rho = \varepsilon \langle \tilde{u},v_y \rangle,
\end{align*}
we finally obtain
\begin{align}
    \label{ineq:psi:final} \Psi(\gamma_{\varepsilon}^{x,x+u},\gamma_{\varepsilon}^{y,y+\tilde{u}})  \geqslant & \Psi(x,y)  - \varepsilon 2\langle \overrightarrow{xy}, \dot{\gamma}^{x,x+ u}_0 \rangle - \varepsilon 2\langle \overrightarrow{yx},\dot{\gamma}^{y,y+ \tilde{u}}_0 \rangle  \\&+\varepsilon^2 |\dot{\gamma}^{x,x+ u}_0- R(\theta(x,y))\dot{\gamma}^{y,y+\tilde{u}}_0|^2
    \nonumber
    +  o(\varepsilon^2),
\end{align}
which is the desired result.

\smallskip

{\bf Step 2}. Now we need to consider the cases where $x$ or $y$ belong to $\partial \cD$. Let us assume that $x \in \partial \cD$ and $y \notin \partial \cD$, other cases will be treated in a similar way. 


{\bf Step 2.a}. We start by considering the case where there exists $\eta>0$ such that $\gamma^{x,y}$ is a straight line on $[0,\eta]$ with $\eta>0$.

{\bf Step 2.a.i}. If, for $r=|u|>0$ small enough, $x+ u \in \bar{\cD}$ and $\gamma^{x,x+ u}$ is also a straight line, then we can use same arguments as in the case 1, up to a slight detail: now $x+ \alpha u_x$ can leave $\bar \cD$, 
but since this point is just used in an intermediate calculation step, we can artificially define the geodesic $\gamma^{x+\alpha u_x,y}$ as a straight line extension of the geodesic $\gamma^{x,y}$ which gives us $\Psi^{1/2}(x+\alpha u_x,I_x^{\varepsilon,u}) = \Psi^{1/2}(x,I_x^{\varepsilon,u})- \alpha$.
Then, \eqref{ineq:psi:final} strays true. 

{\bf Step 2.a.ii}. Otherwise the perturbing geodesic $\gamma^{x,x+u}$, denoted $\gamma$ to lighten the notation, is such that $\gamma_{|[0,\eta]}$ is not a straight line for any $\eta>0$. It means, in particular, that $\dot \gamma_0$ is tangent to $\partial \cD$ in $x$. In this case, all computations of Step 1 hold true until inequality \eqref{ineq:Psi}. We now need to analyse more carefully $\alpha, \beta, \delta$ and $\gamma$. In this case we have
$$\alpha =  \langle \gamma_{\varepsilon}-x,u_x \rangle, \quad \beta =  \langle \gamma_{\varepsilon}-x,v_x \rangle, \quad \delta = \varepsilon \langle \tilde{u},u_y \rangle, \quad \rho = \varepsilon \langle \tilde{u},v_y \rangle.$$
Then we use \eqref{expansion-log}-\eqref{expansion-order2} to get 
$$\alpha =  \epsilon \langle \dot \gamma_0,u_x \rangle + \langle \gamma_{\varepsilon}-x,v\rangle \langle v,u_x \rangle + o(\varepsilon^2),\quad \beta = \varepsilon\langle \dot \gamma_0,v_x\rangle+o(\varepsilon)$$
where $v$ is orthogonal to $\dot \gamma_0$. Since $\dot \gamma_0$ is tangent to $\partial \cD$, we can take $v$ in the normal exterior cone. Then, $\langle v,u_x \rangle \leqslant 0$. Moreover, when we follow the curved geodesic starting from $x$, we must turn such that $\langle \gamma_{\varepsilon}-x,v\rangle \geqslant 0$ for $\varepsilon>0$ small enough.
we finally get
\begin{align}
    \label{ineq:psi:final2} \Psi(\gamma_{\varepsilon},\gamma_{\varepsilon}^{y,y+\tilde{u}})  \geqslant & \Psi(x,y)  - \varepsilon 2\langle \overrightarrow{xy}, \dot{\gamma}_0 \rangle - \varepsilon 2\langle \overrightarrow{yx},\dot{\gamma}^{y,y+ \tilde{u}}_0 \rangle  \\&+\varepsilon^2 |\dot{\gamma}_0- R(\theta(x,y))\dot{\gamma}^{y,y+\tilde{u}}_0|^2
    \nonumber
    +  o(\varepsilon^2)
\end{align}
which is, once again, the desired result.

{\bf Step 2.b}. it remains to addressed the case where $x \in \partial\cD$ and $\gamma^{x,y}_{[0,\eta]}$ is not a straight line for any $\eta>0$. 

{\bf Step 2.b.i}. If, for $r=|u|>0$ small enough, $x+ u \in \bar{\cD}$ and $\gamma^{x,x+ u}$ is also a straight line, then we can use same arguments as in the case 2.a.i., up to a slight detail: if $\alpha>0$, then $x+\alpha u_x \notin \gamma^{x,y}$ even for $\epsilon$ small enough. Nevertheless, the Pythagorean inequality \eqref{Pythagorean-inequality} stays true and the triangular inequality gives us 
$$\Psi^{1/2}(x+\alpha u_x,I_x^{\varepsilon,u})+\alpha \geqslant \Psi^{1/2}(x,I_x^{\varepsilon,u}).$$
Then, \eqref{ineq:psi:final} is still valid. 

{\bf Step 2.b.ii}. Otherwise, the perturbing geodesic, still denoted $\gamma$, is such that $\gamma_{|[0,\eta]}$ is not a straight line for any $\eta>0$. This case is obvious since we have only two possibilities: 
\begin{itemize}
    \item If $\dot \gamma_{0}$ is positively colinear to $\overrightarrow{xy}$, then $\gamma_{\varepsilon} \in \gamma^{x,y}$ for $\varepsilon>0$ small enough and then $\Psi^{1/2}(\gamma_{\varepsilon},I_x^{\varepsilon,u})= \Psi^{1/2}(x,I_x^{\varepsilon,u})-\epsilon \Psi^{1/2}(\gamma_0,\gamma_1)$.
    \item Otherwise $\dot \gamma_{0}$ is negatively colinear to $\overrightarrow{xy}$. In this case, $\gamma$ is just a continuation of $\gamma^{x,y}$ and then $\Psi^{1/2}(\gamma_{\varepsilon},I_x^{\varepsilon,u})= \Psi^{1/2}(x,I_x^{\varepsilon,u})+\epsilon \Psi^{1/2}(\gamma_0,\gamma_1)$.
\end{itemize}
Finally, \eqref{ineq:psi:final2} is still valid.

\smallskip

Thus, we have finished the proof of the directional second-order expansion of $\Psi$ \eqref{ineq:psi:final3}. Since $\Psi$ is $\Gamma$-convex and $(x,y)\mapsto -2(\overrightarrow{xy},\overrightarrow{yx})$ is continuous on $\bar{\cD}\times\bar{\cD}$, we deduce that $\Psi$ is $C^1$, with $\nabla \Psi(x,y) = -2 (\overrightarrow{xy},\overrightarrow{yx})$. Moreover, $\mathfrak{n}(x) \neq \{0\}$ iff $x \in \partial \cD$ and, hence, $\langle \overrightarrow{xy},u\rangle \leqslant 0$ for all $u \in \mathfrak{n}(x)$: it implies that $\Psi$ is a special $\Gamma$-convex function.

It only remains to prove \eqref{ineq:psi:prop}. Let us denote $f(t) = \Psi(\gamma^1_t,\gamma^2_t)$. Since $f$ is convex, $f''$ exists as a positive Radon measure, and $f'$ is a.e. differentiable in $[0,1]$. Thanks to \eqref{ineq:psi:final3}, we have 
$$f''(t) \geqslant 2|\dot{\gamma}^{1}_t- R(\theta(\gamma^1_t,\gamma_t^2))\dot{\gamma}^{2}_t|^2, \quad \text{a.e.}$$
and 
$$f'(t) \geqslant f'(0)+\int_0^t 2|\dot{\gamma}^{1}_s- R(\theta(\gamma^1_s,\gamma_s^2))\dot{\gamma}^{2}_s|^2 \ud s,\quad \forall t \in [0,1].$$
Then, the function
$$t \mapsto \Psi(\gamma_t^1,\gamma_t^2) - \Psi(\gamma_0^1,\gamma_0^2)- \langle \nabla \Psi(\gamma^1_0,\gamma^2_0), (\dot{\gamma}^{1}_0,\dot{\gamma}^{2}_0) \rangle t -  2 \int_0^t (t-s) \left| \dot{\gamma}^1_s -R(\theta(\gamma_{s}^1,\gamma_{s}^2)) \dot{\gamma}^2_s \right|^2 \ud s$$
is a convex function, with zero value and vanishing derivative at $t=0$, hence, it is nonnegative in $[0,1]$, which yields \eqref{ineq:psi:prop}.

\eproof

{{\noindent \bf Proof of Lemma \ref{lem:confluence points}. }}
Thanks to Proposition \ref{prop:geodesics}, we can find $\varepsilon_1>0$ and $\varepsilon_2>0$ such that $t \mapsto \theta(x,\gamma^{x,y}_t)$ is a non-zero monotonic function on $[0,\varepsilon_1]$ and $t \mapsto \theta(\gamma^{x,y}_t,y)$ is a non-zero monotonic function on $[1-\varepsilon_2,1]$. Thanks to Proposition \ref{Prop-Rotation}, these functions are also continuous. We set 
$$t_1 := \inf\{t \in [0,\varepsilon_1]|\theta(x,\gamma_{t_1}) = \frac{1}{3}\theta(x,\gamma_{\varepsilon_1})\}$$ and
$$t_2 := \sup\{t \in [1-\varepsilon_2,1]|\theta(\gamma_{t_2},y) = \frac{1}{3}\theta(\gamma_{\varepsilon_2},y)\}.$$ 
Then, we necessarily have $0<t_1<t_2<1$ and $(\gamma_{t_1},\gamma_{t_2}) \in \partial \cD^2$. Let us consider the ray starting from $\gamma_{t_1}$, in the 
direction negatively colinear with $\dot \gamma_{t_1}$, and let us denote by $x'$ the first intersection of this ray with $\partial \cD$ (excepted $\gamma_{t_1}$). 
Then $[\gamma_{t_1},x']$ slices $\cD$ into two sub-domains $\cD^1$ and $\cD^2$. Since $x \notin [\gamma_{t_1},x']$, we can assume that $x \in \cD^1$ without loss of generality. Since the geodesic triangle with vertices $(z,x',\gamma_{t_1})$ is thin, for any $z \in \bar{\cD^1}$, we must have that $\dot\gamma^{z,\gamma_{t_1}}$ is positively colinear with $\dot \gamma_{t_1}$. In particular, if we glue $\gamma^{z,\gamma_{t_1}}$ and $\gamma^{x,y}_{[t_1,1]}$, we obtain a geodesic, which is necessarily the minimal geodesic $\gamma^{z,y}$ by Proposition \ref{Prop:CAT(0):cD}. Then, we only need to notice that there exists $\eta_1>0$ small enough such that $B(x,\eta_1) \cap \cD = B(x,\eta_1) \cap \cD^1$: for any point $z \in B(x,\eta_1)$, $\gamma^{z,y}$ contains $\gamma_{t_1}$. By the same token, we can find $\eta_2>0$ small enough such that, for any point $z \in B(y,\eta_2)$, the curve $\gamma^{x,z}$ contains $\gamma_{t_2}$, which is sufficient to conclude.
\eproof

\smallskip

A useful corollary of Proposition \ref{Prop-Psi} is given bellow, and it says that $\Gamma$ martingales on $\bar \cD$ are characterized by the global special $\Gamma$-convex functions (recall Remark \ref{rem:global.vs.local.GammaConvex}.3).

{
\begin{Corollary}
\label{cor-global gamma convex function enough}
The assumptions of Lemma \ref{extension} are satisfied under Assumption \ref{ass:3}. In particular, the conclusion of Corollary \ref{global gamma mart char} holds under Assumption \ref{ass:3}.
\end{Corollary}
}

\proof 
For any point $o\in \bar{\cD}$, we define $\psi_o := \Psi(o,.)$. According to Proposition \ref{Prop-Psi}, $\psi_o$ is a  nonnegative special $\Gamma$-convex function on $\bar \cD$ vanishing only at $o$. Moreover, computations in the proof of Proposition $\ref{Prop-Psi}$, in particular $\eqref{ineq:psi:final3}$, gives us that $\nabla^2\psi_o\ge 2I$.
\eproof

\medskip

In the next section, we apply It\^o's formula to function $\Psi$, which is only $C^1$ and defined on a closed set. In order to justify this, we establish a tailor-made It\^o's formula for $C^1$ $\Gamma$-convex functions on $\bar{\cD}$, or on $\bar{\cD}\times \bar{\cD}$. 

\begin{Proposition}
\label{Prop-TaylormadeIto-1}
Let {Assumption \ref{ass:3} hold.}
 Consider a $C^1$ $\Gamma$-convex function $\psi : \bar{\cD} \rightarrow \R$ (resp., $\psi : \bar{\cD} \times \bar{\cD}\rightarrow \R$) and a (Euclidean) semimartingale $(X_t)_{t \in [0,T]}$ with values in $\bar{\cD}$ (resp., $\bar{\cD}\times \bar{\cD}$). Then,
 \begin{align*} 
 \psi(X_t) \geqslant \psi(X_s) +\int_s^t \nabla \psi (X_u) \ud X_u, \quad \forall \, 0 \leqslant s \leqslant t \leqslant T.
 \end{align*}
 More generally, consider any continuous function $S: \bar{\cD} \rightarrow \mathcal{S}^+$ (resp. $S: \bar{\cD}\times \bar{\cD} \rightarrow \mathcal{S}^+$), where $\mathcal{S}^+ \subset \R^{2 \times 2}$ (resp., $\mathcal{S}^+ \subset \R^{4 \times 4}$) is the set of all symmetric positive semidefinite matrices, such that the following holds for all geodesics $\gamma$:
 \begin{equation}
 \label{E1}
    \psi(\gamma_t) \geqslant \psi(\gamma_0)+\langle \nabla \psi(\gamma_0),\dot \gamma_0 \rangle t + \int_0^t (t-s) \dot \gamma_s^\top S(\gamma_{s})\dot \gamma_s \ud s.
 \end{equation}
 Then, we have
 $$\psi(X_t) \geqslant \psi(X_s) +\int_s^t \nabla \psi (X_u) \ud X_u + \frac12 \int_s^t  \langle \ud X_u, S(X_u) \ud X_u \rangle, \quad \forall 0 \leqslant s \leqslant t \leqslant T.$$
\end{Proposition}

\proof
We prove the result only for the first case, i.e. $\psi : \bar{\cD} \rightarrow \R$. The second case follows the same lines.

Let us remark that the first part of the proposition is implied by the second one when we take $S=0$. So let us prove the second part. For this we need to improve inequality~\eqref{E1}. Observe that this inequality implies that for any $C^2$ curve $\varphi_t$ taking values inside $\cD$, we have 
$$
\psi(\varphi_t) \geqslant \psi(\varphi_0)+\langle \nabla \psi(\varphi_0),\dot \varphi_0 \rangle t + \int_0^t (t-s) \left(\dot \varphi_s^\top S(\varphi_{s})\dot \varphi_s+\langle \nabla \psi(\varphi_s),\ddot\varphi_s\rangle\right) \ud s.
$$
Indeed, this can be done by approximating $\varphi_t$  by piecewise affine geodesics inside $\cD$, and   passing to the limit thanks to the continuity of $S$ and $\nabla \psi$. 

Then, any $\Gamma$-geodesic $\gamma$ in $\bar\cD$ can be approximated uniformly up to order~$1$ by $C^2$ curves $\varphi_t^n$ inside $\cD$ such that  $\ddot\varphi^n$ are uniformly bounded and a.e. converges to $\ddot\gamma$. This implies that
$$
\int_0^t (t-s) \langle \nabla \psi(\varphi_s^n),\ddot\varphi_s^n\rangle \ud s \to \int_0^t (t-s) \langle \nabla \psi(\gamma_s),\ddot\gamma_s\rangle \ud s\quad \hbox{as $n\to\infty$}.
$$
Finally~\eqref{E1} together with the continuity of $S$ imply that 
\begin{equation}\label{E2} 
\psi(\gamma_t) \geqslant \psi(\gamma_0)+\langle \nabla \psi(\gamma_0),\dot \gamma_0 \rangle t + \int_0^t (t-s) \left(\dot \gamma_s^\top S(\gamma_{s})\dot \gamma_s +\langle \nabla \psi(\gamma_s),\ddot\gamma_s\rangle\right) \ud s.
\end{equation}

 We consider a uniform grid $s=t_0^n \leqslant...\leqslant t^n_n=t$ with step size $h:=(t-s)/n$. Then for $0 \leqslant i <n$,
\begin{align*}
 \psi(X_{t_{i+1}^n}) =&  \psi(\gamma^{X_{t_i^n},X_{t_{i+1}^n}}_1)\\ \geqslant& \psi(\gamma^{X_{t_i^n},X_{t_{i+1}^n}}_0) + \langle \nabla \psi(\gamma^{X_{t_i^n},X_{t_{i+1}^n}}_0), \overrightarrow{X_{t_i^n}X_{t_{i+1}^n}} \rangle\\ & +\int_0^1 (1-s) (\dot\gamma^{X_{t_i^n},X_{t_{i+1}^n}}_{s})^\top S(\gamma^{X_{t_i^n},X_{t_{i+1}^n}}_{s})\dot\gamma^{X_{t_i^n},X_{t_{i+1}^n}}_{s} \ud s\\
&+\int_0^1 (1-s)\left\langle \nabla \psi(\gamma^{X_{t_i^n},X_{t_{i+1}^n}}_s), \ddot\gamma^{X_{t_i^n},X_{t_{i+1}^n}}_s\right\rangle\ud s.\\
\end{align*}
On the other hand, we have by Taylor formula with reminder
$$
X_{t_{i+1}^n}-X_{t_i^n}=\overrightarrow{X_{t_i^n}X_{t_{i+1}^n}}+\int_0^1(1-s)\ddot\gamma^{X_{t_i^n},X_{t_{i+1}^n}}_s \ud s.
$$
So the previous inequality can be transformed into
\begin{align*}
 \psi(X_{t_{i+1}^n}) =&  \psi(\gamma^{X_{t_i^n},X_{t_{i+1}^n}}_1)\\ \geqslant& \psi(\gamma^{X_{t_i^n},X_{t_{i+1}^n}}_0) + \langle \nabla \psi(\gamma^{X_{t_i^n},X_{t_{i+1}^n}}_0),
X_{t_{i+1}^n}-X_{t_i^n}
 \rangle\\ & +\int_0^1 (1-s) (\dot\gamma^{X_{t_i^n},X_{t_{i+1}^n}}_{s})^\top S(\gamma^{X_{t_i^n},X_{t_{i+1}^n}}_{s})\dot\gamma^{X_{t_i^n},X_{t_{i+1}^n}}_{s}
 \ud s.\\
&+\int_0^1 (1-s)\left\langle \nabla \psi(\gamma^{X_{t_i^n},X_{t_{i+1}^n}}_s)-\nabla \psi(\gamma^{X_{t_i^n},X_{t_{i+1}^n}}_0), \ddot\gamma^{X_{t_i^n},X_{t_{i+1}^n}}_s\right\rangle\ud s.\\
\end{align*}

Then, as usual, we sum previous inequality over $i$ and we pass to the limit in $n$ by using \eqref{expansion-log2}. Observe that the last term in the right converges to $0$ by continuity of $\nabla\psi$ and boundedness of $\ddot\gamma^{X_{t_i^n},X_{t_{i+1}^n}}$.
\eproof

\smallskip

By combining Propositions \ref{Prop-Psi} and \ref{Prop-TaylormadeIto-1}, we obtain the following It\^o's inequality for the function $\Psi$ defined in \eqref{eq.Psi.d2.def}.

\begin{Corollary}
    \label{Ito-formula-Psi}


Let {Assumption \ref{ass:3} hold} and consider two (Euclidean) semimartingales 
$$\ud X^1_t = \ud A^1_t + \beta^1_t \ud W_t\quad \text{and} \quad  \ud X^2_t = \ud A^2_t + \beta^2_t \ud W_t$$ 
with values in $\bar{\cD}$, where, for $i \in \{1,2\}$, $A^i_t$ and $\beta^i$ are progressively measurable processes such that, for all $t \geqslant 0$, {$\int_0^t \ud \mathrm{Var}_{s}(A)<+\infty$} a.s. and $\int_0^t |\beta_s^i|^2 \ud s<+\infty$ a.s.. Then, for all $0 \leqslant s \leqslant t$,  
\begin{align*}
    \Psi(X_t^1,X_t^2) \geqslant& \Psi(X_s^1,X_s^2) - 2\int_s^t \overrightarrow{X_u^1 X_u^2} \ud X_u^1 -  2\int_s^t \overrightarrow{X_u^2 X_u^1 } \ud X_u^2\\
    & + \int_s^t |\beta^1_u -R(\theta(X_u^1,X_u^2))\beta^2_u|^2 \ud u.
\end{align*}
\end{Corollary}

Moreover, we obtain a slight generalization of Corollary \ref{Ito-formula-Psi} which is used in the next section.

\begin{Proposition}
    \label{Ito-formula-Psi-bis}
Let {Assumption \ref{ass:3} hold} and consider two (Euclidean) semimartingales 
$$\ud X^1_t = \ud A^1_t + \beta^1_t \ud W_t\quad \text{and} \quad  \ud X^2_t = \ud A^2_t + \beta^2_t \ud W_t$$ 
with values in $\bar{\cD}$, where, for $i \in \{1,2\}$, $A^i$ and $\beta^i$ are progressively measurable processes such that, for all $t \geqslant 0$, {$\int_0^t \ud \mathrm{Var}_{s}(A)<+\infty$} a.s. and $\int_0^t |\beta_s^i|^2 \ud s<+\infty$ a.s.

Consider also a nonnegative, absolutely continuous and progressively measurable process $B$, such that {$\int_0^t \ud \mathrm{Var}_{s}(B)<+\infty$} a.s. for all $t \geqslant 0$.
Then, for all $0 \leqslant s \leqslant t$,  
\begin{align*}
    B_t\Psi(X_t^1,X_t^2) \geqslant& B_s\Psi(X_s^1,X_s^2) - 2\int_s^t B_u\overrightarrow{X_u^1 X_u^2} \ud X_u^1 -  2\int_s^t B_u\overrightarrow{X_u^2 X_u^1 } \ud X_u^2\\
    & +\int_s^t \Psi(X_u^1,X_u^2) \ud B_u + \int_s^t B_u|\beta^1_u -R(\theta(X_u^1,X_u^2))\beta^2_u|^2 \ud u.
\end{align*}
\end{Proposition}

\medskip

We end this section with a stability and uniqueness result for $\Gamma$-martingales with a prescribed drift and terminal value, in a general continuous filtration.

\begin{Theorem}[Stability and uniqueness for a general continuous filtration]
    \label{Thm-Uniqu1}
    We make the following assumptions:
    \begin{itemize}
        \item Assumption \ref{ass:3} holds.
        \item  We consider the generator  
    $f: \Omega \times [0,T] \times \bar{\cD} \rightarrow \mathbb{R}^d$ such that $t \mapsto f(t,y)$ is progressively measurable for all $y \in \bar{\cD}$.
        \item For all $C>0$, $\E [e^{C\int_0^T |f(s,0)| \ud s}]<+\infty$
        \item $f$ is a Lipschitz function with respect to $y$: there exists $C_{f,y} \ge 0$ such that, for all $t \in [0,T]$, for all $y_1,y_2 \in \bar{\cD}$
        \begin{equation*}
            |f(t,y_1)- f(t,y_2)| \leqslant C_{f,y} |y_1-y_2| \quad a.s.
        \end{equation*}
    \end{itemize}
    Let $(Y_t)_{t \in [0,T]}$ and $(Y'_t)_{t \in [0,T]}$ be two $\Gamma$-martingales in $\bar{\cD}$, respectively, with drifts $f(\cdot,Y_\cdot)$ and $f(\cdot,Y'_\cdot)$ and with terminal values $\xi$ and $\xi'$. Then, for all $p>1$, 
    \[\sup_{t \in [0,T]}\mathbb{E}[\Psi(Y_t,Y'_t)] \leqslant C_p \mathbb{E}[\Psi(\xi,\xi')^p]^{1/p},\] 
    and the above inequality also for $p=1$ if the random variable $\int_0^T |f(s,0)| \ud s$ is bounded. In particular, if $\xi=\xi'$ then $Y=Y'$.
\end{Theorem}
    
    
    \proof
    1. We have $Y_t = \xi + \int_t^T f(s,Y_s) \ud s -\int_t^T \ud K_s - \int_t^T \ud M_s$ for all $0\leqslant t \leqslant T$. By using Proposition \ref{Prop-TaylormadeIto-1} with Proposition 
    \ref{Prop-Psi} to the process $\Psi(0,Y)$, we get that
    \begin{align*}
        & \Psi(0,Y_0)+ \sum_{i=1}^d \int_0^T \ud \langle M^i\rangle_s\\ 
        \leqslant& \Psi(0,Y_T) - \int_0^T \nabla_2 \Psi(0,Y_s) \ud Y_s\\
         \leqslant&\Psi(0,Y_T) + \int_0^T \nabla_2 \Psi(0,Y_s) f(s,Y_s)\ud s - \int_0^T \nabla_2 \Psi(0,Y_s) \ud M_s.
    \end{align*}
    By the standard localisation procedure and by using the Lipschitz property of $f$, we obtain
    $$\E [|\langle M\rangle_T|] \leqslant C + C\E \left[ \int_0^T C_{f,y}|Y_s| + |f(s,0)| \ud s \right] \leqslant C + C\E \left[ e^{\int_0^T |f(s,0)| \ud s} \right]<+\infty,$$
    that is to say, $M$ is a $L^2$ martingale. By the same token, the local martingale part of $Y'$, denoted $M'$, is also a $L^2$ martingale. 

    2. Let us consider $\lambda>0$ a parameter that will be set after. By using a mere generalization of Proposition \ref{Prop-TaylormadeIto-1} with Proposition 
    \ref{Prop-Psi} to the process $\Gamma^{\lambda}_t \Psi(Y_t,Y'_t)$ where $\Gamma^{\lambda}_t = e^{\lambda t+ \lambda \int_0^t |f(s,0)| \ud s}$, we have    
    \begin{align} \nonumber
    & \Gamma^{\lambda}_t\Psi(Y_t,Y'_t)+ \lambda \int_t^T (1+ |f(s,0)|)\Gamma^{\lambda}_s \Psi(Y_s,Y'_s) \ud s \\ \nonumber
     \leqslant & \Gamma^{\lambda}_T\Psi(\xi,\xi') - \int_t^T \Gamma^{\lambda}_s \nabla_1 \Psi(Y_s,Y'_s) \ud Y_s -\int_t^T \Gamma^{\lambda}_s \nabla_2 \Psi(Y_s,Y'_s) \ud Y'_s\\ \nonumber
      \leqslant& \Gamma^{\lambda}_T \Psi(\xi,\xi')+2\int_t^T \Gamma^{\lambda}_s \overrightarrow{Y_s Y'_s} dY_s +2\int_t^T \Gamma^{\lambda}_s \overrightarrow{Y'_s Y_s} dY'_s\\ \nonumber
      \leqslant& \Gamma^{\lambda}_T \Psi(\xi,\xi')-2\int_t^T \Gamma^{\lambda}_s \overrightarrow{Y_s Y'_s} f(s,Y_s) \ud s +2\int_t^T \Gamma^{\lambda}_s \overrightarrow{Y_s Y'_s} \ud M_s \\ \nonumber
      & -2\int_t^T \Gamma^{\lambda}_s \overrightarrow{Y'_s Y_s} f(s,Y'_s) \ud s +2\int_t^T \Gamma^{\lambda}_s \overrightarrow{Y'_s Y_s} \ud M'_s \\ \nonumber
     \leqslant& \Gamma^{\lambda}_T\Psi(\xi,\xi')-2\int_t^T \Gamma^{\lambda}_s \overrightarrow{Y_s Y'_s}\left( f(s,Y_s)-f(s,Y'_s) \right) \ud s +2\int_t^T \Gamma^{\lambda}_s \overrightarrow{Y_s Y'_s} \ud M_s  \\
     & +2\int_t^T \Gamma^{\lambda}_s \overrightarrow{Y'_s Y_s} \ud M'_s+2\int_t^T \Gamma^{\lambda}_s (R(\theta(Y'_s,Y_s)-I)\overrightarrow{Y'_s Y_s} f(s,Y'_s)\ud s.
     \label{Ineg-proof-uniqueness}
    \end{align}
    BDG and H\"older inequalities give us that 
    \begin{align*}
        \E \left[\sup_{0 \leqslant t \leqslant T} \left| \int_0^t \Gamma^{\lambda}_s \overrightarrow{Y_s Y'_s} \ud M_s\right|\right] \leqslant &  C\E \left[(\Gamma^{\lambda}_T)^{1/2} | \langle M_T \rangle|^{1/2}\right] \leqslant C\E \left[\Gamma^{\lambda}_T\right]^{1/2}\E\left[ | \langle M_T \rangle|\right]^{1/2}<+\infty
    \end{align*}
    which implies that $\int_t^. \Gamma^{\lambda}_s \overrightarrow{Y_s Y'_s} \ud M_s$ is a uniformly integrable martingale. By same arguments, $\int_t^.\Gamma^{\lambda}_s \overrightarrow{Y'_s Y_s} \ud M'_s$ is also a uniformly integrable martingale.
    By taking the expectation in \eqref{Ineg-proof-uniqueness}, we can use Proposition \ref{Prop-equivalence-Psi-Euclidean}, inequality \eqref{ineq-Rtheta} and the Lipschitz property of $f$ to obtain
    \begin{align*}
    &\mathbb{E}[\Gamma^{\lambda}_t\Psi(Y_t,Y'_t)] + \lambda \E \left[\int_t^T (1+ |f(s,0)|)\Gamma^{\lambda}_s \Psi(Y_s,Y'_s) \ud s\right] \\  
    \leqslant& \mathbb{E}[\Gamma^{\lambda}_T\Psi(\xi,\xi')]+C\int_t^T \mathbb{E}[\Gamma^{\lambda}_s\Psi(Y_s,Y'_s)]ds \\
    &+C\int_t^T \mathbb{E}[\Gamma^{\lambda}_s\Psi(Y_s,Y'_s) (C_{f,y}|Y'_s|+f(s,0)|)] ds\\
    \leqslant& \mathbb{E}[\Gamma^{\lambda}_T\Psi(\xi,\xi')]+C \E \left[\int_t^T (1+ |f(s,0)|)\Gamma^{\lambda}_s \Psi(Y_s,Y'_s) \ud s\right]
    \end{align*}
    recalling that $\bar{\cD}$ and so $Y$ and $Y'$ are bounded.
    Then we just have to take $\lambda \geqslant C$ in order to get
    $$ \mathbb{E}[\Psi(Y_t,Y'_t)]  \leqslant \mathbb{E}[\Gamma^{\lambda}_t\Psi(Y_t,Y'_t)] \leqslant \mathbb{E}[\Gamma^{\lambda}_T\Psi(\xi,\xi')].$$
    If $\int_0^T |f(s,0)| \ud s$ is bounded, then $\Gamma^{\lambda}_T$ too and the result is proved for $p=1$. Otherwise, we just have to apply Hölder inequality to conclude.
    \eproof

\section{Existence and uniqueness of solutions to reflected BSDEs in simply-connected two-dimension domains, with a Brownian filtration}



In this section, we develop the desired existence and uniqueness result for solutions to reflected BSDEs in the form \eqref{eq reflected bsde}, assuming that the filtration is Brownian and that $\cD$ is a bounded and simply connected subset of $\R^2$. As mentioned in Remark \ref{rem link bsde}, this problem is essentially equivalent to the problem of existence and uniqueness of $\Gamma$-martingales with prescribed drifts and terminal values.

\begin{Assumption}
    \label{ass:brownian-case}
    \begin{enumerate}
    \item We assume that $(\mathcal{F}_t)_{t \geqslant 0}$ is the augmented natural filtration of a $d'$-dimensional Brownian motion $(W_t)_{t \geqslant 0}$.
    \item Assumption \ref{ass:3} holds.
    \item  We consider a generator  
$f: \Omega \times [0,T]\times \bar{\cD} \times \mathbb{R}^{d \times d'} \rightarrow \mathbb{R}^d$ such that $t \mapsto f(t,y,z)$ is progressively measurable for all $y \in \bar{\cD}$ and $z \in \mathbb{R}^{d \times d'}$.
    \item $|f(.,0,0)|^{1/2} \in \widebar{\mathscr{H}^{\infty}}^{\BP{2}}$ 
    \item $f$ is a Lipschitz function with respect to $y$ and $z$: i.e., there exist $C_{f,y},C_{f,z} \ge 0$ such that, for all $t \in [0,T]$, for all $y,y' \in \bar{\cD}$, for all $z,z' \in \mathbb{R}^{d \times d'}$,
    \begin{equation*}
        |f(t,y,z)- f(t,y',z')| \leqslant C_{f,y} |y-y'|+C_{f,z} |z-z'| \quad a.s.
    \end{equation*}
    \end{enumerate}
\end{Assumption}

\begin{Remark}\label{rem:JN.ineq}
    Thanks to \cite{Schachermayer-96}, we know that $|f(.,0,0)|^{1/2} \in \widebar{\mathscr{H}^{\infty}}^{\BP{2}}$ is equivalent to the fact that $\int_ 0^. |f(s,0,0)|^{1/2} \ud W_s$ is a BMO $\varepsilon$-sliceable martingale for all $\varepsilon>0$. Thanks to John-Nirenberg inequality (see e.g. \cite{Kazamaki-94}), it implies in particular that
    \begin{equation}
        \label{moment-expo-f(0)}
    \left| \sup_{t \in [0,T]}\mathbb{E}_t \left[ e^{\lambda \int_t^T |f(s,0,0)|\ud s}\right]\right|_{\infty}<+\infty, \quad \forall \lambda>0.
    \end{equation}
    Let us remark that this assumption is fulfilled for example when  $|f(.,0,0)|^{1/2+\eta} \in \BP{2}$ for $\eta>0$. 
\end{Remark}

\subsection{A priori estimate, stability and uniqueness}

\begin{Proposition}
    \label{Prop-aprioriZ}
    {Let Assumption \ref{ass:brownian-case} hold} and let $(Y,Z,K)$ be a solution to the reflected BSDE \eqref{eq reflected bsde}. Then, there exists a constant $C_{\eqref{eq Prop-aprioriZ}}>0$, only depending on $\cD$, $T$, $C_{f,y}$ and $C_{f,z}$, such that 
     \begin{equation}
        \label{eq Prop-aprioriZ}
        \|Z\|_{\BP{2}}^2 \leqslant C_{\eqref{eq Prop-aprioriZ}}(1+\||f(.,0,0)|^{1/2}\|_{\BP{2}}^2)<+\infty.
     \end{equation}   
\end{Proposition}

\proof 
The proof follows the same strategy as the first step in the proof of Theorem \ref{Thm-Uniqu1}. We apply Corollary \ref{Ito-formula-Psi} to the processes $X^1=Y$ and $X^2=0$: for all $t \in [0,T]$,  
\begin{align*}
 &\Psi(Y_t,0)+\int_t^T |Z_s|^2 \ud s \\
 \leqslant& \Psi(Y_T,0) +2 \int_t^T \overrightarrow{Y_s 0} dY_s\\
 \leqslant&\Psi(Y_T,0) - 2 \int_t^T\overrightarrow{Y_s 0} f(s,Y_s,Z_s)\ud s +2 \int_t^T \overrightarrow{Y_s 0} Z_s \ud W_s\\
 \leqslant&\Psi(Y_T,0) + 2 \int_t^T \Psi(Y_s,0) \left(C_{f,y}|Y_s|+C_{f,z}|Z_s|+|f(s,0,0)|\right)\ud s +2 \int_t^T \overrightarrow{Y_s 0} Z_s \ud W_s.
\end{align*}
By considering a localizing sequence $(\tau_n)_{n \in \mathbb{N}}$ of stopping times, taking the conditional expectation and using the boundedness of $Y$, the linear growth of $f$ and Young inequality, we get
\begin{align*}
  \mathbb{E}_{t} \left[ \int_{t \wedge \tau_n}^{\tau_n} |Z_s|^2 \ud s  \right] \leqslant& C + C\mathbb{E}_t\left[ \int_{t \wedge \tau_n}^{\tau_n} (1+|Z_s|+ |f(s,0,0)|)\ud s\right]\\
  & \leqslant C + \frac{1}{2}\mathbb{E}_t\left[ \int_{t \wedge \tau_n}^{\tau_n} |Z_s|^2 \ud s\right] +C\mathbb{E}_t\left[ \int_{t }^{T}  |f(s,0,0)|\ud s\right]
\end{align*}
which gives us the result by taking $n \rightarrow +\infty$.

\begin{Corollary}
    {Under Assumption \ref{ass:brownian-case},} for any solution $(Y,Z,K)$ of the reflected BSDE \eqref{eq reflected bsde}, the process $\int_ 0^. Z_s \ud W_s$ is a BMO martingale. In particular, we have
    \begin{equation}
    \label{energy-ineq-Z}
    \left| \sup_{t \in [0,T]}\mathbb{E}_t \left[ \left(\int_t^T |Z_s|^2\ud s \right)^{n}\right]\right|_{\infty}\leqslant n! \|Z\|_{\BP{2}}^{2n},\quad  \forall n \in \mathbb{N}.
    \end{equation}
    and
    \begin{equation}
        \label{moment-expo-f(t,y,z)}
    \left| \sup_{t \in [0,T]}\mathbb{E}_t \left[ e^{\lambda \int_t^T |f(s,Y_s,Z_s)| +|Z_s|\ud s}\right]\right|_{\infty}<+\infty, \quad \forall \lambda>0.
    \end{equation}
\end{Corollary}

\proof 
The first result is a direct consequence of Proposition \ref{Prop-aprioriZ}. \eqref{energy-ineq-Z} comes from the energy inequality for BMO martingales. It remains to prove \eqref{moment-expo-f(t,y,z)}. Thanks to the Lipschitz property of $f$, the boundedness $\bar{\cD}$, H\"older inequality and Young inequality, we get, for all $\varepsilon>0$, 
\begin{align*}
\left| \sup_{t \in [0,T]}\mathbb{E}_t \left[ e^{\lambda \int_t^T |f(s,Y_s,Z_s)|+|Z_s|\ud s}\right]\right|_{\infty}&\le C\left| \sup_{t \in [0,T]}\mathbb{E}_t \left[ e^{\lambda C\int_t^T |f(s,0,0)|+|Z_s|\ud s}\right]\right|_{\infty}\\
&\le  C\left| \sup_{t \in [0,T]}\mathbb{E}_t \left[ e^{\lambda C\int_t^T |f(s,0,0)|\ud s}\right]\right|_{\infty} \left| \sup_{t \in [0,T]}\mathbb{E}_t \left[ e^{\lambda C\int_t^T |Z_s|\ud s}\right]\right|_{\infty}\\
&\le  C_{\varepsilon}\left| \sup_{t \in [0,T]}\mathbb{E}_t \left[ e^{\lambda C\int_t^T |f(s,0,0)|\ud s}\right]\right|_{\infty} \left| \sup_{t \in [0,T]}\mathbb{E}_t \left[ e^{\varepsilon\int_t^T |Z_s|^2\ud s}\right]\right|_{\infty}.
\end{align*}
Then, by taking $\varepsilon$ small enough, \eqref{moment-expo-f(0)} and John-Nirenberg inequality (see e.g. \cite{Kazamaki-94}) give us the result.
\eproof

\begin{Proposition}
    \label{prop stability}
    {Let Assumption \ref{ass:brownian-case} hold} and let $(Y,Z,K)$ (resp. $(Y',Z',K')$) be a solution of the reflected BSDE \eqref{eq reflected bsde} with the terminal condition $\xi$ (resp. $\xi'$) and the generator $f$ (resp. $f'$). 
    We assume that
    \begin{itemize}
        \item Assumption \ref{ass:brownian-case} holds with $f'$ in place of $f$,
        \item there exists a progressively measurable process $(\alpha_t)_{t \in [0,T]}$ such that $1+\1_{C_{f,z} >0}|Z'_t|+|f'(t,Y'_t,Z'_t)| \le \alpha_t$ for all $t \in [0,T]$,
        \item for all $\lambda >0$,
        \begin{equation}
            \label{condition alpha}
            \esp{e^{\lambda \int_0^T \alpha_s\ud s}} <+\infty.
        \end{equation}
    \end{itemize}
     We set $\Gamma^{\lambda}_t = e^{\lambda \int_0^t \alpha_s\ud s}$ for all $\lambda>0$. Then, there exists $C>0$ such that, for any $\eta>0$, there exists a constant $\lambda_0 \ge 0$, only depending on $\eta$, $\bar{\cD}$, $T$, $C_{f,y}$ and $C_{f,z}$, such that, for all $\lambda \ge \lambda_0$, we have
    \begin{align}\label{eq strong estim one}
        &\esp{\sup_{t\in[0,T]}\Gamma^{\lambda}_t \Psi(Y_t,Y'_t)} +\esp{\int_0^T\Gamma^{\lambda}_t|Z_t - R(\theta(Y_t,Y'_t))Z'_t|^2 \ud t }\\ \le&  C\esp{\Gamma^{\lambda}_T \Psi(\xi,\xi')}
        + \eta\esp{\int_0^T \Gamma^{\lambda}_t |f(t,Y'_t,Z'_t)-f'(t,Y'_t,Z'_t)|^2 \ud t }     \nonumber  
    \end{align}
    and
    \begin{align}\label{eq cooked for Z}
        & \esp{\int_0^T|Z_t - Z'_t|^2 \ud t }\\
        \le& C\esp{\Gamma^{\lambda}_T \Psi(\xi,\xi')
        + \int_0^T \Gamma^{\lambda}_t |f(t,Y'_t,Z'_t)-f'(t,Y'_t,Z'_t)|^2 \ud t } \nonumber\\
        &+ \nonumber C(1+\||f'(.,0,0)|^{1/2}\|_{\BP{2}}^2)\esp{\Gamma^{\lambda}_T \Psi(\xi,\xi')
        + \int_0^T \Gamma^{\lambda}_t |f(t,Y'_t,Z'_t)-f'(t,Y'_t,Z'_t)|^2 \ud t }^\frac12.
    \end{align}
\end{Proposition}

\proof {\bf Step 1.a.} 
To streamline the notation, we denote $f_t := f(t,Y_t,Z_t)$, $f'_t := f'(t,Y'_t,Z'_t)$ and $\delta f'_t := f(t,Y'_t,Z'_t)-f'(t,Y'_t,Z'_t)$.
Using Proposition \ref{Ito-formula-Psi-bis}, we get, for $0 \le t \le T$,
\begin{align}
    \Gamma^\lambda_T \Psi(Y_T,Y'_T) \ge& \Gamma^\lambda_t \Psi(Y_t,Y'_t) + \lambda \int_t^T \alpha_s \Gamma^\lambda_s \Psi(Y_s,Y'_s) \ud s + 2 \int_t^T \Gamma^\lambda_s \left(\overrightarrow{Y_s Y'_s}f_s + \overrightarrow{Y'_s Y_s} f'_s\right)  \ud s  \nonumber
    \\
    &- 2 \int_t^T \Gamma^\lambda_s \left(\overrightarrow{Y_s Y'_s}Z_s + \overrightarrow{Y'_s Y_s}Z'_s\right)  \ud W_s +\int_t^T\Gamma^\lambda_s|Z_s - R(\theta(Y_s,Y'_s))Z'_s|^2 \ud s. \label{eq starting point one}
\end{align}
We also compute
\begin{align*}
    \overrightarrow{Y_s Y'_s}f_s + \overrightarrow{Y'_s Y_s} f'_s
    &= \overrightarrow{Y_s Y'_s}(f_s-f'_s) + f'_s(\overrightarrow{Y_s Y'_s}+\overrightarrow{Y'_s Y_s})
    \\
    &\ge -\Psi(Y_s,Y'_s)^{1/2}|f_s-f'_s| -|f'_s||I-R(\theta(Y_s,Y'_s))||\overrightarrow{Y_s Y'_s}| \ud s
\end{align*}
where for the last inequality, we use \eqref{eq e rot-angle}. 
Thanks to the Lipschitz property of $f$,  we have
\begin{align*}
    |f_s-f'_s| \le& C_{f,y} |Y_s -Y'_s| + C_{f,z}|Z_s -Z'_s| + |\delta f'_s|\\
    \le & C \Psi(Y_s,Y'_s)^{1/2}+C|Z_s-R(\theta(Y_s,Y'_s))Z'_s|+C_{f,z}|Z'_s||I-R(\theta(Y_s,Y'_s))|+ |\delta f'_s|.
\end{align*}

Then combining the two previous inequality with \eqref{equivalence normes}, \eqref{ineq-Rtheta} and Young inequality, we obtain, for all $\nu>0$,
\begin{align*}
    \overrightarrow{Y_s Y'_s}f_s + \overrightarrow{Y'_s Y_s} f'_s
    \ge -C_{\nu}(1+|f'_s|+1_{C_{f,z} \neq 0}|Z'_s|)\Psi(Y_s,Y'_s) - \nu |\delta f'_s|^2 - \frac14 |Z_s-R(\theta(Y_s,Y'_s))Z'_s|^2.
\end{align*}
Inserting the previous inequality back into \eqref{eq starting point one}, we obtain, for $\lambda$ large enough with respect to the constant $C_{\nu}$ appearing in the previous inequality,
\begin{align}
    \Gamma^{\lambda}_T \Psi(Y_T,Y'_T) \ge& \Gamma^\lambda_t \Psi(Y_t,Y'_t) - 2\nu \int_t^T \Gamma^{\lambda}_s |\delta f'_s|^2 \ud s \nonumber  
    - 2 \int_t^T \Gamma^{\lambda}_s \left(\overrightarrow{Y_s Y'_s}Z_s + \overrightarrow{Y'_s Y_s}Z'_s\right)  \ud W_s\\ & +\frac12 \int_t^T\Gamma^{\lambda}_s|Z_s - R(\theta(Y_s,Y'_s))Z'_s|^2 \ud s.  \label{eq intrem point one}
\end{align}
By using Burkholder-Davis-Gundy inequality and H\"older inequality we compute
\begin{align}
    \nonumber &\esp{ \sup_{t \in [0,T]}
        \left|\int_0^t \Gamma^{\lambda}_s \left(\overrightarrow{Y_s Y'_s}Z_s + \overrightarrow{Y'_s Y_s}Z'_s\right)  \ud W_s \right|
    }\\
    \le& C
    \esp{ \left(\sup_{t \in [0,T]}\Gamma^{\lambda}_t \Psi(Y_t,Y'_t)\right)^{\frac12}
    \left(
    \int_0^T \Gamma^{\lambda}_s
    |Z_s - R(\theta(Y_s,Y'_s))Z'_s|^2
        \ud s
        \right)^\frac12 
    } \nonumber
    \\
    \le& C_{\eqref{eq for later one}}
     \esp{\sup_{t \in [0,T]}\Gamma^{\lambda}_t \Psi(Y_t,Y'_t)}^\frac12\esp{\int_0^T \Gamma^{\lambda}_s|Z_s - R(\theta(Y_s,Y'_s))Z'_s|^2\ud s}^\frac12
    \label{eq for later one}
    \\
    \le&
    C \esp{\Gamma^{2 \lambda}_T }^\frac12\esp{\Gamma^{2\lambda}_T}^\frac12\esp{\left(\int_0^T(|Z_s|^2+|Z'_s|^2)\ud s\right)^2}^\frac12 <+\infty\,, \nonumber
\end{align}
where for the last inequality, we used the boundedness of $\bar{\cD}$, H\"older inequality, \eqref{energy-ineq-Z}. This shows that the local martingale term in \eqref{eq intrem point one} is a true martingale. For later use, let us remark that $C_{\eqref{eq for later one}}$ only depends on the Burkholder-Davis-Gundy constant and the constant appearing in \eqref{ineq-Rtheta}.
We take expectation on both sides of \eqref{eq intrem point one} and get
\begin{align}\label{eq weak estim one}
    &\sup_{t \in [0,T]} \esp{\Gamma^{\lambda}_t\Psi(Y_t,Y'_t)} +\frac12 \esp{\int_0^T\Gamma^{\lambda}_t|Z_t - R(\theta(Y_t,Y'_t))Z'_t|^2 \ud t } \\
    \le& \esp{\Gamma^{\lambda}_T \Psi(\xi,\xi')
    + 2\nu \int_0^T \Gamma^{\lambda}_t |\delta f'_t|^2 \ud t }\;. \nonumber
\end{align}
{\bf Step 1.b.}  Now, from \eqref{eq intrem point one}, we deduce
\begin{align*}
    \Gamma^{\lambda}_T \Psi(Y_T,Y'_T) \ge& \sup_{t\in[0,T]}\Gamma^{\lambda}_t \Psi(Y_t,Y'_t)-2\nu \int_0^T \Gamma^{\lambda}_t |\delta f'_t|^2 \ud t
    \\
   & -4 \sup_{t \in [0,T]} |\int_0^t \Gamma^{\lambda}_s \left(\overrightarrow{Y_s Y'_s}Z_s + \overrightarrow{Y'_s Y_s}Z'_s\right) \ud W_s|.
\end{align*}
Using Burkholder-Davis-Gundy inequality and \eqref{eq for later one}, we obtain that 
\begin{align*}
    &\esp{\sup_{t\in[0,T]}\Gamma^{\lambda}_t \Psi(Y_t,Y'_t)}\\
    \le& \esp{\Gamma^{\lambda}_T \Psi(\xi,\xi') + 2\nu\int_0^T \Gamma^{\lambda}_t |\delta f'_t|^2 \ud t}\\
    &+ 4C_{\eqref{eq for later one}}\esp{\sup_{t\in[0,T]}\Gamma^{\lambda}_t \Psi(Y_t,Y'_t)}^{\frac12}\esp{\int_0^T \Gamma^{\lambda}_s|Z_s - R(\theta(Y_s,Y'_s))Z'_s|^2  \ud s}^\frac12.
\end{align*}
which gives us, thanks to Young inequality,
\begin{align*}
    \frac12 \esp{\sup_{t\in[0,T]}\Gamma^{\lambda}_t \Psi(Y_t,Y'_t)}
    \le& \esp{\Gamma^{\lambda}_T \Psi(\xi,\xi') + 2\nu\int_0^T \Gamma^{\lambda}_t |\delta f'_t|^2 \ud t}\\
    &+ 8C_{\eqref{eq for later one}}^2\esp{\int_0^T \Gamma^{\lambda}_s|Z_s - R(\theta(Y_s,Y'_s))Z'_s|^2  \ud s}.
\end{align*}
Combining the previous inequality with \eqref{eq weak estim one}, we get
\begin{align*}
    &\esp{\sup_{t\in[0,T]}\Gamma^{\lambda}_t \Psi(Y_t,Y'_t)} + \esp{\int_0^T \Gamma^{\lambda}_s|Z_s - R(\theta(Y_s,Y'_s))Z'_s|^2  \ud s}\\
    \le& (4+32C_{\eqref{eq for later one}}^2)\esp{\Gamma^{\lambda}_T \Psi(\xi,\xi')} + (8+64C_{\eqref{eq for later one}}^2)\nu\esp{\int_0^T \Gamma^{\lambda}_t |\delta f'_t|^2 \ud t},
\end{align*}
which gives us \eqref{eq strong estim one} since we can set $\nu$ as small as we want.
\\
{\bf Step 2.} We compute 
 \begin{align}\label{eq starting point 2}
    \esp{\int_0^T|Z_t - Z'_t|^2 \ud t}
    \le 2 \esp{
        \int_0^T \Gamma^{\lambda}_t|Z_t - R(\theta(Y_t,Y'_t))Z'_t|^2 \ud t
    + \int_0^T \Gamma^{\lambda}_t |(I- R(\theta(Y_t,Y'_t)))Z'_t|^2 \ud t
    }.
 \end{align}
 We have 
 \begin{align*}
    \int_0^T\Gamma^{\lambda}_t |(I- R(\theta(Y_t,Y'_t)))Z'_t|^2 \ud t
    &\le C \esp{\int_0^T\Psi(Y_t,Y'_t)|Z'_t|^2 \ud t}
    \\
    &\le 
    C \esp{\sup_{t \in [0,T]}\Psi(Y_t,Y'_t)^2}^\frac12\esp{(\int_0^T|Z'_t|^2 \ud t)^2}^\frac12
    &
    \\
    &\le 2C_{\eqref{eq Prop-aprioriZ}}(1+\||f'(.,0,0)|^{1/2}\|_{\BP{2}}^2)\esp{\sup_{t \in [0,T]}\Psi(Y_t,Y'_t)}^\frac12
 \end{align*}
 where for the last inequality we used \eqref{energy-ineq-Z}, Proposition \ref{Prop-aprioriZ} and the boundedness of $\cD$.
 Let us remark that, thanks to Proposition \ref{Prop-aprioriZ}, $\|Z'\|_{\BP{2}}$ is upper-bounded by a constant that depend only on  Inserting back the previous inequality into  \eqref{eq starting point 2} and using \eqref{eq strong estim one}, we obtain \eqref{eq cooked for Z}.
\eproof

\smallskip

We easily deduce the uniqueness of a solution to the reflected BSDE \eqref{eq reflected bsde} as a corollary of the stability result of Proposition \ref{prop stability}.

\begin{Theorem}[Uniqueness for reflected BSDE with a Brownian filtration]
\label{Thm-Uniqu2}
Let Assumption \ref{ass:brownian-case} hold and consider $(Y,Z,K)$ and $(Y',Z',K')$, two solutions of the reflected BSDE \eqref{eq reflected bsde} with the same $f$ and $\xi$. Then, $(Y,Z,K)=(Y',Z',K')$.
\end{Theorem}

\proof 
Let us apply Proposition \ref{prop stability} with $\alpha_t = 1+ |Z'_t| + |f'(t,Y'_t,Z'_t)|$, $\xi=\xi'$ and $f=f'$. Let us remark that \eqref{condition alpha} is satisfied thanks to \eqref{moment-expo-f(t,y,z)}. Then  
\begin{align*}
    \esp{\sup_{t\in[0,T]} \Psi(Y_t,Y'_t)}  \le \esp{\sup_{t\in[0,T]}\Gamma^{\lambda}_t \Psi(Y_t,Y'_t)} +\esp{\int_0^T\Gamma^{\lambda}_t|Z_t - R(\theta(Y_t,Y'_t))Z'_t|^2 \ud t }\le 0.  
\end{align*}
We obtain $Y=Y'$ and the previous inequality becomes $\esp{\int_0^T\Gamma^{\lambda}_t|Z_t - Z'_t|^2 \ud t }=0$ which gives us $Z=Z'$. The uniqueness of $K$ follows easily.
\eproof

\medskip

We finish this subsection with another stability result that is used in the next subsection.

\begin{Proposition}
    \label{prop stability 2}
    { Consider any adapted processes  $f=(f_s)_{s \in [0,T]}$ and $f'=(f'_s)_{s \in [0,T]}$ with $|f|^{1/2},\,|f'|^{1/2} \in \widebar{\mathscr{H}^{\infty}}^{\BP{2}}$, and let Assumption \ref{ass:brownian-case} hold.}
    Let $(Y,Z,K)$ (resp., $(Y',Z',K')$) be a solution of the reflected BSDE \eqref{eq reflected bsde} with the terminal condition $\xi$ (resp., $\xi'$) and with the generator $f$ (resp., $f'$). 
    For all $\lambda>0$, we set $\Gamma^{\lambda}_t = e^{\lambda \int_0^t |f'(s)|\ud s}$ and assume that $\esp{\Gamma^{\lambda}_T}<+\infty$. Then, there exist $C>0$ and $\lambda_0 \ge 0$, only depending on $\bar{\cD}$ and $T$, such that, for all $\lambda \ge \lambda_0$, we have
    \begin{align}\label{eq strong estim one 2}
        &\esp{ \sup_{t \in [0,T ]} \Psi(Y_t,Y'_t)} +\esp{\int_0^T\Gamma^{\lambda}_t|Z_t - R(\theta(Y_t,Y'_t))Z'_t|^2 \ud t }\\ \le&  C\mathbb{E}\left[ \Gamma^{4\lambda}_T \right]^\frac14 \mathbb{E}\left[ \Psi(\xi,\xi')^4\right]^\frac14
        + C\mathbb{E}\left[ \Gamma^{4\lambda}_T \right]^\frac14 \| |f(.)-f'(.)|^{1/2} \|_{\BP{2}}^2     \nonumber  
    \end{align}
    and
    \begin{align}\label{eq cooked for Z 2}
        & \esp{\int_0^T|Z_t - Z'_t|^2 \ud t }\\
        \le& C\mathbb{E}\left[ \Gamma^{4\lambda}_T \right]^\frac14 \mathbb{E}\left[ \Psi(\xi,\xi')^4\right]^\frac14
        + C\mathbb{E}\left[ \Gamma^{4\lambda}_T \right]^\frac14 \| |f(.)-f'(.)|^{1/2} \|_{\BP{2}}^2 \nonumber\\
        &+ \nonumber C(1+\||f'(.)|^{1/2}\|_{\BP{2}}^2)\left(\mathbb{E}\left[ \Gamma^{4\lambda}_T \right]^\frac14 \mathbb{E}\left[ \Psi(\xi,\xi')^4\right]^\frac14
        +\mathbb{E}\left[ \Gamma^{4\lambda}_T \right]^\frac14 \| |f(.)-f'(.)|^{1/2} \|_{\BP{2}}^2 \right)^\frac12.
    \end{align}
\end{Proposition}

\proof
    We start with same computations as in the proof of Proposition \ref{prop stability}. Using Proposition \ref{Ito-formula-Psi-bis}, we get, for $0 \le t \le T$,
    \begin{align}
        \Gamma^\lambda_T \Psi(Y_T,Y'_T) \ge& \Gamma^\lambda_t \Psi(Y_t,Y'_t) + \lambda \int_t^T |f'(s)| \Gamma^\lambda_s \Psi(Y_s,Y'_s) \ud s + 2 \int_t^T \Gamma^\lambda_s \left(\overrightarrow{Y_s Y'_s}f_s + \overrightarrow{Y'_s Y_s} f'_s\right)  \ud s  \nonumber
        \\
        &- 2 \int_t^T \Gamma^\lambda_s \left(\overrightarrow{Y_s Y'_s}Z_s + \overrightarrow{Y'_s Y_s}Z'_s\right)  \ud W_s +\int_t^T\Gamma^\lambda_s|Z_s - R(\theta(Y_s,Y'_s))Z'_s|^2 \ud s \label{eq starting point one 2}
    \end{align}
    and we also compute
    \begin{align*}
        \overrightarrow{Y_s Y'_s}f_s + \overrightarrow{Y'_s Y_s} f'_s
        &\ge -\Psi(Y_s,Y'_s)^{1/2}|f(s)-f'(s)| -C|f'(s)|\Psi(Y_s,Y'_s) \ud s,
    \end{align*} 
    Where $C$ only depends on $\bar{\cD}$. Inserting the previous inequality back into \eqref{eq starting point one 2}, we obtain, for $\lambda$ large enough,
    \begin{align}
        \Gamma^{\lambda}_T \Psi(Y_T,Y'_T) \ge& \Gamma^\lambda_t \Psi(Y_t,Y'_t) - C \int_t^T \Gamma^{\lambda}_s \Psi(Y_s,Y'_s)^{1/2}|f(s)-f'(s)| \ud s\nonumber\\  
        &- 2 \int_t^T \Gamma^{\lambda}_s \left(\overrightarrow{Y_s Y'_s}Z_s + \overrightarrow{Y'_s Y_s}Z'_s\right)  \ud W_s + \int_t^T\Gamma^{\lambda}_s|Z_s - R(\theta(Y_s,Y'_s))Z'_s|^2 \ud s.  \label{eq intrem point one 2}
    \end{align}
    As proved in the proof of Proposition \ref{prop stability}, see inequality after \eqref{eq for later one}, the stochastic integral in the previous inequality is a martingale.
    Then, we take conditional expectation on both sides of \eqref{eq intrem point one 2} and we get
    \begin{align}
        &\Gamma^{\lambda}_t\Psi(Y_t,Y'_t) + \mathbb{E}_t\left[\int_0^T\Gamma^{\lambda}_t|Z_t - R(\theta(Y_t,Y'_t))Z'_t|^2 \ud t \right] \nonumber\\
        \le& \mathbb{E}_t\left[ \Gamma^{\lambda}_T \Psi(\xi,\xi') 
        + C \int_t^T \Gamma^{\lambda}_s \Psi(Y_s,Y'_s)^{1/2}|f(s)-f'(s)| \ud s \right] \nonumber \\
        \le& \mathbb{E}_t\left[ \Gamma^{\lambda}_T \Psi(\xi,\xi') 
        + C \Gamma^{\lambda}_T \int_t^T  |f(s)-f'(s)| \ud s \right] \label{eq weak estim one 2}
    \end{align}
    By using Cauchy-Schwarz inequality and energy inequality we obtain the following upper-bound:
    \begin{align}
        \mathbb{E}_t\left[ \Gamma^{\lambda}_T\int_t^T |f(s)-f'(s)| \ud s\right] \le 2\mathbb{E}_t\left[ \Gamma^{2\lambda}_T \right]^\frac12 \| |f(.)-f'(.)|^{1/2} \|_{\BP{2}}^2. \label{eq stab 2 - 2}
    \end{align}
    Finally, we put \eqref{eq stab 2 - 2} into \eqref{eq weak estim one 2}, we apply Cauchy-Schwarz for the first term, we take the supremum on $t \in [0,T]$, we take the expectation and we apply Doob maximal inequality Cauchy-Schwarz to obtain \eqref{eq strong estim one 2}.
    The proof of \eqref{eq cooked for Z 2} follows the same lines as for the proof of inequality \eqref{eq cooked for Z}.
\eproof

\subsection{Existence of a \texorpdfstring{$\Gamma$}{Gamma}-martingale with given terminal value and exogenous drift}

We now turn to the existence results. Throughout this subsection, we assume that the filtration is Brownian and use the Kendall's approach \cite{Kendall-90} which consists of a recursive application of the (conditional) Fr\'echet mean.

{
An existence result for $\Gamma$-martingales in $CAT(0)$ spaces, with a prescribed terminal condition and with zero drift, is already proved in \cite{Christiansen-Sturm-08} (see Theorem 2.5). Their strategy follows the same lines as \cite{Kendall-90}: i.e., the authors consider a size-$n$ partition of the time interval, iterate the Fréchet mean over this partition, and show that the resulting stochastic process converges to a limiting process when $n$ tends to $+\infty$. Nevertheless, \cite{Christiansen-Sturm-08} uses a different definition of a $\Gamma$-martingale, which is proven to coincide with the more canonical definition, stated via $\Gamma$-convex functions, only for Riemannian manifolds without a boundary. In the present paper, we consider manifolds with boundaries, hence, we cannot directly apply the existence result of \cite{Christiansen-Sturm-08}, even for a $\Gamma$-martingale without drift.

}

As mentionned, our approach relies heavily on the existence and uniqueness of a Fréchet mean and the corresponding Jensen's inequality, which are known to hold true in $CAT(0)$ spaces, see \cite{Sturm-02}. These properties are summarized in the following proposition.

\begin{Proposition}
    Let Assumption \ref{ass:3} hold, and let $\xi$ be a random variable with values in $\bar{\cD}$. Then, there exists a unique minimizer of
    \begin{align}
        \inf_{x \in \bar{\cD}} \esp{\Psi(x,\xi)},
    \end{align}
    which is denoted by $\cE(\xi)$. Moreover, we have 
    \begin{align}\label{eq the jensen inequality}
        \Psi(\cE(\xi),\cE(\xi')) \le \esp{\Psi(\xi,\xi')},
    \end{align}
    for all random variables $\xi,\xi'$ with values in $\bar{\cD}$.
\end{Proposition}

For the sake of completeness, a proof of the above result, tailor-fitted to our setting, is presented in Section \ref{subse about frechet mean} (note that $\bar{\cD}$ satisfies Assumption \ref{ass:4} by a direct application of Proposition \ref{Prop-Psi}).


\subsubsection{\texorpdfstring{$\Gamma$}{Gamma}-martingales with exogenous drift}

We start by studying the Markovian case. Let us consider an $\R^M$-valued process $X$ defined as the solution of the following SDE:
\begin{align} \label{eq SDE}
    \ud X_t = b(t,X_t) \ud t + \sigma(t,X_t) \ud W_t
\end{align}
where $b$ and $\sigma$ are bounded measurable functions, and $x \mapsto (b(t,x),\sigma(t,x))$ is Lipschitz continuous uniformly in time.
Under these assumptions, there exists a unique strong solution which satisfies:
\begin{align}
    \esp{\sup_{t \in [0,T]}|X_t^{s,x}|^p} \le C_p,\; p \ge 1,
\end{align}
where $(X^{s,x}_t)$ denotes classically the solution of \eqref{eq SDE} starting at time $s$ from $x$ and by convention is constantly equal to $x$ for time before $s$. It is also easily obtained that 
\begin{align}\label{eq estim X - X'}
    \esp{|X_t^{s,x} - X_t^{s,x'}|^2}\le e^{C(t-s)}|x-x'|^2\;,
\end{align}
for some positive constant $C$.\\




\begin{Proposition}\label{prop:Gamma.mtg.drift.Markov}
{Let Assumption \ref{ass:3} hold, and let $X$ be the solution of \eqref{eq SDE} with an arbitrary initial condition $X_0=x\in \R^M$. Consider any bounded Lipschitz functions $F:\R^M \rightarrow \R^d$, $g:\R^M \rightarrow \bar{\cD}$, and let $\xi=g(X_{T})$, $f_t(\omega):=F(X_t(\omega))$, $\omega \in \Omega$.} Then, there exists a $\Gamma$-martingale $Y$ with the drift $f$ and with the terminal value $Y_T = \xi$. Moreover, thanks to Theorem \ref{Thm-Uniqu1}, such a $\Gamma$-martingale is unique, {and $Y_t$ is given by a Lipschitz function of $X_t$}.
\end{Proposition}

\proof 
We construct a solution using the Kendall's approach, which relies on a sequence of approximations. For any $n \in \mathbb{N}^*$ and any partition of the time interval $[0,T]$,
$$\pi^n := \set{0=t_0<t_1<\cdots<t_{n-1}<t_n=T},$$  
we define the function $g^n$ recursively: $g^n(T,x) := g(x)$ and
\begin{align}
& \tilde g^n(t,x) = \mathcal{E}[g^n(t_i,X^{t,x}_{t_i})] \text{ when } t_{i-1} \leqslant t <t_i, \label{eq de approx scheme-1} \\
& g^n(t,x) = R^{x,\tilde g^n(t,x),i}(t),\quad t_{i-1} \leqslant t <t_i,\quad i=1,\ldots,n, \label{eq de approx scheme-2}
\end{align}
where $R^{x,y,i}(t)$ is defined as the unique solution to the reflected ODE:
\begin{align}
& -\ud R^{x,y,i}(t) = F(x)\,\ud t - \ud K^{x,y,i}_t,\quad t\in(t_{i-1},t_i),
\label{eq dyn R-1}
\\
& \ud K^{x,y,i}_t \in \mathfrak{n}(R^{x,y,i}(t))\,\ud\mathrm{Var}_t(K^{x,y,i}),\quad R^{x,y,i}(t_i) = y\, \text{ and } \, {K^{x,y,i}_{t_i} = 0}.
\label{eq dyn R-2}
\end{align}
To ensure the well-posedness of the above ODE and to obtain important estimates used further in the proof, we state the following intermediary result.

\begin{Lemma}\label{le control ode}
    Under Assumption \ref{ass:3} and provided that $F$ is Lipschitz and bounded, for any $n\geq 1$, $i=1,\ldots,n$, $x,y\in\bar\cD$, there exists a unique continuous solution $(R^{x,y,i},K^{x,y,i})$ to the system \eqref{eq dyn R-1}-\eqref{eq dyn R-2}, with $R^{x,y,i}$ taking values in $\bar\cD$.
    Moreover, there exist $\bar C_1,\bar C_2>0$ s.t., for all $n\geq 1$, $i=1,\ldots,n$, $x,x'\in\R^{d'}$, $y,y'\in\bar\cD$, and $t^n_{i-1}\leq t<t'\leq t^n_{i}$, we have
    \begin{align}
    & \Psi(R^{x,y,i}(t),R^{x,y',i}(t)) \leq (1+\bar C_1\,h)\,\Psi(y,y'),
    \label{eq le stab R - psi}
    \\
    & |R^{x,y,i}(t)-R^{x',y,i}(t)| \leq \bar C_1\,h\,(|x-x'|\wedge 1),
    \label{eq le stab R - psi.new}
    \\
    & |R^{x,y,i}(t') - R^{x,y,i}(t)| = |K^{x,y,i}_{t'} - K^{x,y,i}_{t} - F(x)(t'-t)| \leq \bar C_2\, (t'-t).
    \label{eq le stab R - time}
    \end{align}
\end{Lemma}
\proof
The existence and uniqueness statement, as well as the estimate \eqref{eq le stab R - time}, follow directly from Theorem 2.2 in \cite{Lions-Sznitman-84} { noting that their assumption (5) is fulfiled thanks to Proposition \ref{prop:extsphere:interiorcone}}. To show \eqref{eq le stab R - psi}--\eqref{eq le stab R - psi.new}, we recall Proposition \ref{Prop-Psi}, equation \eqref{expansion-log}, as well as the inequalities $|y-y'|^2\leq \Psi(y,y')\leq C_1\,|y-y'|^2$ (see Proposition \ref{Prop-equivalence-Psi-Euclidean}), to obtain
\begin{align*}
& \langle\nabla\Psi(R^{x,y,i}(t),R^{x,y',i}(t)),(dK^{x,y,i}_t,dK^{x,y',i}_t)\rangle \geq 0,\\
& \nabla\Psi(R^{x,y,i}(t),R^{x,y',i}(t))= -2(R^{x,y',i}(t)-R^{x,y,i}(t), R^{x,y,i}(t)-R^{x,y',i}(t))\\
& + O(\Psi(R^{x,y,i}(t),R^{x,y',i}(t))),
\end{align*}
where $|O(\Psi(R^{x,y,i}(t),R^{x,y',i}(t)))|\leq C_1\,\Psi(R^{x,y,i}(t),R^{x,y',i}(t))$ for all $x,y,y',t,i$.
Using the above, we obtain:
\begin{align*}
&d\Psi(R^{x,y,i}(t),R^{x,y',i}(t)) 
= -\langle\nabla\Psi(R^{x,y,i}(t),R^{x,y',i}(t)),(F(x),F(x))\rangle\,\ud t\\
&+\langle\nabla\Psi(R^{x,y,i}(t),R^{x,y',i}(t)),(dK^{x,y,i}_t,dK^{x,y',i}_t)\rangle\\
&\geq 2\langle R^{x,y',i}(t)-R^{x,y,i}(t),F(x)\rangle\,\ud t + 2\langle R^{x,y,i}(t)-R^{x,y',i}(t),F(x)\rangle\,\ud t\\
& - C_2\,\Psi(R^{x,y,i}(t),R^{x,y',i}(t))\,dt=- C_2\,\Psi(R^{x,y,i}(t),R^{x,y',i}(t))\,dt.
\end{align*}
An application of Gronwall's inequality yields \eqref{eq le stab R - psi}.

To verify \eqref{eq le stab R - psi.new}, we recall $|y-y'|^2\leq \Psi(y,y')$ (recall Proposition \ref{Prop-equivalence-Psi-Euclidean}) and proceed as before:
\begin{align*}
&d\Psi^{1/2}(R^{x,y,i}(t),R^{x',y,i}(t))\\
&\geq -\frac{1}{2}\Psi^{-1/2}(R^{x,y,i}(t),R^{x',y,i}(t))\langle\nabla\Psi(R^{x,y,i}(t),R^{x',y,i}(t)),(F(x),F(x'))\rangle\,\ud t\\
&\geq \Psi^{-1/2}(R^{x,y,i}(t),R^{x',y,i}(t))\langle R^{x',y,i}(t)-R^{x,y,i}(t),F(x)-F(x')\rangle\,\ud t \\
& - C_3\,\Psi^{1/2}(R^{x,y,i}(t),R^{x,y',i}(t))\,dt\\
&\geq -C_4\,(|x-x'|\wedge1)\,\ud t - C_3\,\Psi^{1/2}(R^{x,y,i}(t),R^{x,y',i}(t))\,dt.
\end{align*}
Another application of Gronwall's inequality yields \eqref{eq le stab R - psi.new}.
\eproof.

\medskip  
    
We now consider an increasing sequence of dyadic partitions $\pi^n:=\set{t^n_i = ih, 0\le i \le n:=2^k, h := 2^{-k}}$, $k \ge 1$. For the readers convenience, we shall suppress below the dependence on $n$ for the time in the partition and we denote by $(g^n,\tilde{g}^n)$ the scheme built in \eqref{eq de approx scheme-1}-\eqref{eq de approx scheme-2}.
    
    \medskip
    
    {\bf Step 1.} In this step, we show that $\{g^n\}$ is equicontinuous in $[0,T] \times \R^{d'}$. 
    To this end, we notice that $|y-y'|^2\leq \Psi(y,y')\leq C_1\,|y-y'|^2$ (recall Proposition \ref{Prop-equivalence-Psi-Euclidean}) and apply \eqref{eq le stab R - psi} along with \eqref{eq the jensen inequality}  and Young's inequality, to obtain, for $t\in[t_{i-1},t_i]$: 
    \begin{align*}
    &\Psi(g^n(t,x),g^n(t,x')) = \Psi(R^{x,\tilde g^n(t,x),i}(t),R^{x',\tilde g^n(t,x'),i}(t))\\
    &\leq (1+h)\,\Psi(R^{x,\tilde g^n(t,x),i}(t),R^{x,\tilde g^n(t,x'),i}(t)) + (1+h^{-1})\, \Psi(R^{x,\tilde g^n(t,x'),i}(t),R^{x',\tilde g^n(t,x'),i}(t))\\
    & \leqslant (1+C_2\,h)\, \left[\Psi(\tilde g^n(t,x),\tilde g^n(t,x'))\right] + C_3\,h\,|x-x'|^2\\
    &\leqslant (1+C_2\,h)\, \mathbb{E}\left[\Psi\left(g^n(t_i,X^{t,x}_{t_i}),\,g^n(t_i,X^{t,x'}_{t_i})\right)\right] + C_3\,h\,|x-x'|^2.
    \end{align*}
Now, we observe that 
\begin{align*}
    \mathbb{E}\left[\Psi\left(g^n(t_i,X^{t,x}_{t_i}),\,g^n(t_i,X^{t,x'}_{t_i})\right)\right] &=\mathbb{E}\left[\frac{\Psi\left(g^n(t_i,X^{t,x}_{t_i}),\,g^n(t_i,X^{t,x'}_{t_i})\right)}{|X^{t,x}_{t_i}-X^{t,x'}_{t_i}|^2\bone_{\set{X^{t,x}_{t_i}\neq X^{t,x'}_{t_i}}}} |X^{t,x}_{t_i}-X^{t,x'}_{t_i}|^2 \right]
    \\
    &\leqslant \sup_{x''\neq x'''\in\R^{d'}} \frac{\Psi(g^n(t_i,x''),g^n(t_i,x'''))}{|x''-x'''|^2} \esp{|X^{t,x}_{t_i}-X^{t,x'}_{t_i}|^2}
    \\
    &\leqslant (1+ C_4h)\, |x-x'|^2\,\sup_{x''\neq x'''\in\R^{d'}} \frac{\Psi(g^n(t_i,x''),g^n(t_i,x'''))}{|x''-x'''|^2},
\end{align*}
where we used \eqref{eq estim X - X'} to get the last inequality. This leads to 
\begin{align*}
\Psi(g^n(t,x),g^n(t,x'))
    \le& (1+C_5\,h)\, |x-x'|^2\,\sup_{x''\neq x'''\in\R^{d'}} \frac{\Psi(g^n(t_i,x''),g^n(t_i,x'''))}{|x''-x'''|^2}\\
    &+ C_3\,h\,|x-x'|^2,\quad\quad t\in[t_{i-1},t_i].
\end{align*}

Iterating the above, we obtain, for all $t\in[0,T]$ and $n\geq 0$:
    \begin{align*}
    &\sup_{x\neq x'\in\R^{d'}} \frac{\Psi(g^n(t,x),g^n(t,x'))}{|x-x'|^2}
    \leqslant (1+C_5\,h)^n\,\sup_{x\neq x'\in\R^{d'}} \frac{\Psi(g(x),g(x'))}{|x-x'|^2}\\
     &+ C_3\,h\,\sum_{i=0}^{n-1} (1+C_5\,h)^i
     = (1+C_5\,h)^n\,\sup_{x\neq x'\in\R^{d'}} \frac{\Psi(g(x),g(x'))}{|x-x'|^2} 
     + \frac{C_3}{C_5}\,((1+C_5\,h)^n-1).
    \end{align*}
    Recalling again that $|y-y'|^2\leq \Psi(y,y')\leq C_1\,|y-y'|^2$ and using the above along with the Lipschitz property of $g$, we obtain
    \begin{align}
    &|g^n(t,x) - g^n(t,x')|^2 \leq \Psi(g^n(t,x),g^n(t,x'))\leqslant C_6 \,|x-x'|^2.\label{eq.RBSDE.exist.newEq.1}
    \end{align}

\smallskip

Next, for $t_{i-1} \leqslant t<t'<t_i$, by using once again Proposition \ref{Prop uniqueness Frechet} and  $|y-y'|^2\leq \Psi(y,y')$, we obtain
\begin{align}
&|\tilde g^n(t,x)-\tilde g^n(t',x)|^2 \leqslant \Psi(\tilde g^n(t,x),\tilde g^n(t',x))\nonumber \\
& \leqslant \mathbb{E}\Psi\left(g^n(t_i,X^{t,x}_{t_i}),\,g^n(t_i,X^{t',x}_{t_i})\right)
\leqslant C_7\,\mathbb{E}[|X_{t_i}^{t,x}-X_{t_i}^{t',x}|^2] \le C_8\,(t' -t). \label{eq useful comp}
\end{align}

Using Lemma \ref{le control ode} and $|y-y'|^2\leq\Psi(y,y')\leq C_1\,|y-y'|^2$, we obtain
\begin{align}
&|g^n(t,x)-g^n(t',x)|^2 = |R^{x,\tilde g^n(t,x),i}(t) - R^{x,\tilde g^n(t',x),i}(t')|^2\nonumber\\
&\leq 2|R^{x,\tilde g^n(t,x),i}(t) - R^{x,\tilde g^n(t',x),i}(t)|^2
+  2|R^{x,\tilde g^n(t',x),i}(t) - R^{x,\tilde g^n(t',x),i}(t')|^2\nonumber\\
&\leq 2\Psi(R^{x,\tilde g^n(t,x),i}(t),R^{x,\tilde g^n(t',x),i}(t)) + C_9\,(t'-t)^2\label{eq.RBSDE.exist.newEq.2}\\
&\leq C_{10}\Psi(\tilde g^n(t,x),\tilde g^n(t',x)) + C_9\,(t'-t)^2 \leq C_{11}\,(t'-t),\nonumber
\end{align}
where we used \eqref{eq useful comp} to obtain the last inequality.
Finally, for $t_{i-1}\leq t\leqslant t_i\leq t_j\leq t'\leq t_{j+1}$, we have, using the previous inequality, 
\begin{align*}
&|g^n(t,x)-g^n(t',x)|^2 \leq 2\left[|g^n(t_i,x)-g^n(t,x)|^2 \right.\\
&\left. + \sum_{l=i}^{j-1} |g^n(t_l,x)-g^n(t_{l+1},x)|^2 + |g^n(t',x)-g^n(t_j,x)|^2\right]
\leq 2 C_{11}\,(t'-t),
\end{align*}
which completes the proof of equicontinuity of $\{g^n\}$.

\medskip

{\bf Step 2.} 
As $\{g^n\}$ is uniformly bounded and equicontinuous, we apply the Arzela-Ascoli Theorem to extract a subsequence of $\{g^n\}_{n=1}^\infty$ that converges locally uniformly, and denote the limiting function by $G$. It remains to prove that $Y:=G(.,X^{0,x}_.)$ is a $\Gamma$-martingale with drift $f$. (For the remainder of the proof, we omit the superscript $(0,x)$ in $X$).
To this end, we first notice that, by construction, $Y$ is a continuous process adapted to $(\mathcal{F}_t)$ and taking values in $\bar\cD$.
%
Then, according to Corollary \ref{cor-global gamma convex function enough}, in order to show that $Y$ is a $\Gamma$-martingale with drift $f$, it suffices to prove that, for any $0\leq t<t'\leq T$ and any \emph{global} special $\Gamma$-convex function $\psi$, we have
\begin{align}
\EE\left[ \psi(Y_{t'}) - \psi(Y_{t}) 
+\int_{t}^{t' } \nabla \psi(Y_s)\cdot F(X_s) \ud s
\,\vert\,\mathcal{F}_{t}\right]\geq0\quad a.s..
\label{eq.exist.pf.step 2.eq2}
\end{align}
We note that it suffices to verify the above inequality for $t,t'$ that are dyadic rational, and we choose $n_0\ge 1$ large enough, so that $t,t'\in\pi^n$ for all $n\geq n_0$. We introduce $j(n)$ and $j'(n)$, such that $t=t_{j(n)}<t'=t_{j'(n)}$ (where we drop the superscript `$n$' in the elements of a partition, to ease the notation). 
Let us denote $(Y^n,\tilde Y^n) := (g^n,\tilde g^n)(.,X_.)$. Recalling that, a.s., $Y^n_t\rightarrow Y_t$ for all $t\in[0,T]$, we deduce that 
\begin{align}
    &\EE\left[ \psi(Y^{n}_{t'}) - \psi(Y^{n}_{t}) +\int_{t}^{t'} \nabla \psi(Y^n_s)\cdot F(X_s) \ud s \,\vert\,\mathcal{F}_{t}  \right] 
    \label{eq starting point kendall}
\end{align}
converges in $L^p$ (for any $p \ge 1$) to the left hand side of \eqref{eq.exist.pf.step 2.eq2}.
Thus, we aim to estimate \eqref{eq starting point kendall} from below. To do so, we decompose it into the sum of $A_n+\tilde{A}_n$ where 
\begin{align}\label{eq de A-term}
A_n &:= \EE\left[ \sum_{i=j}^{j'-1}\,\EE\left[ \psi(Y^{n}_{t_{i+1}}) - \psi(\tilde{Y}^{n}_{t_i}) 
\,\vert\,\mathcal{F}_{t_i}\right]\,\vert\,\mathcal{F}_{t}\right]
\end{align}
and 
\begin{align} 
    &\tilde A_n := \EE\left[ \sum_{i=j}^{j'-1}\EE\left[ \psi(\tilde{Y}^{n}_{t_i}) - \psi(Y^{n}_{t_i}) 
    +\int_{t_i}^{t_{i+1} } \nabla \psi(Y^n_s)\cdot F(X_s) \ud s
    \,\vert\,\mathcal{F}_{t_i}\right]\,\vert\,\mathcal{F}_{t}\right].
    \label{eq de tilde A}
    \end{align}

{\bf Step 2.a} We first study the $A$-term given in \eqref{eq de A-term} by observing that 
\begin{align*}
    \EE\left[ \psi(Y^{n}_{t_{i+1}}) - \psi(\tilde{Y}^{n}_{t_i}) 
\,\vert\,\mathcal{F}_{t_i}\right] =
\EE\left[ \psi(Y^{n}_{t_{i+1}}) - \psi(\tilde{Y}^{n}_{t_i}) 
\,\vert\,X_{t_i}\right],
\end{align*}
by the Markov property of $X$ (recall that $Y^{n}_{t_{i+1}}$ and $\tilde{Y}^{n}_{t_i}$, respectively, are functions of $X_{t_{i+1}}$ and $X_{t_i}$).
We then compute, for any bounded measurable function $\theta$,
\begin{align*}
    \esp{\psi(\tilde{Y}^{n}_{t_i}) \theta(X_{t_i})}
    &=
    \int{\psi(\tilde{g}^{n}(t_i,x)) \theta(x) \ud \P_{X_{t_i}}(x)}
    \\
    &=\int{\psi(\mathcal{E}[g^n(t_{i+1},X^{t_i,x}_{t_{i+1}})]) \theta(x) \ud \P_{X_{t_i}}(x)}
    \\
    &\le \int{ \esp{\psi(g^n(t_{i+1},X^{t_i,x}_{t_{i+1}}))} \theta(x) \ud \P_{X_{t_i}}(x)}
    \\
    &= \esp{\esp{\psi(g^n(t_{i+1},X^{t_i,X_{t_i}}_{t_{i+1}}))\theta(X_{t_i})\,\vert\,\mathcal{F}_{t_i}}}
    \\
    &= \esp{\psi(g^n(t_{i+1},X_{t_{i+1}}))\theta(X_{t_i})} = \esp{\psi(Y^{n}_{t_{i+1}})\theta(X_{t_i})},
\end{align*}
where we used Proposition \ref{prop Jensen for Frechet} to obtain the inequality, and the flow property of $(X^{t,x})$ to obtain the second-to-last equality. We thus deduce that
\begin{align*}
\EE\left[ \psi(Y^{n}_{t_{i+1}}) - \psi(\tilde{Y}^{n}_{t_i}) \,\vert\,X_{t_i}\right]\ge 0
\end{align*}
and, hence,
\begin{align}\label{eq bound for A-term}
 A_n \ge 0,
\end{align}
according to \eqref{eq de A-term}.

\smallskip

{\bf Step 2.b} We now turn to the $\tilde{A}$-term defined in \eqref{eq de tilde A}. 
Let us bound from below the terms $\esp{\psi(\tilde{Y}^{n}_{t_i}) - \psi(Y^{n}_{t_i}) | \cF_{t_i}}$. First, recalling \eqref{eq de approx scheme-2} and \eqref{eq dyn R-2}, we notice that $g^n(t_i,x) = R^{x,\tilde{g}^n(t_i,x),i+1}(t_i)$ and $\tilde g^n(t_{i},x) = R^{x,\tilde{g}^n(t_i,x),i+1}(t_{i+1})$, and, hence, according to \eqref{eq dyn R-1},
\begin{align*}
    \psi(\tilde{g}^n(t_i,x)) - \psi({g}^n(t_i,x)) 
    &= \int_{t_i}^{t_{i+1}} \nabla \psi(R(s))\left(-F(x) \ud s + \ud K_s\right)
    \\
    & \ge - \int_{t_i}^{t_{i+1}} \nabla \psi(R(s)) F(x) \ud s,
\end{align*}
where we used the defining property of special $\Gamma$-convex functions to obtain the last inequality, and we dropped the superscript of $R$ (here and throughout the remainder of the proof). Denoting by $\omega_{\nabla \psi}$ the modulus of continuity of $\nabla \psi$, we compute
\begin{align*}
    \psi(\tilde{g}^n(t_i,x)) - \psi({g}^n(t_i,x)) 
    \ge - h \nabla \psi({g}^n(t_i,x)) F(x) 
    - |F|_\infty h \,\omega_{\nabla \psi}(\bar{C}_2 h),
\end{align*}
where we used Lemma \ref{le control ode}. We thus conclude that
\begin{align}\label{eq useful A tilde one}
    \esp{\psi(\tilde{Y}^{n}_{t_i}) - \psi(Y^{n}_{t_i}) | \cF_{t_i}}
    \ge - h \nabla \psi(Y^n_{t_i}) F(X_{t_i}) 
    - |F|_\infty h \,\omega_{\nabla \psi}(\bar{C}_2 h).
\end{align}
Next, we denote 
$$
\omega^i_X(h) := \sup_{s\in[t_i,t_{i+1}]}|X_s-X_{t_i}|,
$$
and notice (by examining the SDE satisfies by $X$) that there exist identically distributed random variables $\{\eta_i\}$, whose distribution does not depend on $h$, such that $\omega^i_X(h)\leq C_{12} \,h^{1/2}\,\eta_i$ and $\eta_i$ is independent of $\mathcal{F}_{t_i}$, for each $i=1,\ldots,n$.
Using the latter observation, we obtain:
\begin{align*}
&\EE\left[\int_{t_i}^{t_{i+1} }| \nabla \psi(Y^n_s)\cdot F(X_s) - \nabla \psi(Y^n_{t_i}) F(X_{t_i}) |\ud s\,\vert\,\mathcal{F}_{t_i}\right]\\
&\le C_{13}\,h\, \EE\left[\omega_{\nabla\psi}\left(\sup_{s\in[t_i,t_{i+1}]}|g^n(s,X_{s})-g^n(t_i,X_{t_i})|\right) + \mathrm{Lip}(F)\,(\omega^i_X(h)\wedge C_{14})\,\vert\,\mathcal{F}_{t_i}\right]\\
&\leq C_{15}\,h\, \EE\left[\omega_{\nabla\psi}\left(C_{16}\,(\omega^i_X(h)\wedge1) + \sqrt{C_{11}\,h}\right) + \omega^i_X(h)\wedge 1\,\vert\,\mathcal{F}_{t_i}\right]\\
&\leq C_{17}\,h\, \EE\left[\omega_{\nabla\psi}\left(C_{18}\,((\eta_i\,h^{1/2})\wedge1) + \sqrt{C_{11}\,h})\right)
+ (\eta_i\,h^{1/2})\wedge 1\right],
\end{align*}
where we used \eqref{eq.RBSDE.exist.newEq.1} and \eqref{eq.RBSDE.exist.newEq.2} to obtain the second inequality.

\smallskip

{\bf Step 2.c} We conclude by combining the above display, \eqref{eq useful A tilde one} and \eqref{eq bound for A-term}:
\begin{align*}
&A_n  + \tilde{A}_n \\
&\ge -C_{17}\,\EE\left[\omega_{\nabla\psi}\left(C_{18}\,((\eta_1\,h^{1/2})\wedge1)+\sqrt{C_{11}\,h}\right)
+ (\eta_1\,h^{1/2})\wedge 1\right] - C_{19}\,\omega_{\nabla \psi}(\bar{C}_2 h).
\end{align*}
Setting $h \rightarrow 0$, we deduce that the right hand side of the above converges to zero (e.g., via the monotone convergence theorem). The latter yields \eqref{eq.exist.pf.step 2.eq2}, concluding the proof of the existence statement of the proposition and of the Markovian representation of $Y$. The uniqueness follows from Theorem \ref{Thm-Uniqu1}.
\eproof 

\smallskip

\begin{Remark}
The estimates \eqref{eq le stab R - psi}--\eqref{eq le stab R - psi.new} in Lemma \ref{le control ode} (more precisely, their derivation) provide a refinement of Theorem 2.2 in \cite{Lions-Sznitman-84} by showing that the Skorokhod's map is Lipschitz-continuous (as opposed to $1/2$-H\"older-continuous) with respect to the uniform norm on any set of uniformly Lipschitz paths, provided that $d=2$ and that the domain $\cD$ is simply connected and satisfies Assumptions \ref{ass:main}.
\end{Remark}

\smallskip

\begin{Theorem}
    \label{thm existence bsde exogeneous}
{Let Assumption \ref{ass:brownian-case} holds, let $\xi$ be a $\mathcal{F}_T$-measurable r.v. with values in $\bar{\cD}$, and let $f: \Omega \times [0,T] \rightarrow \mathbb{R}^d$ be progressively measurable with $|f|^{1/2} \in \widebar{\mathscr{H}^{\infty}}^{\BP{2}}$.} Then, there exists a $\Gamma$-martingale $Y$ with the drift $f$ and with the terminal value $Y_T = \xi$. Moreover, thanks to Theorem \ref{Thm-Uniqu1}, such a $\Gamma$-martingale is unique.
\end{Theorem}

{
\proof 
{\bf Step 1} In this step we prove the theorem under additional assumptions. First, we choose an arbitrary $\ell \ge 1$, consider a grid $\Re := \set{0 =: t_0 < \dots < t_\ell :=T}$, and set $b^\top=(1,0)\in\R^{1+d'\ell}$, $\sigma^\top(t) = (\mathrm{diag}_{d'}(0),\mathrm{diag}_{d'}(\bone_{\set{t\le t_1}}),\ldots,\mathrm{diag}_{d'}(\bone_{\set{t\le t_{\ell}}})) \in \cM_{d',1+d'\ell}$, so that $M=1+d'\ell$, $X^0_t=t$ and $X_t^{i+ (j-1)d'} = W^i_{t \wedge t_j}$ for all $i=1,\ldots,d'$, $j=1,\ldots,\ell$, where $X$ solves \eqref{eq SDE}, started from the initial value zero at time zero. Then, we consider bounded Lipschitz functions $g:\,\R^{d'\ell}\rightarrow\bar\cD$ and $F:\,\R^{1+d'\ell}\rightarrow\R^d$, and set $\xi:=g(X_T)=g(W_{t_1},\ldots,W_{t_\ell})$, $f_t:=F(X_t)=F(t,W_{t\wedge t_1},\ldots,W_{t\wedge t_\ell})$. Iterating Proposition \ref{prop:Gamma.mtg.drift.Markov}, we construct a $\Gamma$-martingale with drift $f$ and terminal value $\xi$. Its uniqueness follows from Theorem \ref{Thm-Uniqu1}.

\smallskip

{\bf Step 2} In this step, we consider a general terminal value $\xi$, given by an $\cF_T$-measurable random variable taking values in $\bar\cD$, while keeping the same drift, given by $f_t:=F(t,W_{t\wedge t_1},\ldots,W_{t\wedge t_\ell})$.
The main idea of this step is to approximate $\xi$ by $\xi_n = g_n(W_{s_1},\dots,W_{s_{n}})$, where $g_n$ is Lipschitz. However, a modicum of care is needed here because $g_n$ must take values in a potentially non-convex set $\bar\cD$, making it difficult to apply the standard approximation results. To address this issue, we notice that $\xi$ can be approximated, with arbitrary precision in $\cL^q$, for any $q\geq1$, by an $\cF_T$-measurable random variable $\eta$ taking values in a finite set in $\bar\cD$, denoted by $\{y_1,\ldots,y_k\}$. Next, we connect $y_{i-1}$ to $y_i$ with a Lipschitz curve $\gamma_i:\,[0,1]\rightarrow\bar\cD$, for $i=2,\ldots,k$, and define a new $\cF_T$-measurable random variable $\tilde\eta$, with values in $\{1,\ldots,k\}\subset\R$ as follows: $\tilde\eta=i$ if and only if $\eta=y_i$. Standard approximation results (cf., \cite{Nualart-06}) yield, for any $\varepsilon>0$, the existence of $n\geq 1$, $0\leq s_1<\cdots<s_n\leq T$, and a Lipschitz function $\tilde g$, such that $\tilde g(W_{s_1},\dots,W_{s_{n}})$ is within $\varepsilon$ away from $\tilde\eta$, with respect to $\cL^q$ norm. Without loss of generality, we can assume that $\tilde g$ takes values in $[1,k]$. Then, we define
\begin{align*}
& g(W_{s_1},\dots,W_{s_{n}}):=\gamma_{\lfloor\tilde g(W_{s_1},\dots,W_{s_{n}})\rfloor}(\tilde g(W_{s_1},\dots,W_{s_{n}})-\lfloor\tilde g(W_{s_1},\dots,W_{s_{n}})\rfloor).
\end{align*}
It is easy to see that $g$ is Lipschitz and that there exists a constant $C$ such that
\begin{align*}
& |\eta - g(W_{s_1},\dots,W_{s_{n}})| \leq C\,|\tilde\eta - \tilde g(W_{s_1},\dots,W_{s_{n}})|.
\end{align*}
Collecting the above, we conclude that there exist a sequence of partitions $\cS_n:=\set{0 \leq s^n_1 < \dots < s^n_n \leq T}$ of $[0,T]$ and the associated Lipschitz functions $g_n:\,\R^{d'n}\rightarrow\bar\cD$, such that $\|\xi-\xi^n\|_{\cL^2}$ converges to zero, where $\xi^n:=g_n(W_{s^n_1},\dots,W_{s^n_{n}})$.

Without loss of generality, we assume that $\Re=\set{0 =: t_0 < \dots < t_\ell :=T} \subset \cS_n$ for all $n$ and drop the superscript `$n$' in $s^n_i$. Using the results of Step 1, we obtain a sequence $(Y^n)$ of $\Gamma$-martingales with the drift $f_t = F(t,W_{t_1 \wedge t},\dots,W_{t_\ell \wedge t})$ and the terminal condition $\xi^n$. As $(\xi_n)$ is a Cauchy sequence in $\cL^2$, we apply Proposition \ref{prop stability} with $\alpha_t = 1+ \|f\|_{\mathscr{S}^{\infty}}$, to deduce that $(Y^n)$ is a Cauchy sequence in $\mathscr{S}^2$. Thus $(Y^n)$ converges to a limit denoted $Y \in \mathscr{S}^2$. 
To show that $Y$ is a $\Gamma$-martingale with drift $f$, we first notice that, for all $n\geq1$, all global special $\Gamma$-convex function $\psi$, and all $0\le t \le t' \le T$,
\begin{align*}
    \EE\left[ \psi(Y^{n}_{t'}) - \psi(Y^{n}_{t}) +\int_{t}^{t'} \nabla \psi(Y^n_s)\cdot f_s \,\ud s \,\vert\,\mathcal{F}_{t}  \right] \ge 0.
\end{align*} 
Thanks to the $\mathscr{S}^2$-convergence of $(Y^n)$ to $Y$, we can pass to the limit in the above inequality, concluding that $Y$ is a $\Gamma$-martingale with drift $f$.
Thus, we have constructed a $\Gamma$-martingale with a drift given by $F(t,W_{t\wedge t_1},\ldots,W_{t\wedge t_\ell})$, where $F$ is bounded and Lipschitz, and with a general $\cF_T$-measurable terminal value.

\smallskip

{\bf Step 3} In this step, we consider a general $\xi$ and a progressively measurable bounded $f$. First, we approximate $f$ (in $\HP{2}$) via a piecewise constant adapted process $\sum_{j=1}^{\ell} \eta_j\1_{\set{t_{j-1} \le t < t_{j}}}$, where each $\eta_j$ is $\cF_{t_{j-1}}$-measurable. Then, standard approximation results (cf., \cite{Nualart-06}) imply that every $\eta_j$ can be approximated with arbitrary precision (in $\cL^2$) by $F_j(W_{s_1^j},\ldots,W_{s_k^j})$, with $s^j_i\leq t_{j-1}$. All in all, we obtain an approximation of $f$ by 
\begin{align*}
&\sum_{j=1}^{\ell} F_j(W_{t\wedge s_1^j},\ldots,W_{t\wedge s_k^j})\1_{\set{t_{j-1} \le t < t_{j}}}
=:F(t,W_{t\wedge t_1},\ldots,W_{t\wedge t_n}).
\end{align*}
Thus, applying the results of Step 2, we obtain a sequence $(Y^n)$ of $\Gamma$-martingales with drifts given by $f^n_t:=F_n(t,W_{t^n_1 \wedge t},\dots,W_{t^n_n \wedge t})$ and with terminal condition $\xi$, such that $f^n$ converges to $f$ in $\HP{2}$. Repeating the arguments in the last paragraph of Step 2, we deduce the convergence of $(Y^n)$ to a $\Gamma$-martingale with drift $f$ and terminal value $\xi$.

\smallskip

{\bf Step 4.} In this step, we consider a general $\xi$ and a general $f$ with $|f|^{1/2} \in \widebar{\mathscr{H}^{\infty}}^{\BP{2}}$. First, we assume $|f|^{1/2} \in \mathscr{H}^{\infty}$. 
Approximating $f$ with $f^n = (-n) \vee f  \wedge n$ (component-wise), we notice that $\||f^n-f|^{1/2}\|_{\mathscr{H}^{\infty}}$ and, in turn, $\||f^n-f|^{1/2}\|_{\BP{2}}$ converge to zero. Then, we consider the $\Gamma$-martingale $Y^n$ with drift $f^n$, which exists according to Step 3. 
As $\sup_n\||f^n|^{1/2}\|_{\mathscr{H}^{\infty}}<\infty$, it is easy to see that
\begin{align}\label{eq.sec4.last.step}
\sup_n\EE e^{\lambda\int_0^T |f^n_s|\,ds} < \infty\quad \forall\,\lambda>0.
\end{align}
Thus, we apply Proposition \ref{prop stability 2} to deduce that $(Y^n)$ is a Cauchy sequence in $\mathscr{S}^2$. Using the characterization of $\Gamma$-martingales via global special $\Gamma$-convex functions (as in the last paragraph of Step 2) and the fact that $f^n$ converges to $f$ in $\HP{1}$, we conclude that $(Y^n)$ converges to a $\Gamma$-martingale with drift $f$ and terminal value $\xi$.

Finally, we consider $f$ such that $|f|^{1/2} \in \widebar{\mathscr{H}^{\infty}}^{\BP{2}}$. Then, there exists a sequence $(f^n)$ with $|f^n|^{1/2} \in \mathscr{H}^{\infty}$ for all $n \in \mathbb{N}$ and $\lim_{n \rightarrow + \infty}\| |f^n-f|^{1/2} \|_{\mathscr{B}^2}=0$. Next, we notice that $\||f^n|^{1/2}\|^2_{\mathscr{B}^2}\leq \||f^n-f|^{1/2}\|^2_{\mathscr{B}^2} + \||f|^{1/2}\|^2_{\mathscr{B}^2}$ and the latter converges to a finite number. Thus, $\sup_n \||f^n|^{1/2}\|_{\mathscr{B}^2}<\infty$ and, using the ``slicing" method for BMO martingales and John-Nirenberg inequality (see Remark \ref{rem:JN.ineq}), we deduce \eqref{eq.sec4.last.step}.
Then, we apply Proposition \ref{prop stability 2} to deduce that $(Y^n)$ is a Cauchy sequence in $\mathscr{S}^2$.
Since $f^n$ converges to $f$ in $\HP{1}$, we repeat once again the arguments in the last paragraph of Step 2, to conclude that $(Y^n)$ converges to a $\Gamma$-martingale with drift $f$ and terminal value $\xi$.
\eproof 
} 

\subsection{An Existence result for reflected BSDEs in $\bar{\cD}$}

We are now able to state our existence result for solutions to reflected BSDEs in $\bar{\cD}$.

\begin{Theorem}
    \label{thm existence bsde}
    {Let Assumption \ref{ass:brownian-case} hold, with the associated $f$ and $\xi$.} Then, there exists a solution $(Y,Z,K)$ to the reflected BSDE \eqref{eq reflected bsde} with the generator $f$ and the terminal value $\xi$. Moreover, thanks to Theorem \ref{Thm-Uniqu2}, this solution is unique.
\end{Theorem}

To prove the above theorem, we the following auxiliary result.

\begin{Proposition}
    \label{prop existence bsde intermediary step}
   {Let Assumption \ref{ass:brownian-case} hold, and assume that there exists a process $V \in \mathscr{B}^2$ and a deterministic function $f_v$ such that $f(t,y,z) := f_v(t,y,V_t)$}. Then, there exists a solution $(Y,Z,K)$ to the reflected BSDE \eqref{eq reflected bsde} with the generator $f$ and the terminal value $\xi$. Moreover, thanks to Theorem \ref{Thm-Uniqu2}, this solution is unique.
\end{Proposition}

\proof
    We argue by contradiction. First, we consider $U^1$ and ${U^2}$, two arbitrary continuous adapted processes with values in $\bar{\cD}$. By applying Theorem \ref{thm existence bsde exogeneous}, there exists a unique solution $(Y^1,Z^1,K^1)$ (resp. $({Y^2},{Z^2},{K^2})$) solution to the reflected BSDE \eqref{eq reflected bsde} with exogenous generator $f(.,U^1_.,V_.)$ (resp. $f(.,U^2_.,V_.)$) and same terminal condition $\xi$. We can apply Proposition \ref{prop stability} with $\alpha_t = 1+|f(t,0,V_t)|+C_{f,y}\sup_{y \in \bar{\cD}} |y|$, since we have $C_{f,z}=0$, to get, for $\lambda$ large enough,
    \begin{align*}
        \esp{\sup_{t\in[0,T]}\Gamma^{\lambda}_t |Y^1_t-Y^2_t|^2} & \le\esp{\sup_{t\in[0,T]}\Gamma^{\lambda}_t \Psi(Y^1_t,Y^2_t)} \le \eta C_{f,y} T\esp{\sup_{t\in[0,T]}\Gamma^{\lambda}_t \Psi(U^1_t,U^2_t)}\\
        &\le \eta C_{f,y} T C \esp{\sup_{t\in[0,T]}\Gamma^{\lambda}_t |U^1_t-U^2_t|^2}
    \end{align*}
    where we have used \eqref{equivalence normes} in the last inequality.
    Then, we set $\eta <(C_{f,y} T C)^{-1}$ to get a contraction in a Banach space:  there exists a unique fixed point $Y$ such that there exists a unique solution $(Y,Z,K)$ to the reflected BSDE \eqref{eq reflected bsde} with generator $f(.,Y_.,V_.)$ which gives us the result.
\eproof

\smallskip

We are now able to prove Theorem \ref{thm existence bsde}.

\proof
{\bf Step 1}. Let us start by introducing constants $_{\eqref{cste BMO norm Picard}}$ and $C'_{\eqref{cste BMO norm Picard}}$ given by
\begin{equation}
\label{cste BMO norm Picard}
C_{\eqref{cste BMO norm Picard}} :=2C_{\eqref{eq Prop-aprioriZ}}\left(1+\frac{C_{\eqref{eq Prop-aprioriZ}}C_{f,z}^2T}{2}\right), \quad C'_{\eqref{cste BMO norm Picard}} := C_{\eqref{cste BMO norm Picard}} +C_{\eqref{cste BMO norm Picard}}  \| |f(.,0,0)|^{1/2} \|
\end{equation}
and let us define the following closed subset of a Banach space
\begin{equation*}
    \mathscr{B}^2_{b} := \left\{ Z \in \BP{2} | \|Z\|_{\mathscr{B}^2}^2 \leqslant C'_{\eqref{cste BMO norm Picard}}   \right\}
\end{equation*}
equipped with the equivalent norm, for any weight parameter $\alpha \ge 0$,
\begin{equation*}
    \|Z\|_{\alpha,\mathscr{B}_d^2} := \NL{\infty}{\mathrm{sup}_{t \in [0,T]} \esp{\int_{t}^T e^{\alpha s} |Z_s|^2 \ud s | \cF_{t}}}^\frac12 . \;
\end{equation*}

Let us consider $V \in  \mathscr{B}^2_{b}$. By applying Proposition \ref{prop existence bsde intermediary step}, there exists a unique solution $(Y,Z,K)$ to the reflected BSDE with generator $f_V(t,y,z) := f(t,y,V_t)$. We apply Proposition \ref{Prop-aprioriZ} and Young inequality to get the following estimate on $Z$:
\begin{align*}
    \|Z\|_{\BP{2}}^2 \leqslant& C_{\eqref{eq Prop-aprioriZ}}(1+\||f_V(.,0,0)|^{1/2}\|_{\BP{2}}^2)\\
    \le & C_{\eqref{eq Prop-aprioriZ}}(1+\||f(.,0,0)|^{1/2}\|_{\BP{2}}^2+C_{f,z}\||V|^{1/2}\|_{\BP{2}}^2)\\
    \le & C_{\eqref{eq Prop-aprioriZ}}\left(1+\||f(.,0,0)|^{1/2}\|_{\BP{2}}^2+\frac{C_{\eqref{eq Prop-aprioriZ}}C_{f,z}^2 T}{2}\right)+\frac12 \|V\|_{\BP{2}}^2
     \leqslant  C'_{\eqref{cste BMO norm Picard}} 
\end{align*}
which means that $Z \in \mathscr{B}^2_{b}$. In other words, $\Phi : V \mapsto  Z$ is a function from $\mathscr{B}^2_{b}$ to itself.

\smallskip

{\bf Step 2}. We now prove that $\Phi$ is a contraction which is sufficient to conclude. Let us consider $V^1, V^2 \in \mathscr{B}^2_{b}$ and let us denote $(Y^1,Z^1,K^1)$ and $(Y^2,Z^2,K^2)$ the solutions of the associated reflected BSDEs.
For $\alpha \ge 0$, $\beta\ge 0$ and $\eta > 0$, we denote 
\[
\Gamma_{t_1,t_2} := e^{\int_{t_1}^{t_2} \alpha+\beta(|V^2|+|f(s,0,0)|) + \eta |Z^2_s|^2 \ud s}
.\]
We set $\eta$ such that $4 \eta C'_{\eqref{cste BMO norm Picard}} \leqslant 1$.  Thanks to \eqref{moment-expo-f(0)} and John-Nirenberg inequality, we can show that $\Gamma_{0,T}$ is square integrable for all $\alpha\ge 0$, $\beta \ge 0$. Parameters $\alpha$ and $\beta$ will be set after. For $\varepsilon \in (0,1/2)$, we have
\begin{align}
    &\Psi(Y^1_t,Y^2_t) + \varepsilon \mathbb{E}_t \left[\int_t^T e^{\alpha s}|Z^1_s-Z^2_s|^2 \ud s\right] \nonumber\\
    \leqslant& \Psi(Y^1_t,Y^2_t) + \mathbb{E}_t \left[\int_t^T \Gamma_{t,s}|Z^1_s-R(\theta(Y^1_s,{Y^2_s}))Z^2_s|^2 \ud s\right] + 2\varepsilon \mathbb{E}_t \left[\int_t^T \Gamma_{t,s} \|I -R(\theta(Y^1_s,{Y^2_s}))\|^2|Z^2_s|^2 \ud s\right] \nonumber\\
    \leqslant& \Psi(Y^1_t,Y^2_t) + \mathbb{E}_t \left[\int_t^T \Gamma_{t,s} |Z^1_s-R(\theta(Y^1_s,{Y^2_s}))Z^2_s|^2 \ud s\right] + \varepsilon C \mathbb{E}_t \left[\int_t^T \Gamma_{t,s} \Psi(Y^1_s,Y^2_s)|Z^2_s|^2 \ud s\right]. \label{proof existence eq1}
\end{align}

By applying Proposition \ref{Ito-formula-Psi-bis} and taking the expectation after checking that the stochastic integral term is a martingale, we have
\begin{align}
    \label{proof existence eq2}
    &\Psi(Y^1_t,{Y}^2_t) +  \mathbb{E}_t \left[\int_t^T \Gamma_{t,s} |Z^1_s-R(\theta(Y^1_s,{Y^2_s})){Z}^2_s|^2 \ud s\right] \\ \nonumber
    \leqslant & \mathbb{E}_t \left[\int_t^T \Gamma_{t,s}\Psi^{1/2}(Y^1_s,{Y^2_s}) |V^1_s-V^2_s| \ud s \right]
    \\ \nonumber
    &+ C\mathbb{E}_t \left[\int_t^T \Gamma_{t,s}\Psi^{1/2}(Y^1_s,{Y^2_s}) \|I-R(\theta(Y^1_s,{Y}^2_s))\|(1+|f(s,0,0)|+|{V}^2_s|) \ud s \right]\\ \nonumber
    &- \mathbb{E}_t \left[\int_t^T \Gamma_{t,s}\Psi(Y^1_s,{Y^2_s})(\alpha+\beta|{V}^2_s|+\beta|f(s,0,0)|)\ud s\right] 
    - \eta \mathbb{E}_t \left[\int_t^T \Gamma_{t,s}\Psi(Y^1_s,{Y^2_s})|{Z}^2_s|^2\ud s\right]\\ \nonumber
    \leqslant & C_{\gamma}\mathbb{E}_t \left[\int_t^T \Gamma_{t,s} \Psi(Y^1_s,{Y^2_s})\ud s \right]\\ \nonumber
    &+ \gamma  \mathbb{E}_t \left[ e^{\int_{t}^{T} \beta(|{V}^2_s|+|f(s,0,0)|) + \eta |{Z}^2_s|^2 \ud s} \right]^{1/2}\mathbb{E}_t \left[ \left(\int_t^T e^{\alpha (s-t)}|V^1_s-{V}^2_s|^2 \ud s\right)^2 \right]^{1/2}
    \\ \nonumber
    &+ C\mathbb{E}_t \left[\int_t^T \Gamma_{t,s}\Psi(Y^1_s,{Y^2_s}) (1+|f(s,0,0)|+|{V}^2_s|) \ud s \right]\\ \nonumber
    &-  \mathbb{E}_t \left[\int_t^T \Gamma_{t,s}\Psi(Y^1_s,{Y^2_s})(\alpha+\beta|{V}^2_s|+\beta|f(s,0,0)|)\ud s\right] - \eta \mathbb{E}_t \left[\int_t^T \Gamma_{t,s}\Psi(Y^1_s,{Y^2_s})|{Z}^2_s|^2\ud s\right],
\end{align}
where the parameter $\gamma>0$ comes from the application of Young estimate and will be set after. We use Energy inequality for BMO martingales to get
\begin{align}
    \label{proof existence eq3}
    \mathbb{E}_t \left[ \left(\int_t^T e^{\alpha s}|V^1_s-{V}^2_s|^2 \ud s\right)^2 \right]^{1/2} &\leqslant 
    \sup_{t \in [0,T]} \mathbb{E}_t \left[ \left(\int_t^T e^{\alpha s}|V^1_s-{V}^2_s|^2 \ud s\right)^2 \right]^{1/2} \leqslant C \| V^1-{V^2}\|_{\alpha,\mathscr{B}^2_b}^2 .   
\end{align}
Finally, we plug \eqref{proof existence eq1}, \eqref{proof existence eq2} and \eqref{proof existence eq3} together to get
\begin{align}
    \label{proof existence eq4}
    &\Psi(Y^1_t,{Y}^2_t) + \varepsilon \mathbb{E}_t \left[\int_t^T e^{\alpha (s-t)}|Z^1_s-{Z}^2_s|^2 \ud s\right] \\ \nonumber
    \leqslant& \gamma  \mathbb{E}_t \left[ e^{\int_{t}^{T} \beta(|{V}^2_s|+|f(s,0,0)|) + \eta |{Z}^2_s|^2 \ud s} \right]^{1/2}e^{-\alpha t}C \| V^1-{V^2}\|_{\alpha,\mathscr{B}^2_b}^2 
    \\ \nonumber
    & + (C-\beta)\mathbb{E}_t \left[\int_t^T \Gamma_{t,s}\Psi(Y^1_s,{Y^2_s}) (|f(s,0,0)|+|{V}^2_s|) \ud s \right] +(C_{\gamma}-\alpha) \mathbb{E}_t \left[\int_t^T \Gamma_{t,s}\Psi(Y^1_s,{Y^2_s}) \ud s \right]\\
    & + (\varepsilon C - \eta) \mathbb{E}_t \left[\int_t^T \Gamma_{t,s} \Psi(Y^1_s,{Y}^2_s)|{Z}^2_s|^2 \ud s\right] \nonumber
\end{align}
Now, we set $\beta$ and $\varepsilon$ such that $C-\beta<0$ and $\varepsilon C-\eta<0$. 
Let us remark that John-Nirenberg inequality gives us, once again, that
\begin{equation}
\label{norm gamma term}
    \sup_{t \in [0,T]}  \mathbb{E}_t \left[ e^{\int_{t}^{T} \beta(|{V}^2_s|+|f(s,0,0)|) + \eta |{Z}^2_s|^2 \ud s} \right]^{1/2} \leqslant C_{\eqref{norm gamma term}}
\end{equation}
where $C_{\eqref{norm gamma term}}$ depends on $V$ and $Z$ only through $C'_{\eqref{cste BMO norm Picard}}$. Indeed we have
\begin{align*}
    \sup_{t \in [0,T]}  \mathbb{E}_t \left[ e^{\int_{t}^{T} \beta(|{V}^2_s|+|f(s,0,0)|) + \eta |{Z}^2_s|^2 \ud s} \right] \le & \sup_{t \in [0,T]}  \mathbb{E}_t \left[ e^{\int_{t}^{T} 3\beta|{V}^2_s| \ud s} +e^{\int_{t}^{T} 3|f(s,0,0)| \ud s}+e^{\int_{t}^{T} 3\eta |{Z}^2_s|^2 \ud s}\right]
\end{align*}
Moreover, $3\beta|{V}^2_s|\le \frac{(3\beta)^2C'_{\eqref{cste BMO norm Picard}}}{2}+ \frac{\beta|{V}^2_s|^2}{2C'_{\eqref{cste BMO norm Picard}}}$, so John-Nirenberg inequality gives us that we can take 
\begin{align*}
    C_{\eqref{norm gamma term}}^2
     =  2e^{(3\beta)^2C'_{\eqref{cste BMO norm Picard}}} + \sup_{t \in [0,T]}  \mathbb{E}_t \left[ e^{\int_{t}^{T} 3|f(s,0,0)| \ud s}\right] + \frac{1}{1-3\eta C'_{\eqref{cste BMO norm Picard}}} .
\end{align*}
Then we can set $\gamma$ such that
$$\gamma C_{\eqref{norm gamma term}} C\leqslant \frac{\varepsilon}{2}. $$

Finally, we set $\alpha$ in order to have $C_{\gamma}-\alpha<0$. Thus, \eqref{proof existence eq4} becomes 
\begin{align*}
    &e^{\alpha t}\Psi(Y^1_t,{Y}^2_t) + \varepsilon \mathbb{E}_t \left[\int_t^T e^{\alpha s}|Z^1_s-{Z}^2_s|^2 \ud s\right]
    \leqslant \frac{\varepsilon}{2} \| V^1-{V^2}\|_{\alpha,\mathscr{B}^2_b}^2, \quad \forall t \in [0,T] 
\end{align*}
which leads us to
\begin{equation*}
    \| Z^1-Z^2\|_{\alpha,\mathscr{B}^2_b}^2 \leqslant \frac{1}{2} \| V^1-V^2\|_{\alpha,\mathscr{B}^2_b}^2.
\end{equation*}
So, we have proved that $\Phi$ is a contraction and the Banach's fixed-point theorem allows to conclude the existence part of Theorem \ref{thm existence bsde}. 

\smallskip

{\bf Step 3}. As mentioned, the uniqueness result is a direct application of Theorem \ref{Thm-Uniqu2} thanks to \eqref{moment-expo-f(t,y,z)}.
\eproof

\section{Appendix}
\label{subse about frechet mean}

 
In this appendix, we study the properties of a Fr\'echet mean of a random variable taking values in a bounded domain $ \mathfrak{D}$, which satisfies the regularity property \textbf{(R)} and which is a geodesic space with the geodesic distance function denoted by $d_{\mathfrak{D}}$. Except for Proposition \ref{Prop uniqueness Frechet}, we do not restrict the analysis to $d=2$, {but we do assume the following.

\begin{Assumption}\label{ass:4}
We assume that, for any $(x,y) \in \bar{\mathfrak{D}} \times \bar{\mathfrak{D}}$, the minimizing geodesic between $x$ and $y$ is unique, that, for any $y \in \bar{\mathfrak{D}}$, the function $x \mapsto \Psi(x,y):=d^2_{\mathfrak{D}}(x,y)$ is $C^1$, and that $\nabla\Psi(x,y)=-2\overrightarrow{xy}$ {(with $\overrightarrow{xy}$ defined in Definition \ref{def:re de rot-angle})}.
\end{Assumption}
}



\begin{Definition}
 \label{def Frechet mean}
 Let Assumption \ref{ass:4} hold, and let $\xi$ be a random variable with values in $\bar{\mathfrak{D}}$. A point $q \in \bar{\mathfrak{D}}$ is a Fr\'echet mean of $\xi$ (w.r.t. $\bar{\mathfrak{D}}$) if it is a minimizer of the function 
$$Q(x) = \mathbb{E}\left[ \Psi(x,\xi)\right], \quad x \in \bar{\mathfrak{D}}.$$
\end{Definition}

\smallskip

\begin{Proposition}
 \label{prop Frechet mean}
  Let Assumption \ref{ass:4} hold, and let $\xi$ be a random variable with values in $\bar{\mathfrak{D}}$. Then, there exists at least one Fr\'echet mean $q$ of $\xi$. Moreover, any Fr\'echet mean $q$ of $\xi$ satisfies
  $$\int_{\bar{\mathfrak{D}}} \overrightarrow{qy} \,\mu(dy) =0,$$
  where $\mu$ is the distribution of $\xi$.
\end{Proposition}

\proof
The existence follows from the compactness of $\bar{\mathfrak{D}}$ and the continuity of $\Psi$. To show the second statement, let us consider an arbitrary Fréchet mean $q$ of $\xi$. Then, it is a minimizer of the function 
$$Q(x)= \int_{\mathfrak{D}} \Psi(x,y)\mu(dy).$$
Since $x \mapsto \Psi(x,y)$ is $C^1$ on $\bar{\mathfrak{D}}$, we easily get that $Q$ is differentiable and 
$$\nabla Q(x) = -2\int_{\bar{\mathfrak{D}}} \overrightarrow{xy} \mu(dy).$$
In particular, for any $u \in \mathbb{R}^2$ such that $q+\varepsilon u \in \bar{\mathfrak{D}}$ for all small enough $\varepsilon>0$, we must have
$$\nabla Q(q) \cdot u  = \lim_{\varepsilon \rightarrow + \infty} \frac{Q(x+\varepsilon u)-Q(x)}{\varepsilon} \geqslant 0.$$
If $q \in \mathfrak{D}$, the above holds with any $u \in\mathbb{R}^2$, which implies $\nabla Q(q)=0$. Otherwise, $q \in \partial \mathfrak{D}$, and we can take any $u$ in the dual cone of $-\mathfrak{n}(q)$ (which is not empty according to Proposition \ref{prop:extsphere:interiorcone}), which yields 
$\nabla Q(q)\in - \mathfrak{n}(q)$. To conclude, we just have to remark that 
$$0 \leqslant \nabla Q(q) \cdot \nabla Q(q) = 2\int_{\bar{\mathfrak{D}}} \overrightarrow{qy} \cdot (-\nabla Q(q)) \mu(dy)\leqslant 0$$
since $-\nabla Q(q) \in \mathfrak{n}(q)$.
\eproof

\smallskip

\begin{Proposition}
 \label{prop Jensen for Frechet}
 Let Assumption \ref{ass:4} hold, and consider a $\Gamma$-convex function $\psi : \bar{\mathfrak{D}} \rightarrow \mathbb{R}$, a random variable $\xi$ with values in $\bar{\mathfrak{D}}$, and a Fr\'echet mean $q$ of $\xi$. Then, the following version of Jensen's inequality holds:
 \begin{equation}
    \label{ineq Jensen}
 \psi(q) \leqslant \mathbb{E}[\psi(\xi)].
 \end{equation}
\end{Proposition}



\proof
By the $\Gamma$-convexity assumption on $\psi$, we have that, for any minimizing geodesic $\gamma$ in $\bar{\mathfrak{D}}$,
$$\psi(\gamma_0)+u(\gamma_0) \cdot  \dot \gamma_0
\leqslant \psi(\gamma_1),$$
for all $u(\gamma_0) \in \partial \psi(\gamma_0)$.
Then, by considering the minimizing geodesics between $q$ and the points $y \in \bar{\mathfrak{D}}$, and integrating in $y$ with respect to the distribution $\mu$ of $\xi$, we get
$$\psi(q) + u(q) \cdot \int_{\bar{\mathfrak{D}}} \overrightarrow{qy} \mu(dy) \leqslant \int_{\bar{\mathfrak{D}}} \psi(y) \mu(dy).$$  
Applying Proposition \ref{prop Frechet mean}, we obtain the desired result.
\eproof

\smallskip

As a consequence of the previous proposition we deduce the uniqueness of a Fréchet mean.

\begin{Proposition}
 \label{Prop uniqueness Frechet}
 {Let Assumption \ref{ass:4} hold.
 }
 Then, for any random variables $\xi$ and $\xi'$ with values in $\bar{\mathfrak{D}}$, we have 
 $$\Psi(q_1,q_2) \leqslant \mathbb{E}[\Psi(\xi,\xi')],$$ 
 where $q_1$ (resp., $q_2$) is a Fr\'echet mean of $\xi$ (resp., $\xi'$). In particular, 
 there exists a unique Fr\'echet mean of $\xi$. We denote it $\mathcal{E}[\xi]$.
\end{Proposition}

\proof
We consider two Fréchet means $q_1$ and $q_2$ of $\xi$. Then $(q_1,q_2)$ is a Fréchet mean of $(\xi,\xi)$ on the Euclidean product manifold $\bar{\mathfrak{D}} \times \bar{\mathfrak{D}}$ with a boundary, where the minimizing geodesics are given by the pairs of minimizing geodesics in $\bar{\mathfrak{D}}$. Moreover, we observe that $\bar{\mathfrak{D}} \times \bar{\mathfrak{D}}$ satisfies Assumption \ref{ass:4}. Then, applying Proposition \ref{prop Jensen for Frechet} to the function $\Psi$ on $\bar{\mathfrak{D}} \times \bar{\mathfrak{D}}$, we obtain
$$\Psi(q_1,q_2) \leqslant \mathbb{E}[\Psi(\xi,\xi)]=0,$$
which gives us $q_1=q_2$.
\eproof


\def\cprime{$'$}


\end{document}